%

\magnification=\magstep1
\input amstex
\UseAMSsymbols
\input pictex
\vsize=23truecm
\NoBlackBoxes
\parindent=20pt
  \font\gross=cmbx10 scaled\magstep1 
         
      \font\tt=cmtt10


      \def\supp {\operatorname{supp}}

      \def\arr#1#2{\arrow <1.5mm> [0.25,0.75] from #1 to #2}
      
\def\Rahmen#1%
   {$$\vbox{\hrule\hbox%
                  {\vrule%
                       \hskip0.5cm%
                            \vbox{\vskip0.3cm\relax%
                               \hbox{$\displaystyle{#1}$}%
                                  \vskip0.3cm}%
                       \hskip0.5cm%
                  \vrule}%
           \hrule}$$}


\vglue1truecm

\centerline{\gross The root posets and their rich antichains}
		       		       	          \bigskip
\centerline{Dedicated to Professor Liu Shao-Xue on his 90th Birthday}
	\bigskip 
\centerline{Claus Michael Ringel}
		  \bigskip\bigskip
{\narrower\narrower 
Abstract. Let  $\Delta$ be a (connected) Dynkin diagram of rank $n\ge 2$ and $\Phi_+ = \Phi_+(\Delta)$ 
the corresponding root poset (it consists of all positive roots with respect to 
a fixed root basis). The width of $\Phi_+$ is $n$. We will show that $\Phi_+$ is
``conical'': it is the disjoint union of $n$ solid chains. 

The rich antichains in $\Phi_+$
are the antichains of cardinality $n-1$.
It is well known
that the number of rich antichains is equal to the cardinality of $\Phi_+$.
The set $\Cal R(\Delta)$ of rich antichains in $\Phi_+$ can itself be
considered as a poset which is quite similar, but not always isomorphic,
to $\Phi_+$.

We will show that there always exists a 
unique rich antichain $A$ such that any rich antichain is contained
in the ideal generated by $A$. 
For $\Delta\neq \Bbb E_6$ all roots in $A$ have the same length, 
namely $e_2$, where $e_1 \le e_2 \le \dots \le e_n$ are the exponents of $\Delta.$ 
For  $\Delta = \Bbb E_6$, the antichain $A$ consists of 
four roots of length $e_2 = 4$ and one root of length $5$. \par }
	\bigskip

{\bf 1. Introduction.}
	\medskip 
Let $\Phi$ be a root system of (connected) Dynkin type 
$\Delta = \Bbb A_n, \Bbb B_n,\dots, \Bbb G_2$ of rank $n$.
The {\it root poset} $\Phi_+ = \Phi_+(\Delta)$ is the set of positive
roots in $\Phi$ with respect to some fixed choice of a root
basis; here, one takes the following partial ordering: 
$x \le y$ provided $y-x$ is 
a non-negative linear combination of elements of the root basis.
The root posets and its antichains play an important 
role in many parts of mathematics. 
The width of the root poset $\Phi_+(\Delta)$ is the rank of $\Delta$ (this is
the number of vertices of $\Delta$). 
If $P$ is a poset of width $n$, the antichains of $P$ 
of cardinality $n-1$ will be called the {\it rich} antichains.
It is well known
that the number of rich antichains is equal to number of positive roots 
and the paper is devoted to a study of the rich antichains of $\Phi_+$. 
	\medskip
Following Stanley [S], a finite poset $P$ is said to be {\it graded} provided all maximal chains have cardinality $m$ (the height of $P$).
Given a graded poset $P$, we denote by $P_t$ the set
of elements of height $t$; note that $P_t$ is an antichain. The cardinality of $P_t$
will be denoted by $h_t(P)$. Of course, $h_1(P)$ is just the number of minimal elements
of $P$. 

Given a poset $P$, we call a subposet $P'$ a {\it solid} subposet provided neighbors in $P'$ are neighbors in $P$ (if $x<y$ are neighbors in the subposet $P'$, this interval cannot be refined in $P$). We say that $P$ is
{\it conical} provided $P$ is a graded poset with a unique maximal element and 
is the disjoint union of solid subchains 
such that each of the subchains contains a minimal element (such a set of subchains will be called a {\it conical decomposition}). If $\Cal C$ is a 
conical decomposition $\Cal C$, the number of chains in $\Cal C$ is equal to $h_1(P).$
The cardinalities of the chains in $\Cal C$ will be called the {\it exponents of $P$}.

If $P$ is conical, we have $h_t(P) \ge h_{t+1}(P)$ for all $t\ge 1,$ thus
$(h_1(P),h_2(P),\dots)$ is a (Young) partition and the dual partition 
$(e_1(P),e_2(P),\dots)$ is the sequence of exponents. 

Finally, given a subset $X$ of a poset $P$, we denote by $\Lambda(X)  = \{w\in P
\mid w\le x \ \text{for some}\ x\in X\}$ the ideal,
by $V(X)    = \{y\in P
\mid x\le y \ \text{for some}\ x\in X\}$ 
the coideal generated by $X$. 
	\medskip
Let us consider now the root posets. 
Since any root poset has a maximal element, the definition of the partial ordering of
$\Phi_+$ shows that all maximal chains in $\Phi_+$ have the same cardinality,
thus $\Phi_+$ is a graded poset with a unique maximal element. 
For a positive root $x$, the height of $x$ in $\Phi_+$ is just 
the length of $x$, it is the sum of the coefficients 
when $x$ is written as a linear combination with respect to the root basis (see for
example [B]). 
We write $\Phi_t$ instead of $(\Phi_+)_t$.
	\medskip 
{\bf Theorem 1}. {\it A root poset is conical.}
	\medskip 
Theorem 1 strengthens the well-known assertion that $h_t(\Phi_+) \ge h_{t+1}(\Phi_+)$
for all $t\ge 1$. Actually, the sequence of exponents of $\Phi_+(\Delta)$ as defined above 
is just the usual sequence of exponents as considered in the invariant theory
of Weyl groups (this is the celebrated Shapiro-Kostant
theorem, see for example [R2] Theorem 1.4.1.1). We also should mention that
$$
 h_t(\Phi_+) + h_{g-t+1}(\Phi_+) = n
$$
for all $t$, where $g$ is the Coxeter number for $\Delta$
(see Humphreys [H, Theorem 3.20] and Armstrong [Ar, Theorem 5.4.1]).
	\medskip
Note that $e_1 = 1$ and, for $n\ge 2$, we have $e_2 \ge 2$. It follows that all $\Phi_t$
with $2\le t \le e_2$ are rich antichains, whereas any $\Phi_t$ with $t > e_2$ is
an antichain of cardinality at most $n-2$.
	\medskip 
The set of rich antichains of a poset can be considered as a poset $\Cal R(\Delta)$
with $A \le A'$ provided $V(A) \supseteq V(A')$ (recall that $V(A)$ 
is the coideal generated by $A$, see section 4). 
As we will see, for $\Delta \neq \Bbb E_6$, the antichain $\Phi_{e_2}$
is the largest element in $\Cal R(\Delta)$, whereas for $\Delta = \Bbb E_6,$ there is
a single rich antichain which is larger than $\Phi_{e_2}$.
	\medskip 
For any natural number $t$, let $\Phi_t'$ be the set consisting 
of join-irreducible elements of $\Phi_{t+1}$ and the elements $a\in \Phi_t$
which have no join-irreducible neighbor $a' > a$ in $\Phi_{t+1}.$ Since $\Phi_+$
is conical, the cardinality of $\Phi_t'$ is equal to the cardinality of $\Phi_t.$ 
In particular, for $t = e_2,$ $\Phi_t'$ is again a rich antichain.
	\medskip
{\bf Theorem 2.} {\it Let $\Delta$ be a Dynkin diagram
of rank at least $2$. Then $\Phi_+$ has a unique maximal rich antichain,
namely $\Phi_{e_2}'$. If $\Delta \neq \Bbb E_6,$ then 
$\Phi_{e_2}' = \Phi_{e_2}$, whereas for $\Delta = \Bbb E_6,$ the antichains
$\Phi_{e_2}',\Phi_{e_2}$ differ by one element.}
	\medskip
Here is the value of $e_2= e_2(\Delta)$:
$$
{\beginpicture
\setcoordinatesystem units <1cm,.7cm>
\plot -.6 0.5  9.3 0.5 /
\plot .4 -0.3  .4 1.3 /
\put{$\Delta$} at -.2 1 
\put{$\Bbb A_n$} at  1 1
\put{$\Bbb B_n$} at 2 1
\put{$\Bbb C_n$} at 3 1
\put{$\Bbb D_n$} at 4 1
\put{$\Bbb E_6$} at 5 1
\put{$\Bbb E_7$} at 6 1
\put{$\Bbb E_8$} at 7 1
\put{$\Bbb F_4$} at 8 1
\put{$\Bbb G_2$} at 9 1
\put{$e_2(\Delta)$} at -.2 0
\multiput{$2$} at  1 0  /
\multiput{$3$} at 2 0 3 0  4 0 /
\multiput{$4$} at 5 0 /
\multiput{$5$} at  6 0 8 0   9 0 /
\multiput{$7$} at   7 0  /
\endpicture}
$$

Let us reformulate Theorem 2. In oder to find the rich antichains of $\Phi_+(\Delta)$,
it is sufficient to look at the ideal $\Lambda(\Phi_{e_2}')$, this is a poset with
$n$ minimal and $n\!-\!1$ maximal elements. Here are 
the posets $\Lambda(\Phi_{e_3}')$ for $\Delta = \Bbb A_5, \Bbb B_5, \Bbb D_6,
\Bbb E_6, \Bbb E_7, \Bbb E_8, \Bbb F_4, \Bbb G_2.$

$$
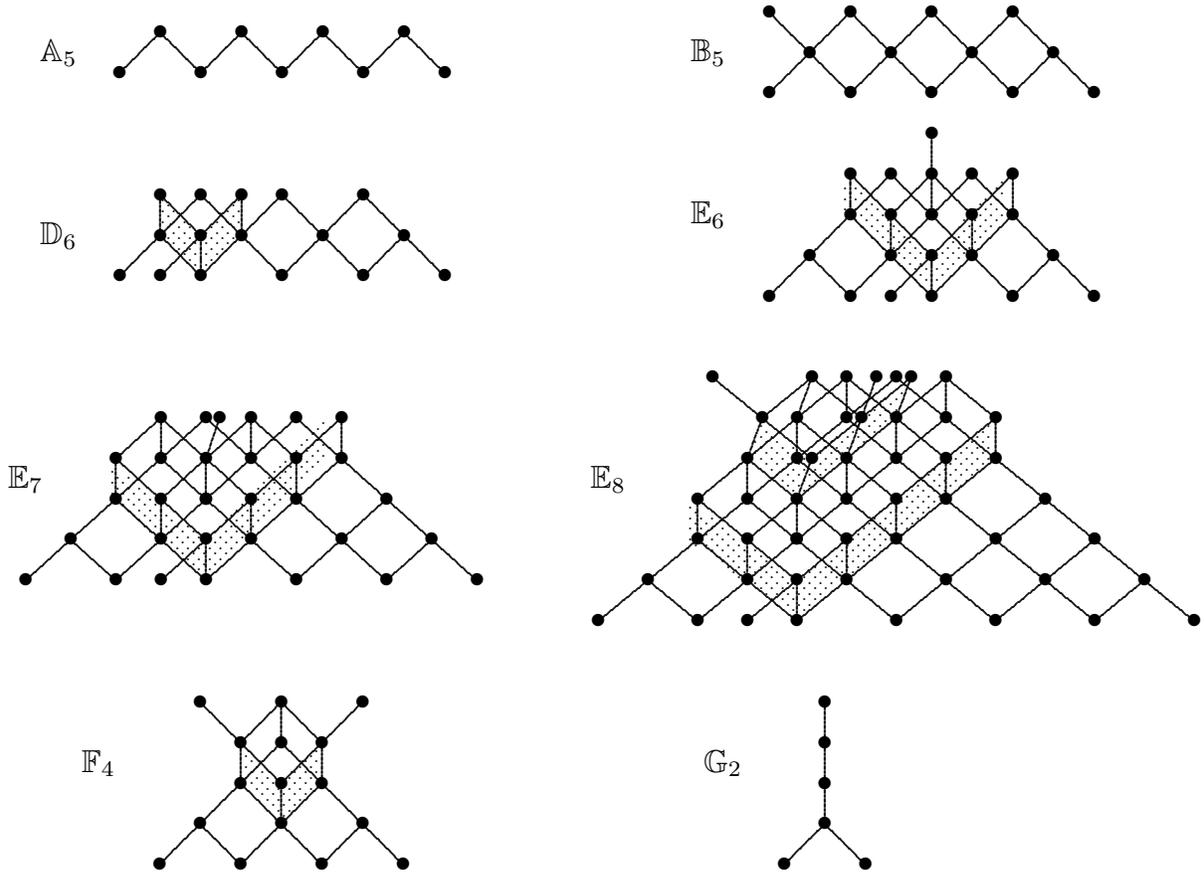
{\beginpicture
\setcoordinatesystem units <.45cm,.45cm>
\put{\beginpicture
\put{$\Bbb A_5$}  at -1.5 0.5
\multiput{$\bullet$} at 0 0  1 1  2 0  3 1   4 0  5 1  6 0  7 1  8 0  /
\plot 0 0  1 1  2 0  3 1   4 0  5 1  6 0  7 1  8 0  /
\endpicture} at 0 0
\put{\beginpicture
\put{$\Bbb B_5$}  at -1.5 1
\multiput{$\bullet$} at 0 0  1 1  2 0  3 1   4 0  5 1  6 0  7 1  8 0  
    2 2  4 2  6 2  0 2   /
\plot 0 0  1 1  2 0  3 1   4 0  5 1  6 0  7 1  8 0  /
\plot 0 2  1 1  2 2  3 1  4 2  5 1  6 2  7 1   /
\endpicture} at 16 0
\put{\beginpicture
\put{$\Bbb D_6$}  at -9.5 1
\multiput{$\bullet$} at 0 0  -1 1  -2 0  -3 1   -4 0  -5 1  -6 0  -7 1  -8 0  -7 0 
    -2 2  -4 2  -6 2   /
\plot 0 0  -1 1  -2 0  -3 1   -4 0  -5 1  -6 0  -7 1  -8 0  /
\plot -1 1  -2 2  -3 1  -4 2  -5 1  -6 2  -7 1   /
\multiput{$\bullet$} at -6 1   -5 2  -7 2  /
\plot -5 1  -5 2  -6 1  -7 2  -7 1 /
\plot -7 0  -6 1  -6 0 /
\setshadegrid span <.6mm>

\vshade   -7 1 2  <,z,,> -6 0 1  <z,,,> -5 1 2 /

\endpicture} at 0 -4.5
\put{\beginpicture
\put{$\Bbb E_6$}  at -1.5 2
\multiput{$\bullet$} at 0 0  1 1  2 0  3 1   4 0  5 1  6 0  7 1  8 0 
    2 2  4 2  6 2  5 2   /
\plot 0 0  1 1  2 0  3 1   4 0  5 1  6 0  7 1  8 0  /
\plot 1 1  2 2  3 1  4 2  5 1  6 2  7 1   /
\multiput{$\bullet$} at  3 2  4 1  4 2  3 0  /
\plot 3 1  3 2  4 1  5 2  5 1 /
\plot 3 0  4 1  4 0 /
\multiput{$\bullet$} at 2 3  3 3  4 3  5 3  6 3   /
\plot 4 2  3 3  2 2  2 3  3 2  4 3  4 2  5 3  6 2  6 3  5 2  4 3  /
\put{$\bullet$} at 4 4
\plot 4 3  4 4 /

\setshadegrid span <.6mm>
\vshade 1.8 2 3  <,z,,> 4 0 1  <z,,,> 5.8 2 3 /

\endpicture} at 16 -4
\put{\beginpicture
\setcoordinatesystem units <.5cm,.45cm>
\put{$\Bbb E_7$}  at -.0 2.5
\multiput{$\bullet$} at 0 0  1 1  2 0  3 1   4 0  5 1  6 0  7 1  8 0 
    2 2  4 2  6 2  5 2   /
\plot 0 0  1 1  2 0  3 1   4 0  5 1  6 0  7 1  8 0  /
\plot 1 1  2 2  3 1  4 2  5 1  6 2  7 1   /
\multiput{$\bullet$} at  3 2  4 1  4 2  3 0  /
\plot 3 1  3 2  4 1  5 2  5 1 /
\plot 3 0  4 1  4 0 /
\multiput{$\bullet$} at 2 3  3 3  4 3  5 3  6 3   /
\plot 4 2  3 3  2 2  2 3  3 2  4 3  4 2  5 3  6 2  6 3  5 2  4 3  /

\put{$\bullet$} at 4.3 4
\plot 4 3  4.3 4 /

\multiput{$\bullet$} at 10 0  9 1  8 2  7 3 /
\plot 6 2  7 3  10 0 /
\plot 8 0  9 1 /
\plot 7 1  8 2 /

\multiput{$\bullet$} at 3 4  4 4  5 4  6 4  7 4  /
\plot 2 3  3 4  4 3  5 4  6 3  7 4  7 3  6 4  5 3  4 4  3 3  3 4 /
\plot 5 3  5 4 /

\setshadegrid span <.6mm>
\vshade 1.8 2 3  <,z,,> 4 0 1  <z,,,> 6.7 3 4 /

\endpicture} at 0 -11
\put{\beginpicture
\setcoordinatesystem units <.55cm,.45cm>
\put{$\Bbb E_8$}  at 0.2 3.5
\multiput{$\bullet$} at 0 0  1 1  2 0  3 1   4 0  5 1  6 0  7 1  8 0 
    2 2  4 2  6 2  5 2   /
\plot 0 0  1 1  2 0  3 1   4 0  5 1  6 0  7 1  8 0  /
\plot 1 1  2 2  3 1  4 2  5 1  6 2  7 1   /
\multiput{$\bullet$} at  3 2  4 1  4 2  3 0  /
\plot 3 1  3 2  4 1  5 2  5 1 /
\plot 3 0  4 1  4 0 /
\multiput{$\bullet$} at 2 3  3 3  4 3  5 3  6 3   /
\plot 4 2  3 3  2 2  2 3  3 2  4 3  4 2  5 3  6 2  6 3  5 2  4 3  /

\put{$\bullet$} at 4.3 4
\plot 4 3  4.3 4 /

\multiput{$\bullet$} at 10 0  9 1  8 2  7 3 /
\plot 6 2  7 3  10 0 /
\plot 8 0  10 2 /
\plot 7 1  9 3 /

\multiput{$\bullet$} at 3 4  4 4  5 4  6 4  7 4  /
\plot 2 3  3 4  4 3  5 4  6 3  7 4  7 3  6 4  5 3  4 4  3 3  3 4 /
\plot 5 3  5 4 /

\multiput{$\bullet$} at 12 0  11 1  10 2  9 3  8 4 /
\plot 7 3  8 4  12 0 /
\plot 10 0  11 1 /

\multiput{$\bullet$} at 3.3 5  4 5  5 5  5.3 5  6 5  7 5  8 5  /
\plot 3 4  4 5  5 4  6 5  7 4  8 5  8 4 /
\plot 4 5  4 4 /
\plot 6 5  6 4 /
\plot 4 4  5 5  6 4  7 5  8 4 /
\plot 3 4  3.3 5  4.3 4  5.3 5  5 4 /

\multiput{$\bullet$} at 2.3 6  4.3 6  5 6  5.6 6  6 6  6.3 6  7 6  /
\plot 5 5  6 6  7 5 /
\plot 4 5  5 6  6 5  7 6  8 5 /
\plot 5 5  5 6 /
\plot 4 5  4.3 6 /
\plot 2.3 6  3.3 5  4.3 6  5.3 5  6.3 6  6 5 /
\plot 5.3 5  5.6 6 /
\plot 7 5  7 6 /

\setshadegrid span <.6mm>
\vshade 1.8 2 3  <,z,,> 4 0 1  <z,,,>  8 4 5 /
\vshade 3 4 4  <z,z,,> 3.3 3.8 5  <z,z,,> 4 3 4.2 <z,z,,>
  4.2  3.2 4 <z,z,,> 6 5 5.7  
 <z,,,> 6.3   6 6 /

\endpicture} at 16 -11
\put{\beginpicture
\put{$\Bbb F_4$}  at -1.5 2.5
\multiput{$\bullet$} at 0 0  1 1  2 0  3 1   4 0  5 1  6 0   
    2 2  4 2  
    3 2  
    2 3  3 3  4 3 
    1 4  3 4  5 4  /
\plot 0 0  1 1  2 0  3 1   4 0  5 1  6 0  /
\plot 1 1  2 2  3 1  4 2  5 1  /
\plot 2 2  3 3  4 2  4 3  3 2  2 3  2 2 /
\plot 1 4  2 3  3 4  4 3  5 4 /
\plot 3 1  3 2 /
\plot 3 3  3 4 /
\setshadegrid span <.6mm>
\vshade 2 2 3  <,z,,> 3 1 2  <z,,,>  4 2 3  /
\endpicture} at 0 -18

\put{\beginpicture
\put{$\Bbb G_2$}  at -1.5 2.5
\put{}  at 7.5 2.5
\multiput{$\bullet$} at 0 0  1 1  2 0  1 2 1 3  1 4    /
\plot 0 0  1 1  2 0  /
\plot 1 1  1 4  /
\endpicture} at 16 -18

\endpicture}
$$
	\bigskip
It is well known that the number of rich antichains in $\Phi_+$ 
is equal to the cardinality of $\Phi_+$. It turns out 
that in general the poset structure of $\Cal R(\Delta)$ is quite 
similar to $\Phi_+(\Delta)$, but these posets are not always isomorphic:
For example, in case $\Bbb F_4$, the poset $\Lambda (\Phi_{e_2}')$ has a non-trivial
symmetry, thus also $\Cal R(\Bbb F_4)$, whereas $\Phi_+(\Bbb F_4)$ has no
non-trivial symmetry. 

	\medskip
{\bf Theorem 3.} {\it The poset $\Cal R(\Delta)$ of rich antichains in 
$\Phi_+(\Delta)$ is isomorphic to $\Phi_+(\Delta)$ if and only if 
$\Delta \neq \Bbb E_8,\ \Bbb F_4.$}
	\bigskip
{\bf Acknowledgment.} 
The paper is dedicated to Professor Liu Shao-Xue. 
His visit to Europe in 1985 was the start of a long lasting and
very fruitful cooperation between China and Germany devoted to
the representation theory of finite dimensional algebras. The core of this theory
are the Dynkin algebras and their indecomposable
representations, indexed by the elements of the corresponding root poset.
It seems to be surprising that these root posets which look like  
innocent combinatorial creatures still provide a lot of mysteries. 
The present paper and its successor [R3] try to illuminate some of their features. 

A first version [R1] of this paper was written in 2013
at SJTU Shanghai. Parts of the results have been presented at the 
57th Annual Meeting of the Australian Mathematical Society 2013 in Sydney, and in 
my ICRA lectures 2014 at Sanya, see [R2]. 
The author is grateful to many mathematicians for comments. 
In particular, he has to thank H\. Thomas for pointing out the references 
Humphreys [H, Theorem 3.20] and Armstrong [Ar, Theorem 5.4.1]).

Since the appendix of this paper has been used already quite frequently, we should
stress that we have changed the labeling of the vertices of $\Bbb B_n,\ \Bbb C_n$,
and $\Bbb D_n$  
in order to focus the attention (for $\Delta \neq \Bbb A_n$) 
to a special vertex which now is always labeled $c$.
	\bigskip\bigskip
{\bf 2. Solid subchain decompositions of $\Phi_+.$}
	\medskip 
We are going to prove Theorem 1: {\it Any root poset has a conical
decomposition.}
	\bigskip
		
Proof of Theorem 1.
The assertion is obvious for the cases $\Bbb A_n,$ $\Bbb B_n$
and $\Bbb G_2$.
Below we show solid subchain decompositions in the special 
cases $\Bbb D_5,\Bbb D_6,$ and for $\Bbb E_6,
\Bbb E_7, \Bbb E_8, \Bbb F_4$. We hope that
the cases $\Bbb  D_5$ and $\Bbb D_6$ show nicely the general rule how to obtain
a solid subchain decomposition in the cases $\Bbb  D_n$ in general. 

	\medskip 
Always, we use solid lines in order to specify the solid subchains.
The largest elements of the solid subchains are encircled. 

$$
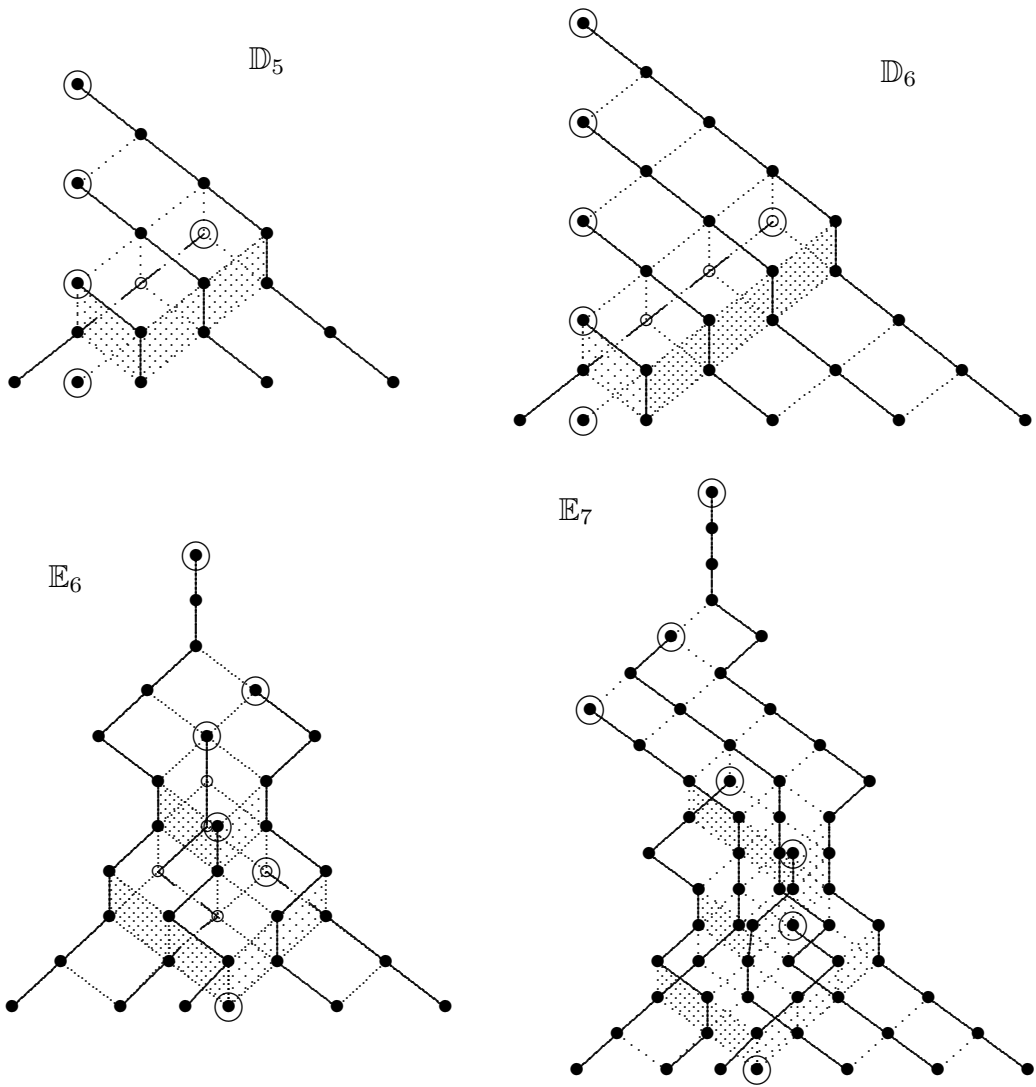
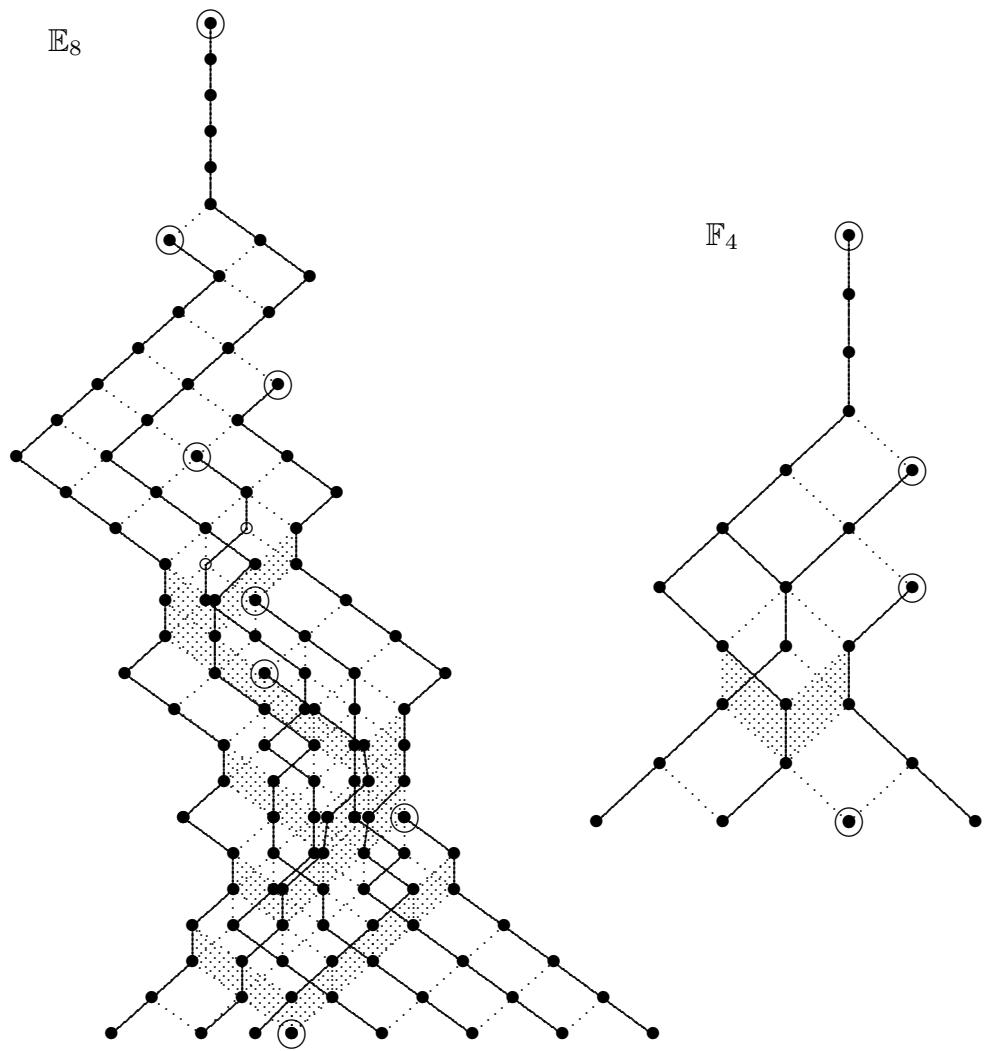
\hbox{\beginpicture
\setcoordinatesystem units <.7cm,.55cm>
\put{\beginpicture
\put{\bf $\Bbb D_5$} at -2 6.5 
\multiput{$\bullet$} at  
  0 0  -2 0  -4 0  -5 0  -6 0
  -1 1  -3 1  -4 1  -5 1
  -2 2  -3 2  -5 2
  -2 3  -4 3
  -3 4  -5 4
  -4 5
  -5 6 /
\multiput{$\bigcirc$} at -3 3  -5 0  -5 2  -5 4  -5 6 /   
          
\plot 0 0  -2 2  -2 3  -5 6 / 
\plot -2 0  -3 1  -3 2  -5 4 /
\plot -4 0  -4 1  -5 2 /
\plot -3 3  -3.3 2.7 /
\plot -3.7 2.3  -4.3 1.7 /
\plot -4.7 1.3  -6 0 /

\setdots <1mm>
\plot -5 4  -4 5 /
\plot -2 2  -4 0  -5 1  -5 2 /
\plot -2 3  -4 1 /
\plot -5 0  -2 3 /
\plot -5 1  -5 2  -3 4 / 
\plot -3 4  -3 3  -5 1 /
\plot -3 1  -4 2  -4 3 /  
\plot -2 2  -3 3 /

\multiput{$\circ$} at -3 3  -4 2 /

\setshadegrid span <.6mm>
\vshade -5 1 2 <,z,,> -4 0 1    <z,,,>  -2 2 3  /

\endpicture} at 0 0

\put{\beginpicture
\setsolid
\plot  0 0  -3 3  -3 4  -7 8 /
\plot -2 0  -4 2  -4 3  -7 6 /
\plot -4 0  -5 1  -5 2  -7 4 /
\plot -6 0  -6 1  -7 2 /
\plot -8 0  -6.7 1.3 /
\plot -6.3 1.7  -5.7 2.3 /
\plot -5.3 2.7  -4.7 3.3 /
\plot -4 4  -4.2 3.8 /
\multiput{$\bigcirc$} at -7 8  -7 6  -7 4  -7 2  -4 4  -7 0 /

\put{$\Bbb  D_6$} at -2 7
\setdots <.7mm> 
\multiput{$\bullet$} at 0 0  -2 0  -4 0  -6 0  -8 0  
     -1 1  -3 1  -5 1  -7 1
     -2 2  -4 2 
     -3 3 
      -6 1
     -5 2  -7 2
     -4 3  -6 3   
     -3 4  -5 4  -7 4 
     -4 5  -6 5 
     -5 6  -7 6 
     -6 7 
     -7 8 
     -7 0 /
               
\plot  0 0  -3 3  -6 0  -7 1  -8 0 /
\plot -1 1  -2 0  -4 2 /
\plot -2 2  -4 0  -5 1 /

 \plot -7 0  -3 4  -7 8 /
\plot -6 1  -7 2  -4 5 /
\plot -5 2  -7 4  -5 6 /
\plot -4 3  -7 6  -6 7 /
\plot -3 3  -3 4 /
\plot -4 2  -4 3 /
\plot -5 1  -5 2 /
\plot -6 0  -6 1 /
\plot -7 1  -7 2 /

\plot  -3 3  -4 4  -7 1 /
\plot -4 4  -4 5 /
\plot -4 2  -5 3  -5 4 /
\plot -5 1  -6 2  -6 3 /
   
\multiput{$\circ$} at -4 4  -5 3  -6 2 /
\setshadegrid span <.7mm>
\vshade -7 1 2  <z,,,>  -6 0 1 /
\setshadegrid span <.5mm>
\vshade  -6 0 1 <z,z,,>  -3 3 4  /

\endpicture} at 9 0
\endpicture}
$$

$$
\hbox{\beginpicture
\setcoordinatesystem units <.7cm,.55cm>
\put{\beginpicture
\setcoordinatesystem units <.6cm,.5cm>
\put{\bf $\Bbb  E_6$} at 1 9.5 
\multiput{$\bigcirc$} at  3.4 10  
   4.5 7  3.6 6  3.8 4  4.7 3    / 
\put{$\bigcirc$} at  4 0  

\setsolid
\plot 3.4 10  3.4 8  1.6 6  2.7 5  2.7 4  1.8 3  1.8 2  0 0 /
\plot  4.5 7  5.6 6  4.7 5  4.7 4  5.8 3  4.9 2  4.9 1  6 0 /
\plot 3.6 6  3.6 4  2.7 3   3 2.7 /
\plot 3.5 2.3  3.8 2  3.5 1.75  /
 
\plot 2 0  3.1 1.3  /

\plot 3.8 4  3.8 3  2.9 2  4 1  /
\plot 4.7 3  5.1 2.65 /
\plot 5.5 2.3  8 0 /

\plot 3.2 0  4 1 /

\setdots <.5mm> 

\multiput{$\circ$} at 3.8 2
       2.7 3  4.7 3
       3.6 4  
       3.6 5 /

\multiput{$\bullet$} at
  0 0  2 0  4 0  6 0  8 0
  0.9 1  2.9 1  4.9 1  6.9 1
  1.8 2  5.8 2  / 
\plot 0 0  1.8 2  /
\plot 5.8 2   8 0 /
\plot 0.9 1  2 0  2.8 1 /
\plot 1.8 2  4 0  5.8 2 /
\plot 6 0  6.9 1 /
\plot  2.7 3  4.9 1 /

\plot 1.8 2  3.6 4  5.8 2 / 
\plot 2.8 1  4.7 3 /
\plot 3.8 2   3.8 3  /
\plot 3.6 4  3.6 5 /
\plot 2.7 4   3.6 5  4.7 4 /
\plot 3.6 5  3.6 6  /

\multiput{$\bullet$} at 4 1
  2.9 2  4.9 2 
  1.8 3  3.8 3  5.8 3
  2.7 4  4.7 4  / %

\setdots <.5mm> 

\plot 4 0  4 1 /
\plot   2.9 1   2.9 2 /
\plot  1.8 2    1.8 3  /
\plot 5.8 2  5.8 3 /
\plot  2.7 3    2.7 4  /
\plot 4.7 3  4.7 4 /
\plot   4 1  5.8 3  4.7 4 /
\plot 2.7 4  
    1.8 3  4 1 / 
\plot 2.9 2  4.7 4 /
\plot  4.9 2   2.7 4  /

\multiput{$\bullet$} at 3.8 4
    2.7 5  4.7 5
    3.6 6 /
\plot  3.8 3  3.8 4 /
\plot 2.7 4  2.7 5  /
\plot 4.7 4  4.7 5 /
\plot 3.8 4
    2.7 5   3.6 6  4.7 5
   3.8 4 /

\setshadegrid span <.5mm>
\vshade 1.8 2 3  <z,z,,> 4 0 1   /
\vshade 2.7 4 5  <z,z,,> 3.8 3 4  /

\setshadegrid span <.7mm>
\vshade  4 0 1  <z,,,> 5.8 2 3 /
\vshade  3.8 3 4  <z,,,> 4.7 4 5 /

\multiput{$\bullet$} at 1.6 6  5.6 6  2.5 7  4.5 7  3.4 8  / 
\plot 2.7 5  1.6 6  3.4 8  5.6 6  4.7 5 /
\plot 2.5 7  3.6 6  4.5 7 / 

\multiput{$\bullet$} at 3.4 9  3.4 10  / 
\plot 3.4 8  3.4 10 /
  
\put{$\bullet$} at 3.2 0 
\plot 4 0  4 1 /

\endpicture} at 0 0

\put{\beginpicture
\setcoordinatesystem units <.5cm,.4cm>

\multiput{$\bullet$} at 3.8 2
       2.7 3  4.7 3      
       4.5 5   4.5 6 
       3.6 5 
       4.5 7  3.4 8 /
\multiput{$\bullet$} at
  0 0  2 0  4 0  6 0  8 0  10 0  
  0.9 1  2.9 1  4.9 1  6.9 1  8.9 1
  1.8 2  5.8 2  7.8 2
  6.7 3 /
  
\setdots <1mm>

\plot 0 0  1.8 2  /
\plot 5.8 2   8 0  8.9 1 /
\plot 0.9 1  2 0  2.8 1 /
\plot 1.8 2  4 0  5.8 2 /
\plot 4.9 1  6 0  7.8 2 /
\plot 5.8 2  6.7 3  10 0 /

\plot 1.8 2  3.6 4  5.8 2 / 
\plot 3.6 4  4.5 5  6.7 3 /
\plot 2.8 1  5.6  4 /
\plot  2.7 3  4.9 1 /
\plot 3.8 2   3.8 3  /
\plot 3.6 4  3.6 5 /
\plot 2.7 4   3.6 5  4.7 4 /
\plot 3.6 5  3.6 6  /
\plot 5.6 4  5.6 5 /
\plot 4.5 5  4.5 6 /
\plot 3.6 5  4.5 6  5.6 5 / 
\plot  4.5 6  4.5 7 /
\plot 2.5 7  3.4 8  5.6 6 /
\plot 3.6 6  4.5 7 / 
\plot 4.5 7  4.5 8 /

\plot 3.4 8  3.4 9 /

\multiput{$\bullet$} at 4 1
  2.9 2  4.9 2 
  1.8 3  3.8 3  5.8 3
   5.6 5 / 
\plot 4 0  4 1 /
\plot 2.9 1   2.9 2 /
\plot 4.9 1  4.9 2 /
\plot 1.8 2    1.8 3  /
\plot 5.8 2  5.8 3 /
\plot 2.7 3    2.7 4  /
\plot 4.7 3  4.7 4 /
\plot 4 1  5.8 3  4.7 4 /
\plot 2.7 4  
    1.8 3  4 1 / 
\plot 2.9 2  4.7 4 /
\plot  4.9 2   2.7 4  /
\plot 6.7 3  6.7 4  5.8 3  /
\plot 6.7 4  5.6 5  4.7 4 /
\plot 5.6 5  5.6 6 /

\multiput{$\bullet$} at 
    2.7 5  4.8 5
    3.6 6 /

\plot  3.8 3  3.9 4 /
\plot 2.7 4  2.7 5  /
\plot 4.8 4  4.8 5 /
\plot 3.9 4
    2.7 5   3.6 6  4.8 5
   3.9 4 /

\plot  5.6 6  4.8 5 /

\multiput{$\bullet$} at   0.3 10  1.4 9  2.5 8   3.6 7  4.8 6  
     1.2 11   2.3 10  3.4 9  4.5 8   5.6 7 
     2.1 12  3.2 11  4.3 10  5.4 9  6.5 8   
     3. 13  4.1 12
 /

\plot 2.1 12  3 13  4.1 12  3.2 11 /
\plot 0.3 10  4.7 6  6.5 8  2.1 12  0.3 10 /
\plot 1.2 11  5.6 7 /

\plot 1.4 9  3.2 11 /
\plot 2.5 8  4.3 10 /
\plot 3.6 7  5.4 9 /

\plot 4.8 5  4.7 6 /
\plot 3.6 6  3.6 7 /
\plot 5.6 6  5.6 7 /
\plot 2.5 7  2.5 8 /

\setsolid
\plot 0.3 10  3.6 7  3.6 5 / 
\plot 2.1 12  1.2 11 4.5 8  4.5 5 /
\plot 3 16   3 13  4.1 12  3.2 11  6.5 8  5.6 7  5.6 5 / 
\plot 2.7 5  1.6 6  3.4 8 /

\setdots <1mm>
\multiput{$\bullet$} at  3  14  3 15  3 16  /
\plot 3  13  3  16 /

\multiput{$\bullet$} at 1.6 6   5.6 6  2.5 7  /

\plot 2.7 5  1.6 6   2.5 7 /

\plot 2.5 7  3.6 6 /

\put{$\bullet$} at 3.2 0 

\multiput{$\bullet$} at  2.7 4  3.6 4  3.9 4  4.8 4  5.6 4  6.7 4 
  /
\multiput{$\bigcirc$} at  3 16  2.1 12  0.3 10  3.4 8   4 0 
   4.8 6   4.8 4 /
\setsolid

\plot  2.7 4  2.7 5 / 
\plot  3.6 5  3.6 4  /  
\plot 4.5 5  5.6 4 /
\plot 3.9 4  4.8 5  4.8 6 /   
\plot  5.6 5  6.7 4 /

\plot 2.7 4  1.8 3  2.9 2  2.9 1  2 0 /
\plot 3.6 4  0 0 /
\plot 3.9 4  3.8 3  3.8 2  6 0 /
\plot 4.7 4  5.8 3  4 1  3.2 0 /
\plot 5.6 4  4.7 3  8 0 /
\plot 6.7 4  6.7 3  10 0 /

\setshadegrid span <.5mm>
\vshade 1.8 2 3  <z,z,,> 4 0 1   /
\vshade 2.7 4 5  <z,z,,> 3.8 3 4  /
\vshade 2.5 7 8  <z,z,,> 4.8 5 6  /

\setshadegrid span <.7mm>
\vshade  4 0 1  <z,,,> 6.8 3 4 /
\vshade  3.8 3 4  <z,,,> 5.6 5 6 /
\vshade  4.8 5 6  <z,z,,> 5.6 6 7 /
\put{$\Bbb E_7$} at 0 15.5

\endpicture} at 9 0
\endpicture}
$$

$$
\hbox{\beginpicture
\setcoordinatesystem units <.5cm,.4cm>
\put{\beginpicture
\put{\bf $\Bbb E_8$} at -1 27.5

\put{$\bullet$} at 4 0 

\multiput{$\bigcirc$} at  
   2.2 28  1.3 22   3.7 18  1.9 16    3.2 12  3.4 10
    6.5 6  
   4 0   /
\multiput{$\bullet$} at    1.6 6  3.6 6   4.5 6  
   4.8 6   5.4 6    5.7 6   6.5 6  
  /

\setsolid
\plot 1.6 6  2.5 7  /
\plot 3.6 6  3.6 7 /
\plot  4.8 6  5.7 7 /
\plot   5.7 6  6.5 7 /
\plot   5.4 7  5.4 6  /  
\plot  4.5 7  4.5 6  /

\plot 1.6 6  2.7 5 /
\plot 3.6 6  3.6 5 /
\plot 4.5 6  4.5 5 /
\plot 4.8 6  4.7 5 /
\plot 5.4 6  6.5 5 /
\plot 5.7 6  5.6 5 /
\plot 6.5 6  7.6 5 /


\plot 2.7 5  2.7 4  1.8 3  1.8 2  0 0 /
\plot 3.6 5  4.7 4  4.7 3  8 0 /
\plot 4.5 5  3.6 4  2.7 3  6 0 /
\plot 4.7 5  3.8 4  3.8 3  2.9 2  2.9 1  2 0  /
\plot 6.5 5  5.6 4  10 0  /
\plot 5.6 5  6.7 4  4 1  3.2 0  /
\plot 7.6 5  7.6 4  12 0 /

\setdots <1mm>

\multiput{$\bullet$} at 3.8 2
       2.7 3  4.7 3
       3.6 4  5.6 4
       4.5 5  6.5 5
       5.4 6  
       3.6 5  
       4.5 6   
       5.4 7 
       4.5 7 
       3.4 8 
        5.4 8 
       4.3 9 
      2.1 12 
 3.2 11  4.3 10  5.4 9   /

\multiput{$\bullet$} at
  0 0  2 0  4 0  6 0  8 0  10 0  
  0.9 1  2.9 1  4.9 1  6.9 1  8.9 1
  1.8 2  5.8 2  7.8 2
  6.7 3 
  7.6 4  8.7 3  9.8 2  10.9 1  12 0 
   7.6 5  6.5 6  6.5 7 
/  
\plot 0 0  1.8 2  /
\plot 5.8 2   8 0  9.8 2 /
\plot 0.9 1  2 0  2.8 1 /
\plot 1.8 2  4 0  5.8 2 /
\plot 4.9 1  6 0  8.7 3 /
\plot 5.8 2  6.7 3  10 0  10.9 1 /
\plot  7.6 4  12 0 /
\plot 6.7 3  7.6 4  7.6 5  6.5 6  5.6 5 /
\plot 6.5 7  5.7 6 /
\plot 6.5 6  6.5 8 /
\plot 6.7 4 7.6 5 /

\plot 1.8 2  3.6 4  5.8 2 / 
\plot 3.6 4  4.5 5  6.7 3 /
\plot 2.8 1  5.6  4 /
\plot  2.7 3  4.9 1 /
\plot 3.8 2   3.8 3  /
\plot 3.6 4  3.6 5 /
\plot 2.7 4   3.6 5  4.7 4 /
\plot 3.6 5  3.6 6  /
\plot 5.6 4  5.6 5 /
\plot 4.5 5  4.5 6 /
\plot 3.6 5  4.5 6  5.6 5 / 
\plot  4.5 6  4.5 7 /
\plot 2.5 7  3.4 8  5.7 6 /
\plot 3.6 6  4.5 7 / 
\plot 4.5 7  4.5 8 /
\plot 3.4 8  3.4 9 /

\plot 7.6 4  5.4 6  4.5 5 /
\plot 6.5 6  6.5 5  5.6 4 /
\plot 4.5 6  5.4 7  6.5 6 /
\plot 5.4 7  5.4 6 /

\plot 5.4 7  5.4 9 /
\plot 4.5 7  5.4 8   6.5 7 /
\plot 3.4 8  4.3 9  5.4 8 /
\plot 4.3 9  4.3 10 /

\multiput{$\bullet$} at 4 1
  2.9 2  4.9 2 
  1.8 3  3.8 3  5.8 3
  2.7 4  4.7 4 
   6.7 4  5.6 5 / 
\plot 4 0  4 1 /
\plot 2.9 1   2.9 2 /
\plot 4.9 1  4.9 2 /
\plot 1.8 2    1.8 3  /
\plot 5.8 2  5.8 3 /
\plot 2.7 3    2.7 4  /
\plot 4.7 3  4.7 4 /
\plot 4 1  5.8 3  4.7 4 /
\plot 2.7 4  
    1.8 3  4 1 / 
\plot 2.9 2  4.7 4 /
\plot  4.9 2   2.7 4  /
\plot 6.7 3  6.7 4  5.8 3  /
\plot 6.7 4  5.6 5  4.7 4 /
\plot 5.6 5  5.7 6 /

\multiput{$\bullet$} at 3.8 4
    2.7 5  4.7 5
    3.6 6 /
\plot  3.8 3  3.8 4 /
\plot 2.7 4  2.7 5  /
\plot 4.7 4  4.7 5 /
\plot 3.8 4
    2.7 5   3.6 6  4.7 5
   3.8 4 /

\plot  5.7 6  4.7 5 /

\multiput{$\bullet$} at   0.3 10  1.4 9  2.5 8   3.6 7  4.8 6  
     1.2 11   2.3 10  3.4 9  4.5 8   5.7 7  6.5 8  /
\setsolid
\plot 2.3 12  3.1 13 /

\plot 2.5 7  2.5 8  0.3 10  1.2 11  1.2 13 
  -2.1 16  2.4 21  1.3 22 /   
\plot 3.6 7  4.5 8  2.3 10  2.3 12  
 /  
\plot 4.5 7  3.4 8  4.3 9  4.3 10  2.1 12  2.1 13  3 14  
    3 15  1.9 16  /  
\plot 5.4 7  5.4 10  3.2 12 / 
\plot 3.2 13 
   -.1 16   4.4 21   2.2 23  2.2 28 
   /  
\plot 5.7 7  5.6 8  3.4 10 /  
\plot  6.5 7  6.5 9  7.4 10  4.1 13  4.1 14  5 15 
   2.8 17  3.7 18  /  
\setdots <1mm>

\plot   1.2 11 0.3 10  4.8 6  6.5 8 /
\plot 1.2 11  5.7 7 /
\plot 1.4 9  2.3 10 /
\plot 2.5 8  3.4 9  /
\plot 3.6 7  4.5 8  /

\plot 4.7 5  4.8 6 /
\plot 3.6 6  3.6 7 /
\plot 5.7 6  5.7 7 /
\plot 2.5 7  2.5 8 /
\plot  1.2 11  2.1 12  6.5 8  /
\plot  2.3 10  3.2 11  3.2 12  /
\plot  3.4 9   4.3 10  4.3 11 /
\plot  4.5 8   5.4 9   5.4 10 /
\plot 2.1 12  2.1 13 /
\plot 1.2 12  2.1 13  3.2 12 /

\plot 3 14  3 15 /
\plot 2.1 14  2.1 13  3 14  4.1 13 /

\multiput{$\circ$} at 2.1 13  3.0 14  /
\multiput{$\bullet$} at    
    1.2 12  2.3 11  3.4 10  4.5 9  5.6 8 
    3.2 12  4.3 11  5.4 10  6.5 9 
    4.1 13  5.2 12  6.3 11  7.4 10 /
\plot 1.2 13  1.2 12  5.6 8  7.4 10  4.1 13  2.3 11  2.3 12  /
\plot 3.2 13  3.2 12  6.5 9 /
\plot  3.4 10  5.2 12 /
\plot 4.5 9  6.3 11 /

\plot 1.2 12  1.2 11 /
\plot 2.3 11  2.3 10 /
\plot 3.4 10  3.4 9 /
\plot 4.5  9  4.5 8 /
\plot 5.6  8  5.7 7 /
\plot 6.5  9  6.5 8 /

\plot 4.1 13  4.1 14 /

\multiput{$\bullet$} at 1.6 6   5.7 6  2.5 7  / 
\plot 2.7 5  1.6 6   2.5 7 /
\plot 2.5 7  3.6 6 /

\put{$\bullet$} at 3.2 0 

\multiput{$\bullet$} at 
  2.3 12  3.2 13  4.1 14  5 15 
  1.2 13  2.1 14  3 15  3.9 16 
  0.1 14  1 15  1.9 16  2.8 17  3.7 18 
  -1 15  -.1 16    .8 17  1.7 18  2.6 19  3.5 20  4.4 21 
  -2.1 16 -1.2 17 -.3 18  0.6 19  1.5 20    2.4 21  3.3 22               
                             1.3 22  2.2 23 /
\plot 3.3 22  -2.1 16  2.3 12  5 15  0.6 19 /

\multiput{$\bullet$} at  2.2 24  2.2  25  2.2 26  2.2 27  2.2 28 /
\plot 2.2 23  2.2 28  /
\plot 3.5 20  1.3 22  2.2 23  4.4 21  -1 15 /
\plot 0.1 14  3.7 18  1.5 20 /
\plot 1.2 13  3.9 16 /
\plot 3.2 13 -1.2 17 /
\plot 4.1 14  -.3 18 /
\setshadegrid span <.5mm>
\vshade 1.8 2 3  <,z,,> 4 0 1  <z,,,>  7.6 4 5 /
\vshade 2.7 4 5  <,z,,> 3.8 3 4  <z,,,> 6.5   6 7 /
\vshade 2.5 7 8  <,z,,> 4.7 5 6  <z,,,> 6.5   7 8 /
\vshade 1.2 11 12  <,z,,> 5.6 7 8  <z,,,> 6.5 8 9 /
\vshade 1.2 12 13  <,z,,> 2.3 11 12  <z,,,> 4.1 13 14 /
 
\endpicture} at 0 0
\put{\beginpicture
\setcoordinatesystem units <.7cm,.65cm>
\put{$\Bbb F_4$} at 2 10 
\put{$\bullet$} at 4 0 

\multiput{$\bullet$} at       1 4  3 4  5 4
  /
\multiput{$\bigcirc$} at       4 0  5 4  5 6  4 10 
 /

\plot 4 10  4 7  1 4  3 2   3 1  2 0 / 
\plot 5 6  3 4  3 3  0 0 /
\plot 2 5  3 4 /
\plot  5 4  4 3  4 2  6 0 /

\setdots <1mm>

\multiput{$\bullet$} at 
     0 0  2 0  4 0  6 0 
     1 1  3 1  5 1
     2 2  3 2  4 2 
     2 3  4 3
     2 5  4 5
     3 6  5 6
     4 7  4 8  4 9  4 10  /
\multiput{$\bullet$} at 3 3 /
\plot  0 0  2 2  4 0 5 1  /
\plot 1 1  2 0  4 2  6 0 /
\plot 3 1  3 2 /
\plot 2 2  2 3  4 5 /
\plot 4 2  4 3  2 5 /
\plot 3 2  5 4  3 6  1 4  3 2 /
\plot 3 6  4 7  5 6  4 5 /
\plot 4 7  4 10 /

\plot 2 2  3 3  4 2 /
\plot 3 3  3 4 /

\setshadegrid span <.5mm>
\vshade 2 2 3  <,z,,> 3 1 2  <z,,,>  4 2 3  /

\endpicture} at 10 0
\endpicture}
$$

	\medskip 
The largest elements of the solid subchains are always encircled (thus
these are suitable roots of height $\epsilon_i$, where $(\epsilon_1,\dots,\epsilon_n)$
is the exponent partition). 
\hfill$\square$

	\bigskip\bigskip 
Actually, the existence of a solid subchain decomposition 
concerns a local property, namely it concerns the bipartite
subposets $\Phi_{t,t+1}$ of all roots of height $t$ and $t+1$, we call these
subposets the {\it steps} of $\Phi_+$. As we know, we have
$|\Phi_{t+1}| \le |\Phi_{t}|$.
	\medskip 
The essential assertion of Theorem 1
is that for any $t\ge 1$,
there is a matching for $\Phi_{t,t+1}$ (an
injective map $f: \Phi_{t+1} \to \Phi_{t}$ such that $f(y) < y$ for all
$y \in \Phi_{t+1}$).

Namely, if we want to construct a conical decomposition, 
we may start at the top of the poset $\Phi_+$ and go down. If the 
subchains have reached the layer $\Phi_{t+1}$, we have to look at
$\Phi_{t,t+1}$ and we have to continue the path downwards inside a matching.
For example, in case $\Bbb  B_3$, starting with the maximal element $z$, the next two
choices for a solid chain containing $z$ are arbitrary, but then in $\Phi_{2,3}$ we have to be careful:
$$
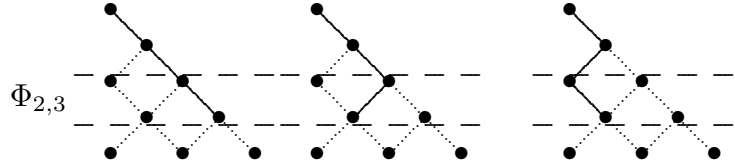
{\beginpicture
\setcoordinatesystem units <.4cm,.4cm> 
\put{\beginpicture
\multiput{$\bullet$} at 0 0  -1 1  -2 2  -3 3  -4 4  -2 0  -3 1  -4 2  -4 0  /
\setdots <.6mm>
\plot 0 0  -4 4 /
\plot -1 1  -2 0  -4 2  -3 3 /
\plot -2 2  -4 0 /
\setsolid 
\plot -4 4  -1 1 /
\setdashes <2mm>
\plot .5 0.8  -5 0.8 /
\plot .5 2.2  -5 2.2 /
\put{$\Phi_{2,3}$} at -6 1.5
\endpicture} at 0 0
\put{\beginpicture
\multiput{$\bullet$} at 0 0  -1 1  -2 2  -3 3  -4 4  -2 0  -3 1  -4 2  -4 0  /
\setdots <.6mm>
\plot 0 0  -4 4 /
\plot -1 1  -2 0  -4 2  -3 3 /
\plot -2 2  -4 0 /
\setsolid 
\plot -4 4  -2 2  -3 1 /
\setdashes <2mm>
\plot .5 0.8  -5 0.8 /
\plot .5 2.2  -5 2.2 /
\endpicture} at 7 0
\put{\beginpicture
\multiput{$\bullet$} at 0 0  -1 1  -2 2  -3 3  -4 4  -2 0  -3 1  -4 2  -4 0  /
\setdots <.6mm>
\plot 0 0  -4 4 /
\plot -1 1  -2 0  -4 2  -3 3 /
\plot -2 2  -4 0 /
\setsolid 
\plot -4 4  -3 3  -4 2  -3 1 /
\setdashes <2mm>
\plot .5 0.8  -5 0.8 /
\plot .5 2.2  -5 2.2 /
\endpicture} at 14 0
\endpicture}
$$
The choice in the middle does not work, since $\Phi_{2,3}$ has just one
matching namely:
$$
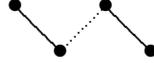
{\beginpicture
\setcoordinatesystem units <.5cm,.5cm> 
\multiput{$\bullet$} at -1 1  -2 2  -3 1  -4 2 /
\setdots <.6mm>
\plot -2 2  -3 1 /
\setsolid 
\plot -1 1  -2 2 /
\plot -3 1  -4 2 /
\endpicture}
$$
	\bigskip
{\bf Remark.}
As we have mentioned, for a conical poset $P$, the sequence $(h_1(P),h_2(P),\dots)$
is a Young partition. Let us stress that there are graded posets with $(h_1(P),h_2(P),\dots)$
being a Young partition, which cannot be written as the disjoint union of solid subchains 
such that each of the subchains contains a minimal element. Here is an example
(note that its width is greater than $h_1(P)$):
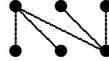
$$
{\beginpicture
\setcoordinatesystem units <.5cm,.5cm> 
\multiput{$\bullet$} at 0 0  1 0  2 0  0 1  1 1  2 1  /
\plot 0 0  0 1  2 0  2 1 /
\plot 1 0  0 1 /
\plot 1 1  2 0 /
\endpicture}
$$

	\bigskip\bigskip

{\bf 3. The rich antichains of $\Phi_+$.} 
	\medskip
This section is devoted to the proof of Theorem 2.

We say that a poset $P$ is {\it bipartite} provided any element of $P$ is minimal
or maximal, but not both. If $P$ is bipartite, let $P_1$ be the set
of minimal elements and $P_2$ the set of maximal elements. A {\it matching} of $P$
is an injective map $f\:P_2 \to P_1$ with $f(y) < y$ for all $y\in P_2$. 

We say that $P$ is an {\it $M$-poset} provided $P$ is bipartite with $|P_1| = |P_2|+1$
and such that for any $y\in P_2$ there is a matching $f$ such that the unique
element $x\in P_1$ which is not in the image of $f$ satisfies $x < y.$
Note that if $P$ is an $M$-poset, then no element in $P_2$ is join-irreducible.
	\medskip 

If $T$ is a finite tree, we define $M(T)$ as the incidence poset of $T$
(with $M(T)_1$ the set
of vertices, $M(T)_2$ the set of edges, such that $x < y$ provided $x$ is a vertex
on the edge $y$). 
	\medskip
For example:
$$
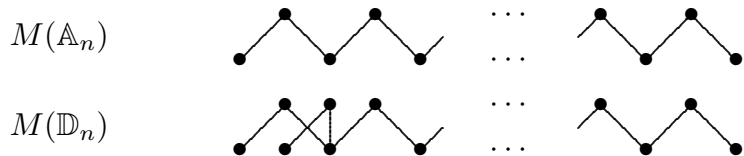
{\beginpicture
\setcoordinatesystem units <.5cm,.5cm>
\put{\beginpicture
\multiput{$\bullet$} at 0 0  2 0  4 0  8 1  10 1  1 1  3 1  9 0  11 0 /
\multiput{$\cdots$} at 6 0  6 1 /
\plot 0 0  1 1  2 0  3 1  4 0  4.5 0.5 /
\plot 11 0  10 1  9 0  8 1  7.5 0.5 /
\put{$M(\Bbb A_n)$} at -4 0.5
\endpicture} at 0 0
\put{\beginpicture
\multiput{$\bullet$} at 0 0  2 0  4 0  8 1  10 1   1 1  3 1  9 0  11 0  1 0  2 1 /
\multiput{$\cdots$} at 6 0  6 1 /
\plot 0 0  1 1  2 0  3 1  4 0  4.5 0.5 /
\plot 11 0  10 1  9 0  8 1  7.5 0.5 /
\plot 1 0  2 1  2 0 /
\put{$M(\Bbb D_n)$} at -4 0.5

\endpicture} at 0 -2
\endpicture}
$$
Here, $n$ is the number vertices of $T$, thus the number of minimal elements
of the poset $M(T)$. 

	\bigskip 
{\bf Proposition 1.} {\it If $T$ is a finite tree, then $M(T)$ is an $M$-poset.}
	\medskip
Proof. Assume that $T$ has $n$ vertices, thus $n-1$ edges. It follows that
$M(T)$ is bipartite with $M(T)_2 = n-1$ and $M(T)_1 = n$. 
Let $y$ be an edge, and $x$ a vertex on $y$. For every edge $y'$, let $f(y')$
be the vertex on $y'$ with maximal distance to $x$. Then $f$ is a matching, and 
$x$ is not in the image of $f$. Of course, $x < y.$ 
		\medskip
	\medskip
{\bf Corollary.} {\it For any root poset $\Phi_+ = \Phi_+(\Delta)$, the poset
$\Phi_{12}$ is an $M$-poset.}
	\medskip
Proof: Obviously, $\Phi_{12} = M(\Delta).$
	\bigskip 
{\bf Proposition 2.} {\it Assume that $\Delta$ has rank at least $3$, let
$\Phi = \Phi(\Delta)$ and $e = e_2(\Delta).$
If $\Delta \neq \Bbb D_4,\ \Bbb E_6,$ then $\Phi_{e,e+1}$ is an $M$-poset.}

	\medskip
We claim that $\Phi_{e,e+1}$ is obtained from a poset of the
form $P = M(T)$ by adding for some pairs $x\in P_1,\ y\in P_2$ 
the relation $x < y$. Since $M(T)$ is an $M$-poset, also $\Phi_{e,e+1}$ has
to be an $M$-poset. Here are drawings of the various cases. On the right, we
exhibit $P = M(T)$. In the drawings, the solid lines show the relations 
in $M(T)$, the added relations are dashed. 
$$
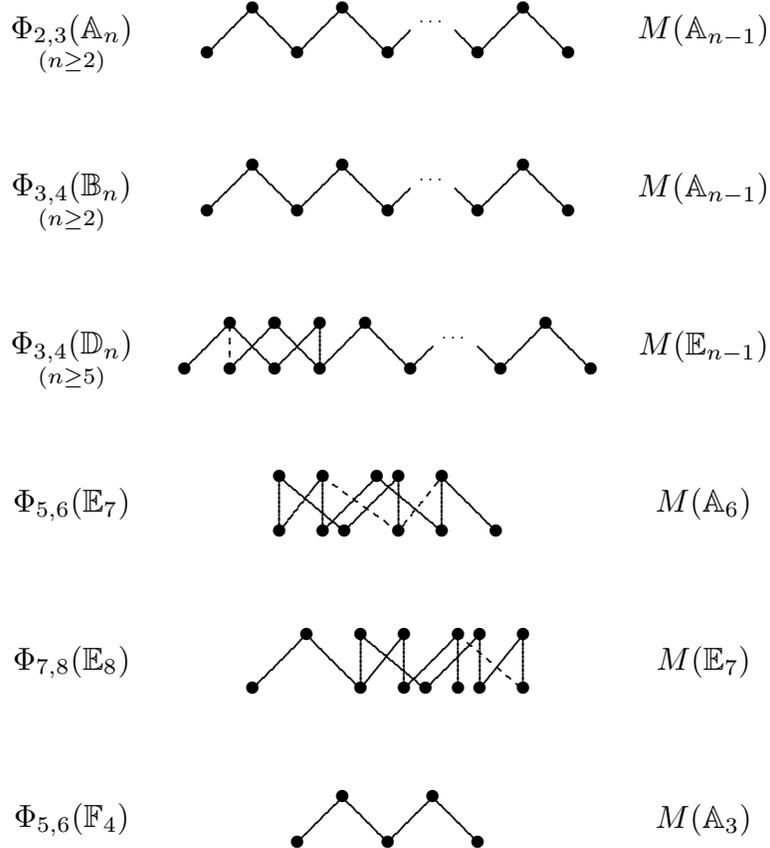
\hbox{\beginpicture
\setcoordinatesystem units <.7cm,.7cm>

\put{\beginpicture
\setcoordinatesystem units <.5cm,.5cm>
\multiput{$\bullet$} at 0 0  2 0  4 0  6 0  8 0 /
\multiput{$\bullet$} at 1 1  3 1  7 1  /
\plot 0 0  1 1  2 0  3 1  4 0  4.5 0.5 /
\plot 5.5 0.5  6 0  7 1  8 0 /
\setdots <1mm>
\plot 4.75 0.7  5.25 0.7 /
\endpicture} at 0 5
\put{$\Phi_{2,3}(\Bbb A_n)$} at -5 5
\put{$\ssize (n\ge 2)$} at -5 4.5
\put{$M(\Bbb A_{n-1})$} at 5 5 

\put{\beginpicture
\setcoordinatesystem units <.5cm,.5cm>
\multiput{$\bullet$} at 0 0  2 0  4 0  6 0  8 0 /
\multiput{$\bullet$} at 1 1  3 1  7 1  /
\plot 0 0  1 1  2 0  3 1  4 0  4.5 0.5 /
\plot 5.5 0.5  6 0  7 1  8 0 /
\setdots <1mm>
\plot 4.75 0.7  5.25 0.7 /
\endpicture} at 0 2.5
\put{$\Phi_{3,4}(\Bbb B_n)$} at -5 2.5
\put{$\ssize (n\ge 2)$} at -5 2
\put{$M(\Bbb A_{n-1})$} at 5 2.5 

\put{\beginpicture
\setcoordinatesystem units <.5cm,.5cm>
\multiput{$\bullet$} at -1 0  0 0  1 0  2 0  4 0  6 0  8 0 /
\multiput{$\bullet$} at  0 1  1 1  2 1  1 1  3 1  7 1  /
\plot 0 0  1 1  2 0  3 1  4 0  4.5 0.5 /
\plot 5.5 0.5  6 0  7 1  8 0 /
\plot -1 0  0 1  1 0  2 1  2 0 /
\setdots <1mm>
\plot 4.75 0.7  5.25 0.7 /
\setdashes <.7mm>
\plot 0 0  0 1 /
\endpicture} at 0 0
\put{$\Phi_{3,4}(\Bbb D_n)$} at -5 0
\put{$\ssize (n\ge 5)$} at -5 -.5
\put{$M(\Bbb E_{n-1})$} at 5 0

\put{\beginpicture
\setcoordinatesystem units <.6cm,.6cm>
\multiput{$\bullet$} at 0 0  0.8 0  1.2 0  2.2 0  3 0  4 0 /
\multiput{$\bullet$} at 0 1  0.8 1  1.8 1  2.2 1  3 1  /
\plot 4 0  3 1  3 0  1.8 1 
  0.8 0  0.8 1  0 0  0 1  1.2 0  2.2 1  2.2 0  /
\setdashes <.7mm>
\plot 0.8 1  2.2 0  3 1 /
\endpicture} at 0 -2.5
\put{$\Phi_{5,6}(\Bbb E_7)$} at -5 -2.5
\put{$M(\Bbb A_6)$} at 5 -2.5 
\put{\beginpicture
\setcoordinatesystem units <.6cm,.6cm>
\multiput{$\bullet$} at 0 0  2 0  2.8 0  3.2 0  3.8 0  4.2 0  5 0 /
\multiput{$\bullet$} at 1 1  2 1  2.8 1  3.8 1  4.2 1  5 1   /
\plot 0 0  1 1  2 0  2 1  3.2 0  4.2 1  4.2 0  5 1  5 0 /
\plot 3.8 0  3.8 1  2.8 0  2.8 1  2 0 /
\setdashes <.7mm>
\plot 3.8 1  5 0  /
\endpicture} at 0 -5
\put{$\Phi_{7,8}(\Bbb E_8)$} at -5 -5
\put{$M(\Bbb E_7)$} at 5 -5 
\put{\beginpicture
\setcoordinatesystem units <.5cm,.5cm>
\multiput{$\bullet$} at 0 0  2 0  4 0  /
\multiput{$\bullet$} at 1 1  3 1   /
\plot 0 0  1 1  2 0  3 1  4 0 /
\endpicture} at 0 -7.5
\put{$\Phi_{5,6}(\Bbb F_4)$} at -5 -7.5
\put{$M(\Bbb A_3)$} at 5 -7.5 
\endpicture}
$$
	\medskip 
Proof of Theorem 2. If the rank $n$ of $\Delta$ is equal to 2, then
$\Cal R(\Delta) = \Phi_+$. Thus, we can assume that $n\ge 3.$ For $\Delta =
\Bbb D_4$, $\Phi_4(\Delta)$ is just the maximal element of $\Phi_+$, and
this element cannot belong to any antichain of cardinality greater than 1. 
Thus, we assume now that $\Delta\neq \Bbb D_4.$

The case $\Bbb E_6$ is special (since $\Phi_5(\Bbb E_6)$ contains a join-irreducible
element; in particular, $\Phi_{45}(\Bbb E_6)$ is not an $M$-poset) and has to be
treated separately. 
Thus, we assume that $\Delta$ has rank at least $3$ and is different from
$\Bbb D_4,\ \Bbb E_6$. This means that we can apply Proposition 2.

Let $z$ be a root of length $t \ge e+1.$ Let $\Phi(z)$ be the set of positive
roots which are not comparable with $z$. We claim that $\Phi(z)$ has width at most
$n-3$ (this shows that $z$ cannot belong to an antichain of cardinality $n-1$).

Let $\Cal C$ be a conical decomposition of $\Phi_+$. Let $C\in \Cal C$ with 
$z\in C.$ Since the length of $z$ is at least $e+1$, the chain $C$
contains a root $y$ of length $e+1$. Since $\Phi_{e,e+1}$ is an $M$-poset,
there is a matching $f\:\Phi_{e+1} \to \Phi_e$ such that the unique element $x$
$\in \Phi_e$ which does not belong to the image of $f$ satisfies $x < y.$

Let $v$ be the root of length 2 which belongs to $C$.
Since $\Phi_{12}$ is an $M$-poset, there is a matching $g\: \Phi_2 \to \Phi_1$
such that the simple root $u$ which 
does not belong to the image of $g$, satisfies $u < v$.
Let $\Cal D$ be obtained from
the decomposition $\Cal C$ by replacing in any chain $C'\in \Cal C$ of cardinality at
least 2 its simple root by the simple root $g(v')$, where $v'$ is the root of length
$2$ in $C'$ and using the singleton $\{u\}$ as the unique chain of cardinality 1
in $\Cal D$. Thus, $\Cal D$ is again a conical decomposition of $\Phi_+$.

Now, let $y'\in \Phi_{e+1}$. If $y'$ belongs to $D\in \Cal D$,
let $J_{y'}$ be the set of elements $z'\in D$ with $y' \le z'.$ 

Similarly, assume that $x'\in \Phi_e$ belongs to $D\in \Cal D$ and let 
$I_{x'}$ be the set of elements $w\in D$ with $w \le x'.$ 
Since $f(y') < y'$, we see that the union $I_{f(y')}\cup J_{y'}$ is a chain.

We have decomposed $\Phi_+$ into the chains $I_{x'}, J_{y'}$
with $x'\in \Phi_e$ and $y'\in \Phi_{e+1}$, as well as the singleton $\{u\}$.
We know that the elements in $J_y,\ I_{f(y)},\ I_x$ as well as the element
$u$ all are comparable with $z$, thus  they belong to $\Phi(z)$.

It follows that $\Phi\setminus \Phi(z)$ is covered by the chains
$I_{f(y')}\cup J_{y'}$, 
where the elements $y'$ are the elements of $\Phi_{e+1} \setminus \{y\}.$ 
The cardinality
of $\Phi_{e+1}$ is $n-2$, thus there are $n-3$ elements of the form
$y'$. This shows that $\Phi\setminus \Phi(z)$ is covered by $n-3$ chains.

	\bigskip
It remains to look at the case $\Delta = \Bbb E_6$. We denote by $a'$ the
join-irreducible element of $\Phi_5$, and by $a$ its neighbor in $\Phi_4.$
Thus $\Phi'_4$ is obtained from $\Phi_4$ by replacing $a$ by $a'$. 

Let us look at $\Phi_{56}(\Bbb E_6)$. There is the following matching of 
$\Phi_{56}(\Bbb E_6)$:
$$
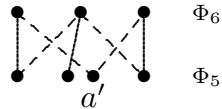
\hbox{\beginpicture
\setcoordinatesystem units <.7cm,.7cm>
\multiput{$\bullet$} at 0 0  0.8 0  1.2 0  2 0  0 1  1 1  2 1 /
\plot 0 0  0 1 /
\plot 0.8 0  1 1 /
\plot 2 0  2 1 /
\setdashes <1mm>
\plot 0 0  1 1  2 0 /
\plot 0 1  1.2 0  2 1 /
\put{$a'$} at 1.2 -0.3
\put{$\ssize \Phi_6$} at 3 1
\put{$\ssize \Phi_5$} at 3 0
\endpicture}
$$
It follows that there is a conical decomposition $\Cal C$ of $\Phi_+$
such that $a'$ is the maximal element of one of the chains. Also, we can assume that
the simple root corresponding to the branching vertex $c$ is the singleton in
$\Cal C$. 

Let $z$ be a root of length at least 5 and different from $a'$. We want to show
that $z$ does not belong to a rich antichain. As above, we denote by $\Phi(z)$ 
the set of positive roots which are not comparable with $z$. We claim that $\Phi(z)$ has width at most $3\ (= n-3)$. Note that 
the simple root corresponding to the branching vertex $c$ belongs to
$\Phi(z)$. Thus, only the five chains in $\Cal C$ which are not singletons 
have to be considered. 

Here is $\Phi_{45}$. Any element of $\Phi_5\setminus\{a'\}$
has been connected to two elements of $\Phi_4.$
$$
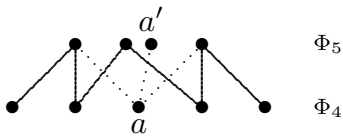
\hbox{\beginpicture
\setcoordinatesystem units <.7cm,.7cm>
\put{\beginpicture
\put{$a'$} at 2.2 1.4
\put{$a$} at 2 -.3
\put{$\ssize \Phi_5$} at 5 1
\put{$\ssize \Phi_4$} at 5 0
\multiput{$\bullet$} at 0 0  1 0  2.2 1  3 0  4 0 /
\multiput{$\bullet$} at 2 0  1 1  1.8 1  3 1 /
\plot 0 0  1 1  1 0  1.8 1  3 0  3 1  4 0 /
\setdots <1mm>
\plot 2 0  2.2 1 /
\plot 1 1  2 0  3 1 /
\endpicture} at 0 0
\endpicture}
$$
As above, we use $\Cal C$ in order to define chains $J_y$ for $y\in \Phi_5$
(with $J_{a'} =\{a'\}$) as well as chains $I_x$ for $x\in \Phi_4.$ 
It follows that the elements outside of $\Phi(z)$ belong to three
chains which are combined from three of the chains of the form $J_y$ and three 
of the chains of the form $I_x$.

This completes the proof. \hfill$\square$
	\bigskip\bigskip
{\bf 4. The poset $\Cal R(\Delta)$ of rich antichains of $\Phi_+(\Delta)$.}
     \medskip	  
If $P$ is a finite poset and $t$ a non-negative integer, let $\Cal A_t(P)$ be the
set of antichains in $P$ of cardinality $t$. For a poset $P$ of width $n$,
we call the antichains of cardinality $n-1$ the rich antichains. Thus, for $\Delta$
a Dynkin diagram of rank $n$, we write 
$\Cal R(\Delta) = \Cal A_{n-1}(\Phi_+(\Delta))$.

Given a Dynkin diagram $\Delta$ of rank $n$ with root poset $\Phi_+$, one knows that
$|\Cal A_t(\Phi_+)| = |\Cal A_{n-t}(\Phi_+)|$, for $0 \le t \le n$, see [At]. In particular, 
since $A_1(\Phi_+) = \Phi_+$, we always have 
$$
 |\Cal R(\Delta)| = |\Phi_+|,
$$
and one may ask whether
it is possible to recover also the partial ordering of $\Phi_+(\Delta)$  
by looking at the set $\Cal R(\Delta)$ of rich antichains in $\Phi_+(\Delta).$
   \medskip
Let us point out, that {\it it is not possible in general 
to recover the partial ordering
of $\Phi_+$  by looking at the set of rich antichains.} 
Namely, consider the Dynkin type $\Bbb F_4$. According to Theorem 2,
$\Cal R(\Bbb F_4)$ is the poset of rich antichains in $\Lambda(\Phi_5(\Bbb F_4))$.
Now the poset $\Lambda(\Phi_5(\Bbb F_4))$ has an automorphism $\phi$
of order 2, and $\phi$ induces a non-trivial automorphism on 
$\Cal R(\Bbb F_4)$
(since obviously there are rich antichains which are not invariant under $\phi$),
whereas $\Phi_+(\Bbb F_4)$ itself has no non-trivial automorphisms.
	\bigskip
Note that given a poset $P$, 
there are several ways to consider $\Cal A_t(P)$ as
a poset: given two antichains $x,y$, we write $x\le_\Lambda y$ provided $x$ lies in the
ideal generated by $y$; similarly, we write $x\le_V y$ provided $y$ lies in the
coideal generated by $x$. Of course, there is also the possibility to combine the
two partial orderings, namely to set
$x \le y$ provided both $x\le_\Lambda y$ and $x\le_V y$, thus provided both
$\Lambda(x) \subseteq \Lambda(y)$ and $V(x) \supseteq V(y).$ Actually, it turns out
that for $\Cal R(\Delta)$, there is the following observation: 
         \bigskip
	 {\bf Proposition.}
{\it Let $A$ and $B$ be rich antichains of $\Phi_+$, with
$A \le_V B$. Then we can label the elements of $A$ and $B$ as $A = \{a_1,\dots,a_{n-1}\},$
$B = \{b_1,\dots,b_{n-1}\}$ with $a_i \le b_i$ for $1\le i \le n-1$. As a consequence,
$A \le_\Lambda B$.}
   \medskip
Proof. Let $\Phi_+ = \Phi_+(\Delta)$, where $\Delta$ is a (connected) Dynkin diagram of rank $n$.
We assume that $A,B$ are rich antichains of $\Phi_+$ with $A \le_V B.$ Let $z$ be a simple
root which does not belong to $A$. Since $B \subseteq V(A)$, it follows that $z$ does not belong to
$B$. According to the Addendum of Theorem 1, there is a conical decomposition $\Cal C$ such that
$\{z\}$ is one of the elements of $\Cal C$. Let $C_1,\dots,C_{n-1}$ be the remaining chains in
$\Cal C.$ Since $z$ does not belong to $A$, the elements of $A$ belong to the chains $C_i$
with $1\le i \le n-1$, the elements $a$ of $A$ can be labeled as $a_1,\dots,a_{i-1}$ with
$a_i\in C_i$, for $1\le i \le n-1$. Similarly, the elements of $B$
can be labeled as $b_1,\dots,b_{i-1}$ with
$b_i\in C_i$, for $1\le i \le n-1$. We claim that $a_i \le b_i$ for all $i$. Otherwise, there is some
$i$ with $a_i > b_i$. Since $B \subseteq V(A)$, there is some $a_j\in A$ with $b \ge a_j.$
But $a_i > b_i \ge a_j$ is impossible, since $A$ is an antichain. This is the first assertion.

Of course, if $A = \{a_1,\dots,a_{n-1}\},$
$B = \{b_1,\dots,b_{n-1}\}$ with $a_i \le b_i$ for all $i$, then  $A \le_\Lambda B$ (as well as
$A \le_V B$).
\hfill$\square$
	\medskip
We always consider $\Cal R(\Delta)$ as a poset using the partial ordering $\le$ (having in mind that this is the same as the partial ordering $\le_V$). 
     \bigskip
{\bf Remark.} If $A$ and $B$ are rich antichains of $\Phi_+$ with
$A \le_\Lambda B$, then we may not have $A \le_V B.$ Here is the first example: $\Delta = \Bbb A_3$ and
$$
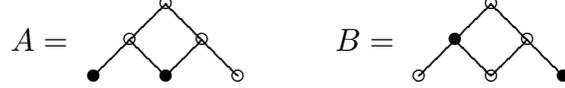
{\beginpicture
\setcoordinatesystem units <.4cm,.4cm>
\put{\beginpicture
\put{$A =$} at -1.5 1
\multiput{$\circ$} at 1 1  2 2  3 1  4 0 /
\multiput{$\bullet$} at 0 0  2 0 /
\plot 0 0  2 2  4 0 /
\plot 1 1  2 0  3 1 /

\endpicture} at 0 0
\put{\beginpicture
\put{$B =$} at -1.5 1
\multiput{$\circ$} at  2 2  3 1  0 0  2 0   /
\multiput{$\bullet$} at 1 1  4 0 /
\plot 0 0  2 2  4 0 /
\plot 1 1  2 0  3 1 /

\endpicture} at 9 0
\endpicture}
$$
	   \bigskip
{\bf Proof of Theorem 3.} 
We are going to present only the essential steps. But we hope that the arguments which
we provide shed some further light on the structure of the individual root posets.
In case $\Delta$ has rank 2, the rich antichains of $\Phi_+$ are just the singletons,
thus in this case $\Phi_+(\Delta) = \Cal R(\Delta)$. In particular, we do not have to
consider the case $\Bbb G_2.$ 
        \medskip
If $\Delta$ is a Dynkin diagram and $t\ge 1$, we denote by $\Phi(\Delta,t)$ the set
of positive roots with coefficients bounded by $t$ such that at least one coefficient 
is equal to $t$. Of course, $\Phi_+$ is the disjoint union of the subsets
$\Phi(\Delta,t)$ with $t\ge 1.$
Note that $\Phi(\Delta,1)$ is always an ideal of $\Phi_+$ (the ideal generated by
the minimal sincere root: all its coefficients are equal to $1$), whereas the
remaining non-empty subposets $\Phi(\Delta,t)$ are intervals.
It should be stressed that for $\Delta$ equal to $\Bbb B_n$ and $\Bbb C_n$, the 
subsets $\Phi(\Delta,2)$ of $\Phi_+$ are the same despite the fact that for any
element $x$ of $\Phi_+$ its coefficients may 
depend on whether we deal with $\Bbb B_n$ or $\Bbb C_n$.
	\medskip
We are trying to exhibit maps $R = R(\Delta,t)\:\Phi(\Delta,t) \to \Cal R(\Delta)$ which 
combine to a bijection $\Phi_+ \to \Cal R(\Delta).$ All our maps $R(\Delta,t)$ will be injective
and the images of $R(\Delta,t)$ and $R(\Delta,t')$ will be disjoint for $t\neq t'$. 
Thus, if for some $\Delta$, the maps $R(\Delta,t)$ are defined for all $t$ (this will
be the case for $\Delta = \Bbb A_n, \Bbb B_n, \Bbb C_n, \Bbb D_n, \Bbb E_6$ and $\Bbb F_4$),
then the union of the images has to be all of $\Cal R(\Delta)$, since it is well-known
that $|\Cal R(\Delta)| = |\Phi_+|$. Thus, we do not bother to verify surjectivity assertions. 
	\bigskip
{\bf The bijections $R = R(\Delta,1)\:\Phi(\Delta,1) \to \Cal R(\Phi_{12}).$}
	\medskip
Let $r$ be in $\Phi(\Delta,1)$, thus $r$ is uniquely determined by its support
$\supp r$, this is the set of simple roots $x$ with $x \le r,$ it is a connected subdiagram of
$\Delta$. We define $R(r)$ as follows:
$$
 R(r) = R(r)_1\cup (r)_2, \quad \text{where}\quad R(r)_1 = \{x\in \Phi_1\mid x \not\le r\},\
R(r)_2  = \{y\in \Phi_2\mid y\le r\}.
$$
It is easy to see that $R(r)$ is an antichain: both $R(r)_1$ and $R(r)_2$ are antichains,
and if $x\in R(r)_1$ and $y\in R(r)_2$, then $x \not< y$ (since $x < y$ would 
imply $x < y \le r$, but $x \not\le r$). Also, $R(r)$ is rich (namely, if
$r$ has length $m$ and $\Delta$ has rank $n$, then $R(r)_1$ has cardinality $n-m$ and
$R(r)_2$ has cardinality $m-1$).
	\medskip
Finally, we claim: {\it $R$ preserves and reflects the partial orderings} (if
$r,r'$ are roots in $\Phi(\Delta,1)$, then $r\le r'$ if and only if $R(r') \le_V R(r)).$ 

Here is the proof. Of course, $r \le r'$ if and only if $\supp r \subseteq \supp r'.$
First, assume that $\supp r \subseteq \supp r'$. We want to show that $R(r') \subseteq 
V(R(r)).$ We have $R(r')_1 \subseteq R(r)_2$, thus, it remains to show that any element
$y\in R(r')_2$ belongs to $V(R(r)).$ Let $a,b$ be the support of $y$. If $a$ is not
in the support of $r$, then $a \in R(r)$, thus $y \in V(R(r))$. Similarly, if
$b$ is not in the support of $r$, then $y \in V(R(r))$. Finally, if both $a,b$
are in the support of $r$, then $y\in R(r)_2 \subseteq V(R(r))$. 

Second, assume that $R(r') \subseteq V(R(r)).$ We have to show that $\supp r \subseteq \supp r'.$
Let $x$ belong to $\supp r$ and assume that $x$ does not belong to $\supp r'.$
Then $x \in R(r')_1 \subseteq V(R(r)).$ Since $x$ is a minimal element, it follows that
$x\in R(r)$, thus $x$ belongs to $R(r)_1$ and therefore $x\notin \supp r,$ a contradiction.  
\hfill$\square$
	\bigskip
{\bf Case} $\Delta = \Bbb A_n$. In this case, $\Phi_+ = \Phi(\Delta,1)$, thus $R = R(\Delta,1)$
provides an isomorphism between $\Phi_+(\Bbb A_n)$ and $\Cal R(\Bbb A_n).$
	\bigskip
Let us now consider the posets $\Phi(\Delta,t)$ with $t\ge 2$. Whereas for $t=1$, we were able to provide a general recipe, we now have to work case by case. 
	\medskip
{\bf Case} $\Delta = \Bbb B_n$. The set $\Phi(\Bbb B_n,2)$ is just the coideal
of $\Phi_+$ generated by the maximal root of $\Bbb B_2$ (considered as a subdiagram
of $\Bbb B_n$). 

We start with the poset $P = \Phi(\Bbb B_n)_{13}$, since the rich
antichains of $\Phi_+(\Bbb B_n)$ are contained in $P$.
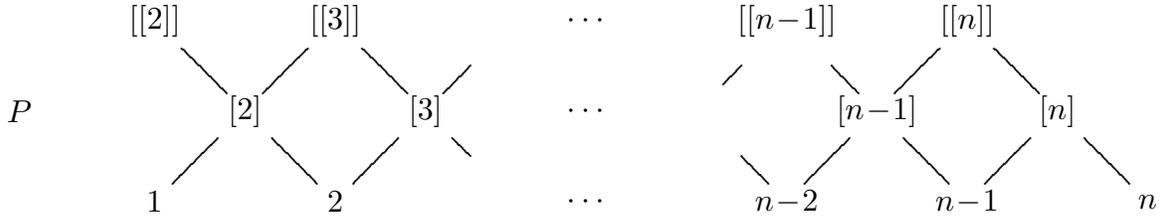
$$
{\beginpicture
\setcoordinatesystem units <1cm,1cm>
\put{$P$} at -12.5 1
\put{$n$} at 0 0
\put{$n\!-\!1$} at -2 0
\put{$n\!-\!2$} at -4 0
\put{$2$} at -9 0
\put{$1$} at -11 0
\put{$[n]$} at -1 1
\put{$[n\!-\!1]$} at -3 1
\put{$[3]$} at -8 1
\put{$[2]$} at -10 1
\put{$[[n]]$} at -2 2
\put{$[[n\!-\!1]]$} at -4 2
\put{$[[3]]$} at -9 2
\put{$[[2]]$} at -11 2
\multiput{$\cdots$} at -6.2 0  -6.2 1  -6.2 2 /
\plot -0.2 0.2 -0.7 0.7 /
\plot -1.3 0.7 -1.8 0.2 /
\plot -2.2 0.2 -2.7 0.7 /
\plot -3.3 0.7 -3.8 0.2 /
\plot -4.2 0.2 -4.5 0.5 /

\plot -7.5 0.5 -7.7 0.7 /
\plot -8.3 0.7 -8.8 0.2 /
\plot -9.2 0.2 -9.7 0.7 /
\plot -10.3 0.7 -10.8 0.2 /

\plot -1.2 1.2 -1.7 1.7 /
\plot -2.3 1.7 -2.8 1.2 /
\plot -3.2 1.2 -3.5 1.5 /

\plot -4.7 1.3 -4.5 1.5 /

\plot -7.5 1.5 -7.8 1.2 /
\plot -8.2 1.2 -8.7 1.7 /
\plot -9.3 1.7 -9.8 1.2 /
\plot -10.2 1.2 -10.7 1.7 /

\endpicture}
$$
The rich antichains in $P_{12}$ form the poset $\Cal R(\Bbb A_n)$; its elements
are of the form $R(r)$ where $r$ is a root with coefficients $0$ and $1$.
We have to determine the additional rich antichains (it is obvious that such an antichain
has to contain the vertex $[[2]]$). 

We consider the roots of $\Bbb B_n$ of the form
$$
 r = r(i,j) = (2,\dots,2,1,\dots,1,0,\dots,0)
$$ 
with $i$ coefficients being $2$, then $j-i$ coefficients being $1$ 
and finally $n-j$ coefficients
equal to $0$, where $0\le i < j \le n$. The roots of the form $r(i,j)$ with 
$1 \le i < j \le n$
are just the elements of $\Phi(\Bbb B_n,2)$. 

For $0\le i < j \le n$, we may also consider the subset
$$
 R(i,j) = \{[[2]],\dots,[[i+1]],[i+2],\dots,[j],j+1,\dots, n\}.
$$
of $\Phi_+(\Bbb B_n)$ (with $i$ roots of length 2, $j-i-1$ roots of length 1, and 
$n-j$ roots of length 1).
It is obvious that $R(i,j)$ is a rich antichain. 

We define $R = R(\Bbb B_n,2)\:\Phi(\Bbb B_n,2) \to \Cal R(\Bbb B_n)$ by $R(r(i,j)) = R(i,j)$. 
	\medskip 
Here is the case $n=5$. On the left, we draw the set of indices $(i,j)$. The map
$R$ provides the following bijection:
$$
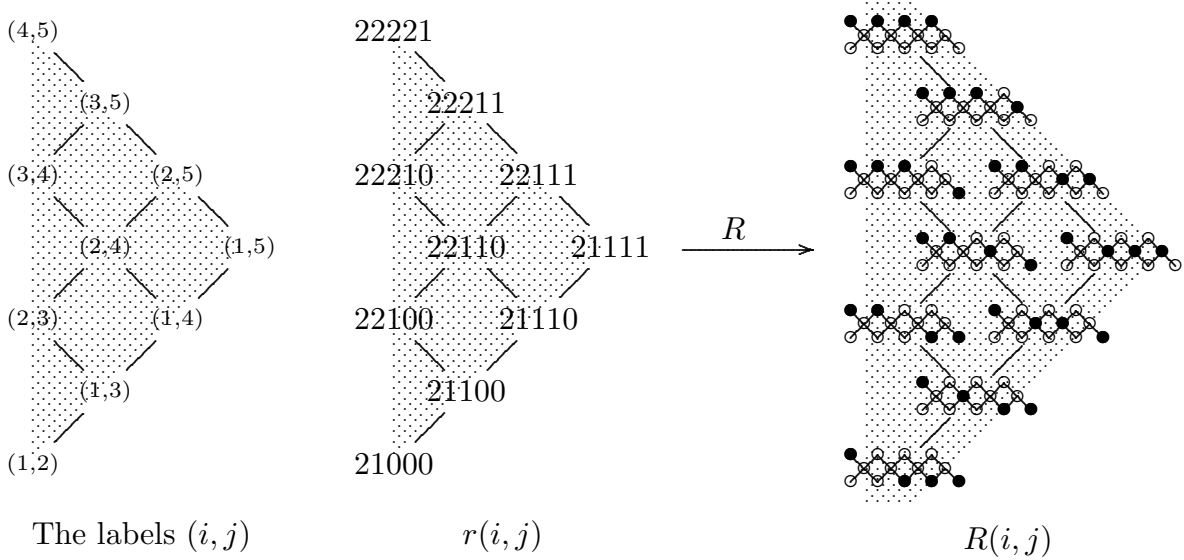
{\beginpicture
\setcoordinatesystem units <.8cm,.8cm>
\put{\beginpicture
\put{$\ssize (1,2)$} at 0 0 
\put{$\ssize (1,3)$} at 1 1 
\put{$\ssize (1,4)$} at 2 2 
\put{$\ssize (1,5)$} at 3 3 
\put{$\ssize (2,3)$} at 0 2 
\put{$\ssize (2,4)$} at 1 3 
\put{$\ssize (2,5)$} at 2 4 
\put{$\ssize (3,4)$} at 0 4  
\put{$\ssize (3,5)$} at 1 5
\put{$\ssize (4,5)$} at 0 6
\plot 0.3 0.3  0.7 0.7 /
\plot 0.3 2.3  0.7 2.7 /
\plot 0.3 4.3  0.7 4.7 /
\plot 1.3 1.3  1.7 1.7 /
\plot 1.3 3.3  1.7 3.7 /
\plot 2.3 2.3  2.7 2.7 /

\plot 0.3 1.7  0.7 1.3 /
\plot 0.3 3.7  0.7 3.3 /
\plot 0.3 5.7  0.7 5.3 /
\plot 1.3 2.7  1.7 2.3 /
\plot 1.3 4.7  1.7 4.3 /
\plot 2.3 3.7  2.7 3.3 /

\setshadegrid span <.6mm>
\vshade 0 0 6 <z,z,,> 3 3 3   /

\put{The labels $(i,j)$} at 1.5 -1
\endpicture} at 0 0
\put{\beginpicture
\put{$21000$} at 0 0 
\put{$21100$} at 1 1 
\put{$21110$} at 2 2 
\put{$21111$} at 3 3 
\put{$22100$} at 0 2 
\put{$22110$} at 1 3 
\put{$22111$} at 2 4 
\put{$22210$} at 0 4  
\put{$22211$} at 1 5
\put{$22221$} at 0 6
\plot 0.3 0.3  0.7 0.7 /
\plot 0.3 2.3  0.7 2.7 /
\plot 0.3 4.3  0.7 4.7 /
\plot 1.3 1.3  1.7 1.7 /
\plot 1.3 3.3  1.7 3.7 /
\plot 2.3 2.3  2.7 2.7 /

\plot 0.3 1.7  0.7 1.3 /
\plot 0.3 3.7  0.7 3.3 /
\plot 0.3 5.7  0.7 5.3 /
\plot 1.3 2.7  1.7 2.3 /
\plot 1.3 4.7  1.7 4.3 /
\plot 2.3 3.7  2.7 3.3 /
\put{$r(i,j)$} at 1.5 -1
\setshadegrid span <.6mm>
\vshade 0 0 6 <z,z,,> 3 3 3   /
\endpicture} at 5 0

\put{\beginpicture
\put{
 \beginpicture
 \setcoordinatesystem units <.15cm,.15cm>
 \multiput{$\circ$} at 
   0 0  2 0  4 0  6 0  8 0  1 1  3 1  5 1  7 1  0 2  2 2  4 2  6 2  /
 \plot 0 0  1 1  2 0  3 1  4 0  5 1  6 0  7 1  8 0 /
 \plot 0 2  1 1  2 2  3 1  4 2  5 1  6 2  7 1 /
\multiput{$\bullet$} at 0 2  4 0  6 0  8 0 /
 \endpicture
 } at 0 0 
\put{
 \beginpicture
 \setcoordinatesystem units <.15cm,.15cm>
 \multiput{$\circ$} at 
   0 0  2 0  4 0  6 0  8 0  1 1  3 1  5 1  7 1  0 2  2 2  4 2  6 2  /
 \plot 0 0  1 1  2 0  3 1  4 0  5 1  6 0  7 1  8 0 /
 \plot 0 2  1 1  2 2  3 1  4 2  5 1  6 2  7 1 /
\multiput{$\bullet$} at 0 2  3 1  6 0  8 0 /
 \endpicture
} at 1 1 
\put{ 
\beginpicture
 \setcoordinatesystem units <.15cm,.15cm>
 \multiput{$\circ$} at 
   0 0  2 0  4 0  6 0  8 0  1 1  3 1  5 1  7 1  0 2  2 2  4 2  6 2  /
 \plot 0 0  1 1  2 0  3 1  4 0  5 1  6 0  7 1  8 0 /
 \plot 0 2  1 1  2 2  3 1  4 2  5 1  6 2  7 1 /
\multiput{$\bullet$} at 0 2  3 1  5 1  8 0 /
 \endpicture
} at 2 2 
\put{
 \beginpicture
 \setcoordinatesystem units <.15cm,.15cm>
 \multiput{$\circ$} at 
   0 0  2 0  4 0  6 0  8 0  1 1  3 1  5 1  7 1  0 2  2 2  4 2  6 2  /
 \plot 0 0  1 1  2 0  3 1  4 0  5 1  6 0  7 1  8 0 /
 \plot 0 2  1 1  2 2  3 1  4 2  5 1  6 2  7 1 /
\multiput{$\bullet$} at 0 2  3 1  5 1  7 1  /
 \endpicture
} at 3 3 
\put{
 \beginpicture
 \setcoordinatesystem units <.15cm,.15cm>
 \multiput{$\circ$} at 
   0 0  2 0  4 0  6 0  8 0  1 1  3 1  5 1  7 1  0 2  2 2  4 2  6 2  /
 \plot 0 0  1 1  2 0  3 1  4 0  5 1  6 0  7 1  8 0 /
 \plot 0 2  1 1  2 2  3 1  4 2  5 1  6 2  7 1 /
\multiput{$\bullet$} at 0 2  2 2  6 0  8 0 /
 \endpicture
} at 0 2 
\put{
 \beginpicture
 \setcoordinatesystem units <.15cm,.15cm>
 \multiput{$\circ$} at 
   0 0  2 0  4 0  6 0  8 0  1 1  3 1  5 1  7 1  0 2  2 2  4 2  6 2  /
 \plot 0 0  1 1  2 0  3 1  4 0  5 1  6 0  7 1  8 0 /
 \plot 0 2  1 1  2 2  3 1  4 2  5 1  6 2  7 1 /
\multiput{$\bullet$} at 0 2  2 2  5 1  8 0 /
 \endpicture
} at 1 3 
\put{
 \beginpicture
 \setcoordinatesystem units <.15cm,.15cm>
 \multiput{$\circ$} at 
   0 0  2 0  4 0  6 0  8 0  1 1  3 1  5 1  7 1  0 2  2 2  4 2  6 2  /
 \plot 0 0  1 1  2 0  3 1  4 0  5 1  6 0  7 1  8 0 /
 \plot 0 2  1 1  2 2  3 1  4 2  5 1  6 2  7 1 /
\multiput{$\bullet$} at 0 2  2 2  5 1  7 1  /
 \endpicture
} at 2 4 
\put{
 \beginpicture
 \setcoordinatesystem units <.15cm,.15cm>
 \multiput{$\circ$} at 
   0 0  2 0  4 0  6 0  8 0  1 1  3 1  5 1  7 1  0 2  2 2  4 2  6 2  /
 \plot 0 0  1 1  2 0  3 1  4 0  5 1  6 0  7 1  8 0 /
 \plot 0 2  1 1  2 2  3 1  4 2  5 1  6 2  7 1 /
\multiput{$\bullet$} at 0 2  2 2  4 2  8 0 /
 \endpicture
} at 0 4  
\put{
 \beginpicture
 \setcoordinatesystem units <.15cm,.15cm>
 \multiput{$\circ$} at 
   0 0  2 0  4 0  6 0  8 0  1 1  3 1  5 1  7 1  0 2  2 2  4 2  6 2  /
 \plot 0 0  1 1  2 0  3 1  4 0  5 1  6 0  7 1  8 0 /
 \plot 0 2  1 1  2 2  3 1  4 2  5 1  6 2  7 1 /
\multiput{$\bullet$} at 0 2  2 2  4 2  7 1 /
 \endpicture
} at 1 5
\put{
 \beginpicture
 \setcoordinatesystem units <.15cm,.15cm>
 \multiput{$\circ$} at 
   0 0  2 0  4 0  6 0  8 0  1 1  3 1  5 1  7 1  0 2  2 2  4 2  6 2  /
 \plot 0 0  1 1  2 0  3 1  4 0  5 1  6 0  7 1  8 0 /
 \plot 0 2  1 1  2 2  3 1  4 2  5 1  6 2  7 1 /
\multiput{$\bullet$} at 0 2  2 2  4 2  6 2  /
 \endpicture
} at 0 6
\plot 0.3 0.3  0.7 0.7 /
\plot 0.3 2.3  0.7 2.7 /
\plot 0.3 4.3  0.7 4.7 /
\plot 1.3 1.3  1.7 1.7 /
\plot 1.3 3.3  1.7 3.7 /
\plot 2.3 2.3  2.7 2.7 /

\plot 0.3 1.7  0.7 1.3 /
\plot 0.3 3.7  0.7 3.3 /
\plot 0.3 5.7  0.7 5.3 /
\plot 1.3 2.7  1.7 2.3 /
\plot 1.3 4.7  1.7 4.3 /
\plot 2.3 3.7  2.7 3.3 /
\put{$R(i,j)$} at 1.5 -1

\setshadegrid span <.6mm>
\vshade -.5 -.5 6.5 <z,z,,> 0.2 -.5 6.5 <z,z,,> 3.7 3 3   /

\endpicture} at 12 0
\arr{7.5 0.5}{9.3 0.5}
\put{$R$} at 8.2 0.8
\endpicture}
$$
One checks without difficulties that the map $R$ preserves and
reflects the respective partial orderings. 

It remains to see in which way the partial ordered sets 
$\Phi(\Bbb B_n,1)$ and $\Phi(B_n,2)$ as well as 
the images of $R(\Bbb B_n,1)\:\Phi(\Bbb B_n,1) \to \Cal R(\Bbb B_n)$ and 
$R(\Bbb B_n,2)\:\Phi(\Bbb B_n,2) \to \Cal R(\Bbb B_n)$ are connected.
This concerns the following neighbors, drawn by dotted lines:
$$
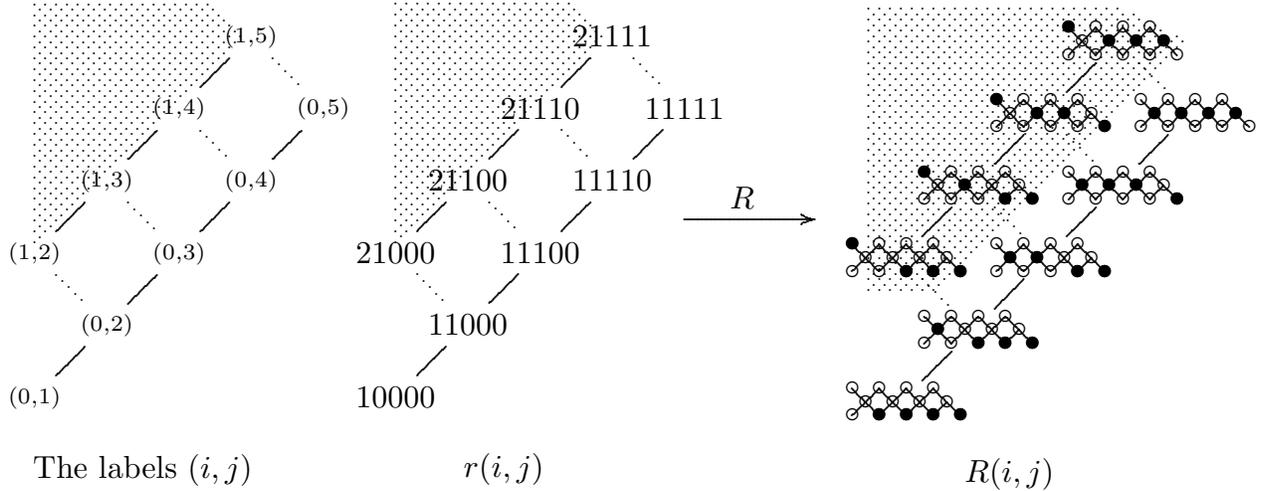
{\beginpicture
\setcoordinatesystem units <.8cm,.8cm>
\put{\beginpicture
\put{$\ssize (1,2)$} at 0 0 
\put{$\ssize (1,3)$} at 1 1 
\put{$\ssize (1,4)$} at 2 2 
\put{$\ssize (1,5)$} at 3 3 

\put{$\ssize (0,1)$} at 0 -2 
\put{$\ssize (0,2)$} at 1 -1 
\put{$\ssize (0,3)$} at 2 0
\put{$\ssize (0,4)$} at 3 1  
\put{$\ssize (0,5)$} at 4 2
\plot 0.3 0.3  0.7 0.7 /
\plot 1.3 1.3  1.7 1.7 /
\plot 2.3 2.3  2.7 2.7 /

\plot 0.3 -1.7	 0.7 -1.3 /
\plot 1.3 -.7  1.7 -.3 /
\plot 2.3 0.3  2.7 0.7 /
\plot 3.3 1.3  3.7 1.7 /
\setdots <1mm>
\plot 0.3 -0.3  0.7 -0.7 /
\plot 1.3  0.7  1.7  0.3 /
\plot 2.3  1.7  2.7  1.3 /
\plot 3.3  2.7  3.7  2.3 /

\put{The labels $(i,j)$} at 1.5 -3

\setshadegrid span <.6mm>
\vshade 0 0 3.5 <z,z,,> 2.5 2.5 3.5 <z,z,,> 3 3 3   /
\endpicture} at 0 0
\put{\beginpicture

\put{$21000$} at 0 0 
\put{$21100$} at 1 1 
\put{$21110$} at 2 2 
\put{$21111$} at 3 3 

\put{$10000$} at 0 -2 
\put{$11000$} at 1 -1 
\put{$11100$} at 2 0
\put{$11110$} at 3 1  
\put{$11111$} at 4 2
\plot 0.3 0.3  0.7 0.7 /
\plot 1.3 1.3  1.7 1.7 /
\plot 2.3 2.3  2.7 2.7 /

\plot 0.3 -1.7	 0.7 -1.3 /
\plot 1.3 -.7  1.7 -.3 /
\plot 2.3 0.3  2.7 0.7 /
\plot 3.3 1.3  3.7 1.7 /
\setdots <1mm>
\plot 0.3 -0.3  0.7 -0.7 /
\plot 1.3  0.7  1.7  0.3 /
\plot 2.3  1.7  2.7  1.3 /
\plot 3.3  2.7  3.7  2.3 /

\put{$r(i,j)$} at 1.5 -3
\setshadegrid span <.6mm>
\vshade 0 0 3.5 <z,z,,> 2.5 2.5 3.5 <z,z,,> 3 3 3   /
\endpicture} at 5 0

\put{\beginpicture
\put{
 \beginpicture
 \setcoordinatesystem units <.15cm,.15cm>
 \multiput{$\circ$} at 
   0 0  2 0  4 0  6 0  8 0  1 1  3 1  5 1  7 1  0 2  2 2  4 2  6 2  /
 \plot 0 0  1 1  2 0  3 1  4 0  5 1  6 0  7 1  8 0 /
 \plot 0 2  1 1  2 2  3 1  4 2  5 1  6 2  7 1 /
\multiput{$\bullet$} at 0 2  4 0  6 0  8 0 /
 \endpicture
 } at 0 0 
\put{
 \beginpicture
 \setcoordinatesystem units <.15cm,.15cm>
 \multiput{$\circ$} at 
   0 0  2 0  4 0  6 0  8 0  1 1  3 1  5 1  7 1  0 2  2 2  4 2  6 2  /
 \plot 0 0  1 1  2 0  3 1  4 0  5 1  6 0  7 1  8 0 /
 \plot 0 2  1 1  2 2  3 1  4 2  5 1  6 2  7 1 /
\multiput{$\bullet$} at 0 2  3 1  6 0  8 0 /
 \endpicture
} at 1 1 
\put{ 
\beginpicture
 \setcoordinatesystem units <.15cm,.15cm>
 \multiput{$\circ$} at 
   0 0  2 0  4 0  6 0  8 0  1 1  3 1  5 1  7 1  0 2  2 2  4 2  6 2  /
 \plot 0 0  1 1  2 0  3 1  4 0  5 1  6 0  7 1  8 0 /
 \plot 0 2  1 1  2 2  3 1  4 2  5 1  6 2  7 1 /
\multiput{$\bullet$} at 0 2  3 1  5 1  8 0 /

 \endpicture
} at 2 2 
\put{
 \beginpicture
 \setcoordinatesystem units <.15cm,.15cm>
 \multiput{$\circ$} at 
   0 0  2 0  4 0  6 0  8 0  1 1  3 1  5 1  7 1  0 2  2 2  4 2  6 2  /
 \plot 0 0  1 1  2 0  3 1  4 0  5 1  6 0  7 1  8 0 /
 \plot 0 2  1 1  2 2  3 1  4 2  5 1  6 2  7 1 /
\multiput{$\bullet$} at 0 2  3 1  5 1  7 1  /
 \endpicture
} at 3 3

\put{
 \beginpicture
 \setcoordinatesystem units <.15cm,.15cm>
 \multiput{$\circ$} at 
   0 0  2 0  4 0  6 0  8 0  1 1  3 1  5 1  7 1  0 2  2 2  4 2  6 2  /
 \plot 0 0  1 1  2 0  3 1  4 0  5 1  6 0  7 1  8 0 /
 \plot 0 2  1 1  2 2  3 1  4 2  5 1  6 2  7 1 /
\multiput{$\bullet$} at 2 0  4 0  6 0  8 0 /
 \endpicture
 } at 0 -2
\put{
 \beginpicture
 \setcoordinatesystem units <.15cm,.15cm>
 \multiput{$\circ$} at 
   0 0  2 0  4 0  6 0  8 0  1 1  3 1  5 1  7 1  0 2  2 2  4 2  6 2  /
 \plot 0 0  1 1  2 0  3 1  4 0  5 1  6 0  7 1  8 0 /
 \plot 0 2  1 1  2 2  3 1  4 2  5 1  6 2  7 1 /
\multiput{$\bullet$} at 1 1  4 0  6 0  8 0 /
 \endpicture
 } at 1 -1
\put{
 \beginpicture
 \setcoordinatesystem units <.15cm,.15cm>
 \multiput{$\circ$} at 
   0 0  2 0  4 0  6 0  8 0  1 1  3 1  5 1  7 1  0 2  2 2  4 2  6 2  /
 \plot 0 0  1 1  2 0  3 1  4 0  5 1  6 0  7 1  8 0 /
 \plot 0 2  1 1  2 2  3 1  4 2  5 1  6 2  7 1 /
\multiput{$\bullet$} at 1 1  3 1  6 0  8 0 /
 \endpicture
 } at 2 0
\put{
 \beginpicture
 \setcoordinatesystem units <.15cm,.15cm>
 \multiput{$\circ$} at 
   0 0  2 0  4 0  6 0  8 0  1 1  3 1  5 1  7 1  0 2  2 2  4 2  6 2  /
 \plot 0 0  1 1  2 0  3 1  4 0  5 1  6 0  7 1  8 0 /
 \plot 0 2  1 1  2 2  3 1  4 2  5 1  6 2  7 1 /
\multiput{$\bullet$} at 1 1  3 1  5 1  8 0 /
 \endpicture
 } at 3 1 
\put{
 \beginpicture
 \setcoordinatesystem units <.15cm,.15cm>
 \multiput{$\circ$} at 
   0 0  2 0  4 0  6 0  8 0  1 1  3 1  5 1  7 1  0 2  2 2  4 2  6 2  /
 \plot 0 0  1 1  2 0  3 1  4 0  5 1  6 0  7 1  8 0 /
 \plot 0 2  1 1  2 2  3 1  4 2  5 1  6 2  7 1 /
\multiput{$\bullet$} at 1 1  3 1  5 1  7 1  /
 \endpicture
 } at 4 2

\plot 0.3 0.3  0.7 0.7 /
\plot 1.3 1.3  1.7 1.7 /
\plot 2.3 2.3  2.7 2.7 /

\plot 0.3 -1.7	 0.7 -1.3 /
\plot 1.3 -.7  1.7 -.3 /
\plot 2.3 0.3  2.7 0.7 /
\plot 3.3 1.3  3.7 1.7 /
\setdots <1mm>
\plot 0.3 -0.3  0.7 -0.7 /
\plot 1.3  0.7  1.7  0.3 /
\plot 2.3  1.7  2.7  1.3 /
\plot 3.3  2.7  3.7  2.3 /

\put{$R(i,j)$} at 1.5 -3
\setshadegrid span <.6mm>
\vshade -.5 -.5 3.5 <z,z,,> 0.5 -.5 3.5 <z,z,,> 2.5 1.5 3.5 <z,z,,> 
   3 2 3.5 <z,z,,>  4 3 3   /

\endpicture} at 12 0

\arr{7 0.5}{8.8 0.5}
\put{$R$} at 7.8 0.8

\endpicture}
$$
It turns out that the bijections $R(\Bbb B_n,1)$ and $R(\Bbb B_n,2)$
combine to provide an isomorphism $\Phi_+(\Bbb B_n) \to \Cal R(\Bbb B_n).$
   \bigskip
{\bf Case} $\Bbb D_n.$ This case is similar  to $\Bbb B_n$.
The set $\Phi(\Bbb D_n,2)$ is just the coideal
of $\Phi_+$ generated by the maximal root of $\Bbb D_4$.
For $0\le i < j \le n-2$, let $r = r(i,j)$ be the root 
of $\Bbb D_n$ of the form
$$
 r = r(i,j) = \left(\smallmatrix 1 \cr
                            & 2 & \cdots & 2 & 1 \cdots & 1 & 0 \cdots & 0 \cr
                          1 \endsmallmatrix \right)
$$ 
with $i$ coefficients being $2$ and $j-i+2$ coefficients being 1 (to be precise:
with coefficients $1$ on the short arms, and $n-j-2$ coefficients
equal to $0$ on the long arm).   
The roots of the form $r(i,j)$ with $1 \le i < j \le n-2$
are just the elements of $\Phi(\Bbb D_n,2)$. 

For $0\le i < j \le n-2$, we also consider a corresponding subset $R(i,j)$ of 
$\Phi_+(\Bbb D_n)$. For $i = 0$, $R(i,j)$ is the rich antichain 
which contain the roots 
$\left(\smallmatrix 1 \cr
                & 1 & 0 & \cdots & 0 \cr
              0 \endsmallmatrix\right)$ and 
$\left(\smallmatrix 0 \cr
                & 1 & 0 & \cdots & 0 \cr
              1 \endsmallmatrix\right)$, $j-1$ additional roots of height 2
and $n-j-2$ simple roots. For $i > 0$, $R(i,j)$ contains the $\Phi(\Bbb D_4)_3$,
$j-i-1$ roots of height 2 and $n-j-2$ simple roots. 
We define $R\:\Phi(\Bbb D_n,2) \to \Cal R(\Bbb D_n)$ by $R(r(i,j)) = R(i,j)$. 
	\medskip 
Here is the case $n=7$. On the left, we draw the set of indices $(i,j)$. The map
$R$ provides the following bijection:

$$
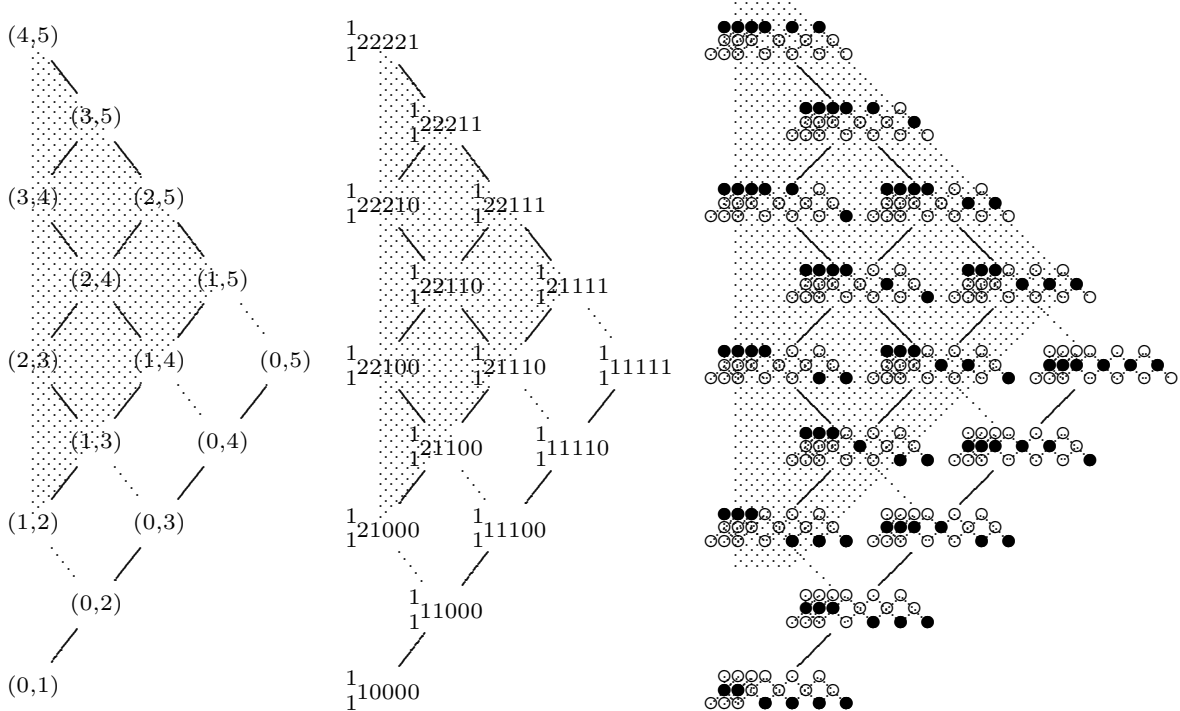
{\beginpicture
\setcoordinatesystem units <.7cm,.9cm>
\put{\beginpicture
\setcoordinatesystem units <.7cm,.9cm>
\put{$\ssize (1,2)$} at 0 0 
\put{$\ssize (1,3)$} at 1 1 
\put{$\ssize (1,4)$} at 2 2 
\put{$\ssize (1,5)$} at 3 3 
\put{$\ssize (2,3)$} at 0 2 
\put{$\ssize (2,4)$} at 1 3 
\put{$\ssize (2,5)$} at 2 4 
\put{$\ssize (3,4)$} at 0 4  
\put{$\ssize (3,5)$} at 1 5
\put{$\ssize (4,5)$} at 0 6
\plot 0.3 0.3  0.7 0.7 /
\plot 0.3 2.3  0.7 2.7 /
\plot 0.3 4.3  0.7 4.7 /
\plot 1.3 1.3  1.7 1.7 /
\plot 1.3 3.3  1.7 3.7 /
\plot 2.3 2.3  2.7 2.7 /

\plot 0.3 1.7  0.7 1.3 /
\plot 0.3 3.7  0.7 3.3 /
\plot 0.3 5.7  0.7 5.3 /
\plot 1.3 2.7  1.7 2.3 /
\plot 1.3 4.7  1.7 4.3 /
\plot 2.3 3.7  2.7 3.3 /

\setshadegrid span <.6mm>
\vshade 0 0 6 <z,z,,> 3 3 3   /

\put{$\ssize (0,1)$} at 0 -2 
\put{$\ssize (0,2)$} at 1 -1 
\put{$\ssize (0,3)$} at 2 0
\put{$\ssize (0,4)$} at 3 1  
\put{$\ssize (0,5)$} at 4 2

\plot 0.3 -1.7	 0.7 -1.3 /
\plot 1.3 -.7  1.7 -.3 /
\plot 2.3 0.3  2.7 0.7 /
\plot 3.3 1.3  3.7 1.7 /
\setdots <1mm>
\plot 0.3 -0.3  0.7 -0.7 /
\plot 1.3  0.7  1.7  0.3 /
\plot 2.3  1.7  2.7  1.3 /
\plot 3.3  2.7  3.7  2.3 /

\put{The labels $(i,j)$} at 1.5 -3.5

\endpicture} at -.5 0
\put{\beginpicture
\setcoordinatesystem units <.7cm,.9cm>
\put{${1\atop1}\!\ssize 21000$} at 0 0 
\put{${1\atop1}\!\ssize 21100$} at 1 1 
\put{${1\atop1}\!\ssize 21110$} at 2 2 
\put{${1\atop1}\!\ssize 21111$} at 3 3 
\put{${1\atop1}\!\ssize 22100$} at 0 2 
\put{${1\atop1}\!\ssize 22110$} at 1 3 
\put{${1\atop1}\!\ssize 22111$} at 2 4 
\put{${1\atop1}\!\ssize 22210$} at 0 4  
\put{${1\atop1}\!\ssize 22211$} at 1 5
\put{${1\atop1}\!\ssize 22221$} at 0 6
\plot 0.3 0.3  0.7 0.7 /
\plot 0.3 2.3  0.7 2.7 /
\plot 0.3 4.3  0.7 4.7 /
\plot 1.3 1.3  1.7 1.7 /
\plot 1.3 3.3  1.7 3.7 /
\plot 2.3 2.3  2.7 2.7 /

\plot 0.3 1.7  0.7 1.3 /
\plot 0.3 3.7  0.7 3.3 /
\plot 0.3 5.7  0.7 5.3 /
\plot 1.3 2.7  1.7 2.3 /
\plot 1.3 4.7  1.7 4.3 /
\plot 2.3 3.7  2.7 3.3 /
\setshadegrid span <.6mm>
\vshade 0 0 6 <z,z,,> 3 3 3   /

\put{${1\atop1}\!\ssize 10000$} at 0 -2 
\put{${1\atop1}\!\ssize 11000$} at 1 -1 
\put{${1\atop1}\!\ssize 11100$} at 2 0
\put{${1\atop1}\!\ssize 11110$} at 3 1  
\put{${1\atop1}\!\ssize 11111$} at 4 2

\plot 0.3 -1.7	 0.7 -1.3 /
\plot 1.3 -.7  1.7 -.3 /
\plot 2.3 0.3  2.7 0.7 /
\plot 3.3 1.3  3.7 1.7 /
\setdots <1mm>
\plot 0.3 -0.3  0.7 -0.7 /
\plot 1.3  0.7  1.7  0.3 /
\plot 2.3  1.7  2.7  1.3 /
\plot 3.3  2.7  3.7  2.3 /

\put{$r(i,j)$} at 1.5 -3.5

\endpicture} at 5 0

\put{\beginpicture
\setcoordinatesystem units <.9cm,.9cm>
\put{
 \beginpicture
 \setcoordinatesystem units <.15cm,.15cm>
 \multiput{$\circ$} at 
   0 0  2 0  4 0  6 0  8 0  1 1  3 1  5 1  7 1  0 2  2 2  4 2  6 2  /
\setdots <.5mm>
 \plot 0 0  1 1  2 0  3 1  4 0  5 1  6 0  7 1  8 0 /
 \plot 0 2  1 1  2 2  3 1  4 2  5 1  6 2  7 1 /
\multiput{$\bullet$} at 0 2  4 0  6 0  8 0 /

\multiput{$\bullet$} at 0 2  4 0  6 0  8 0 /
\multiput{$\circ$} at  -2 0  -1 1 -1 0  -1 2   0 1  1 2 /
\multiput{$\bullet$} at -1 2  1 2 /
\plot -2 0  0 2 /
\plot -1 1  -1 2  0 1  1 2  1 1 /
\plot -1 0  0 1  0 0  -1 1  /

 \endpicture
 } at 0 0 
\put{
 \beginpicture
 \setcoordinatesystem units <.15cm,.15cm>
 \multiput{$\circ$} at 
   0 0  2 0  4 0  6 0  8 0  1 1  3 1  5 1  7 1  0 2  2 2  4 2  6 2  /
\setdots <.5mm>
 \plot 0 0  1 1  2 0  3 1  4 0  5 1  6 0  7 1  8 0 /
 \plot 0 2  1 1  2 2  3 1  4 2  5 1  6 2  7 1 /
\multiput{$\bullet$} at 0 2  3 1  6 0  8 0 /
\multiput{$\circ$} at  -2 0  -1 1 -1 0  -1 2   0 1  1 2 /
\multiput{$\bullet$} at -1 2  1 2 /
\plot -2 0  0 2 /
\plot -1 1  -1 2  0 1  1 2  1 1 /
\plot -1 0  0 1  0 0  -1 1  /

 \endpicture
} at 1 1 
\put{ 
\beginpicture
 \setcoordinatesystem units <.15cm,.15cm>
 \multiput{$\circ$} at 
   0 0  2 0  4 0  6 0  8 0  1 1  3 1  5 1  7 1  0 2  2 2  4 2  6 2  /
\setdots <.5mm>
 \plot 0 0  1 1  2 0  3 1  4 0  5 1  6 0  7 1  8 0 /
 \plot 0 2  1 1  2 2  3 1  4 2  5 1  6 2  7 1 /
\multiput{$\bullet$} at 0 2  3 1  5 1  8 0 /
\multiput{$\circ$} at  -2 0  -1 1 -1 0  -1 2   0 1  1 2 /
\multiput{$\bullet$} at -1 2  1 2 /
\plot -2 0  0 2 /
\plot -1 1  -1 2  0 1  1 2  1 1 /
\plot -1 0  0 1  0 0  -1 1  /

 \endpicture
} at 2 2 
\put{
 \beginpicture
 \setcoordinatesystem units <.15cm,.15cm>
 \multiput{$\circ$} at 
   0 0  2 0  4 0  6 0  8 0  1 1  3 1  5 1  7 1  0 2  2 2  4 2  6 2  /
\setdots <.5mm>
 \plot 0 0  1 1  2 0  3 1  4 0  5 1  6 0  7 1  8 0 /
 \plot 0 2  1 1  2 2  3 1  4 2  5 1  6 2  7 1 /
\multiput{$\bullet$} at 0 2  3 1  5 1  7 1  /
\multiput{$\circ$} at  -2 0  -1 1 -1 0  -1 2   0 1  1 2 /
\multiput{$\bullet$} at -1 2  1 2 /
\plot -2 0  0 2 /
\plot -1 1  -1 2  0 1  1 2  1 1 /
\plot -1 0  0 1  0 0  -1 1  /

 \endpicture
} at 3 3 
\put{
 \beginpicture
 \setcoordinatesystem units <.15cm,.15cm>
 \multiput{$\circ$} at 
   0 0  2 0  4 0  6 0  8 0  1 1  3 1  5 1  7 1  0 2  2 2  4 2  6 2  /
\setdots <.5mm>
 \plot 0 0  1 1  2 0  3 1  4 0  5 1  6 0  7 1  8 0 /
 \plot 0 2  1 1  2 2  3 1  4 2  5 1  6 2  7 1 /
\multiput{$\bullet$} at 0 2  2 2  6 0  8 0 /
\multiput{$\circ$} at  -2 0  -1 1 -1 0  -1 2   0 1  1 2 /
\multiput{$\bullet$} at -1 2  1 2 /
\plot -2 0  0 2 /
\plot -1 1  -1 2  0 1  1 2  1 1 /
\plot -1 0  0 1  0 0  -1 1  /

 \endpicture
} at 0 2 
\put{
 \beginpicture
 \setcoordinatesystem units <.15cm,.15cm>
 \multiput{$\circ$} at 
   0 0  2 0  4 0  6 0  8 0  1 1  3 1  5 1  7 1  0 2  2 2  4 2  6 2  /
\setdots <.5mm>
 \plot 0 0  1 1  2 0  3 1  4 0  5 1  6 0  7 1  8 0 /
 \plot 0 2  1 1  2 2  3 1  4 2  5 1  6 2  7 1 /
\multiput{$\bullet$} at 0 2  2 2  5 1  8 0 /
\multiput{$\circ$} at  -2 0  -1 1 -1 0  -1 2   0 1  1 2 /
\multiput{$\bullet$} at -1 2  1 2 /
\plot -2 0  0 2 /
\plot -1 1  -1 2  0 1  1 2  1 1 /
\plot -1 0  0 1  0 0  -1 1  /

 \endpicture
} at 1 3 
\put{
 \beginpicture
 \setcoordinatesystem units <.15cm,.15cm>
 \multiput{$\circ$} at 
   0 0  2 0  4 0  6 0  8 0  1 1  3 1  5 1  7 1  0 2  2 2  4 2  6 2  /
\setdots <.5mm>
 \plot 0 0  1 1  2 0  3 1  4 0  5 1  6 0  7 1  8 0 /
 \plot 0 2  1 1  2 2  3 1  4 2  5 1  6 2  7 1 /
\multiput{$\bullet$} at 0 2  2 2  5 1  7 1  /
\multiput{$\circ$} at  -2 0  -1 1 -1 0  -1 2   0 1  1 2 /
\multiput{$\bullet$} at -1 2  1 2 /
\plot -2 0  0 2 /
\plot -1 1  -1 2  0 1  1 2  1 1 /
\plot -1 0  0 1  0 0  -1 1  /

 \endpicture
} at 2 4 
\put{
 \beginpicture
 \setcoordinatesystem units <.15cm,.15cm>
 \multiput{$\circ$} at 
   0 0  2 0  4 0  6 0  8 0  1 1  3 1  5 1  7 1  0 2  2 2  4 2  6 2  /
\setdots <.5mm>
 \plot 0 0  1 1  2 0  3 1  4 0  5 1  6 0  7 1  8 0 /
 \plot 0 2  1 1  2 2  3 1  4 2  5 1  6 2  7 1 /
\multiput{$\bullet$} at 0 2  2 2  4 2  8 0 /
\multiput{$\circ$} at  -2 0  -1 1 -1 0  -1 2   0 1  1 2 /
\multiput{$\bullet$} at -1 2  1 2 /
\plot -2 0  0 2 /
\plot -1 1  -1 2  0 1  1 2  1 1 /
\plot -1 0  0 1  0 0  -1 1  /

 \endpicture
} at 0 4  
\put{
 \beginpicture
 \setcoordinatesystem units <.15cm,.15cm>
 \multiput{$\circ$} at 
   0 0  2 0  4 0  6 0  8 0  1 1  3 1  5 1  7 1  0 2  2 2  4 2  6 2  /
\setdots <.5mm>
 \plot 0 0  1 1  2 0  3 1  4 0  5 1  6 0  7 1  8 0 /
 \plot 0 2  1 1  2 2  3 1  4 2  5 1  6 2  7 1 /
\multiput{$\bullet$} at 0 2  2 2  4 2  7 1 /
\multiput{$\circ$} at  -2 0  -1 1 -1 0  -1 2   0 1  1 2 /
\multiput{$\bullet$} at -1 2  1 2 /
\plot -2 0  0 2 /
\plot -1 1  -1 2  0 1  1 2  1 1 /
\plot -1 0  0 1  0 0  -1 1  /

 \endpicture
} at 1 5
\put{
 \beginpicture
 \setcoordinatesystem units <.15cm,.15cm>
 \multiput{$\circ$} at 
   0 0  2 0  4 0  6 0  8 0  1 1  3 1  5 1  7 1  0 2  2 2  4 2  6 2  /
\setdots <.5mm>
 \plot 0 0  1 1  2 0  3 1  4 0  5 1  6 0  7 1  8 0 /
 \plot 0 2  1 1  2 2  3 1  4 2  5 1  6 2  7 1 /
\multiput{$\bullet$} at 0 2  2 2  4 2  6 2  /
\multiput{$\circ$} at  -2 0  -1 1 -1 0  -1 2   0 1  1 2 /
\multiput{$\bullet$} at -1 2  1 2 /
\plot -2 0  0 2 /
\plot -1 1  -1 2  0 1  1 2  1 1 /
\plot -1 0  0 1  0 0  -1 1  /

 \endpicture
} at 0 6
\plot 0.3 0.3  0.7 0.7 /
\plot 0.3 2.3  0.7 2.7 /
\plot 0.3 4.3  0.7 4.7 /
\plot 1.3 1.3  1.7 1.7 /
\plot 1.3 3.3  1.7 3.7 /
\plot 2.3 2.3  2.7 2.7 /

\plot 0.3 1.7  0.7 1.3 /
\plot 0.3 3.7  0.7 3.3 /
\plot 0.3 5.7  0.7 5.3 /
\plot 1.3 2.7  1.7 2.3 /
\plot 1.3 4.7  1.7 4.3 /
\plot 2.3 3.7  2.7 3.3 /

\put{
 \beginpicture
 \setcoordinatesystem units <.15cm,.15cm>
 \multiput{$\circ$} at 
   0 0  2 0  4 0  6 0  8 0  1 1  3 1  5 1  7 1  0 2  2 2  4 2  6 2  /
\setdots <.5mm>
 \plot 0 0  1 1  2 0  3 1  4 0  5 1  6 0  7 1  8 0 /
 \plot 0 2  1 1  2 2  3 1  4 2  5 1  6 2  7 1 /
\multiput{$\bullet$} at 2 0  4 0  6 0  8 0 /

\multiput{$\bullet$} at 0 1  2 0  4 0  6 0  8 0 /
\multiput{$\circ$} at  -2 0  -1 1 -1 0  -1 2   0 1  1 2 /
\multiput{$\bullet$} at -1 1  0 1 /
\plot -2 0  0 2 /
\plot -1 1  -1 2  0 1  1 2  1 1 /
\plot -1 0  0 1  0 0  -1 1  /

 \endpicture
 } at 0 -2
\put{
 \beginpicture
 \setcoordinatesystem units <.15cm,.15cm>
 \multiput{$\circ$} at 
   0 0  2 0  4 0  6 0  8 0  1 1  3 1  5 1  7 1  0 2  2 2  4 2  6 2  /
\setdots <.5mm>
 \plot 0 0  1 1  2 0  3 1  4 0  5 1  6 0  7 1  8 0 /
 \plot 0 2  1 1  2 2  3 1  4 2  5 1  6 2  7 1 /
\multiput{$\bullet$} at 1 1  4 0  6 0  8 0 /
\multiput{$\circ$} at  -2 0  -1 0  -1 2   0 1  1 2 /
\multiput{$\bullet$} at -1 1  0 1  1 1 /
\plot -2 0  0 2 /
\plot -1 1  -1 2  0 1  1 2  1 1 /
\plot -1 0  0 1  0 0  -1 1  /

 \endpicture
 } at 1 -1
\put{
 \beginpicture
 \setcoordinatesystem units <.15cm,.15cm>
 \multiput{$\circ$} at 
   0 0  2 0  4 0  6 0  8 0  1 1  3 1  5 1  7 1  0 2  2 2  4 2  6 2  /
\setdots <.5mm>
 \plot 0 0  1 1  2 0  3 1  4 0  5 1  6 0  7 1  8 0 /
 \plot 0 2  1 1  2 2  3 1  4 2  5 1  6 2  7 1 /
\multiput{$\bullet$} at 1 1  3 1  6 0  8 0 /

\multiput{$\circ$} at  -2 0  -1 1 -1 0  -1 2   0 1  1 2 /
\multiput{$\bullet$} at -1 1  0 1 /
\plot -2 0  0 2 /
\plot -1 1  -1 2  0 1  1 2  1 1 /
\plot -1 0  0 1  0 0  -1 1  /
 \endpicture
 } at 2 0
\put{
 \beginpicture
 \setcoordinatesystem units <.15cm,.15cm>
 \multiput{$\circ$} at 
   0 0  2 0  4 0  6 0  8 0  1 1  3 1  5 1  7 1  0 2  2 2  4 2  6 2  /
\setdots <.5mm>
 \plot 0 0  1 1  2 0  3 1  4 0  5 1  6 0  7 1  8 0 /
 \plot 0 2  1 1  2 2  3 1  4 2  5 1  6 2  7 1 /
\multiput{$\bullet$} at 1 1  3 1  5 1  8 0 /
\multiput{$\circ$} at  -2 0  -1 1 -1 0  -1 2   0 1  1 2 /
\multiput{$\bullet$} at -1 1  0 1 /
\plot -2 0  0 2 /
\plot -1 1  -1 2  0 1  1 2  1 1 /
\plot -1 0  0 1  0 0  -1 1  /
 \endpicture
 } at 3 1 
\put{
 \beginpicture
 \setcoordinatesystem units <.15cm,.15cm>
 \multiput{$\circ$} at 
   0 0  2 0  4 0  6 0  8 0  1 1  3 1  5 1  7 1  0 2  2 2  4 2  6 2  /
\setdots <.5mm>
 \plot 0 0  1 1  2 0  3 1  4 0  5 1  6 0  7 1  8 0 /
 \plot 0 2  1 1  2 2  3 1  4 2  5 1  6 2  7 1 /
\multiput{$\bullet$} at 1 1  3 1  5 1  7 1  /
\multiput{$\circ$} at  -2 0  -1 1 -1 0  -1 2   0 1  1 2 /
\multiput{$\bullet$} at -1 1  0 1 /
\plot -2 0  0 2 /
\plot -1 1  -1 2  0 1  1 2  1 1 /
\plot -1 0  0 1  0 0  -1 1  /

 \endpicture
 } at 4 2 

\plot 0.3 -1.7	 0.7 -1.3 /
\plot 1.3 -.7  1.7 -.3 /
\plot 2.3 0.3  2.7 0.7 /
\plot 3.3 1.3  3.7 1.7 /
\setdots <1mm>
\plot 0.3 -0.3  0.7 -0.7 /
\plot 1.3  0.7  1.7  0.3 /
\plot 2.3  1.7  2.7  1.3 /
\plot 3.3  2.7  3.7  2.3 /

\put{$R(i,j)$} at 1.5 -3.5
\setsolid
\plot -2 -3.4  -2 -3.6 /
\arr{-2 -3.5}{-0.5 -3.5}
\put{$R$} at -1.3 -3.2

\setshadegrid span <.6mm>
\vshade -.5 -.5 6.5 <z,z,,> 0.2 -.5 6.5 <z,z,,> 3.7 3 3   /

\endpicture} at 11.5 0

\endpicture}
$$
As in the case $\Bbb B_n$, we also see in the case $\Bbb D_n$ 
that the bijections $R(\Bbb D_n,1)$ and $R(\Bbb D_n,2)$
combine to provide an isomorphism $\Phi_+(\Bbb D_n) \to \Cal R(\Bbb D_n).$
   \bigskip
{\bf The cases $\Delta = \Bbb E_6$ and $\Delta = \Bbb E_7.$} Here, we have
to consider (besides $\Phi(\Delta,1)$) the subsets $\Phi(\Delta,2)$ and
$\Phi(\Delta,3)$, and finally also $\Phi_(\Bbb E_7,4)$. Note that the
posets $\Phi(\Bbb E_6,2)$, $\Phi(\Bbb E_7,2)$ as well
as $\Phi(\Bbb E_7,3)$ are products of the form $[3]\times [3],\
[3]\times W(4),$ and $[2]\times [4]$, respectively, where $[m]$ denotes
the chain of cardinality $m$ and where $W(4)$ is the poset 
$$
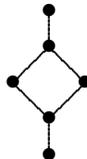
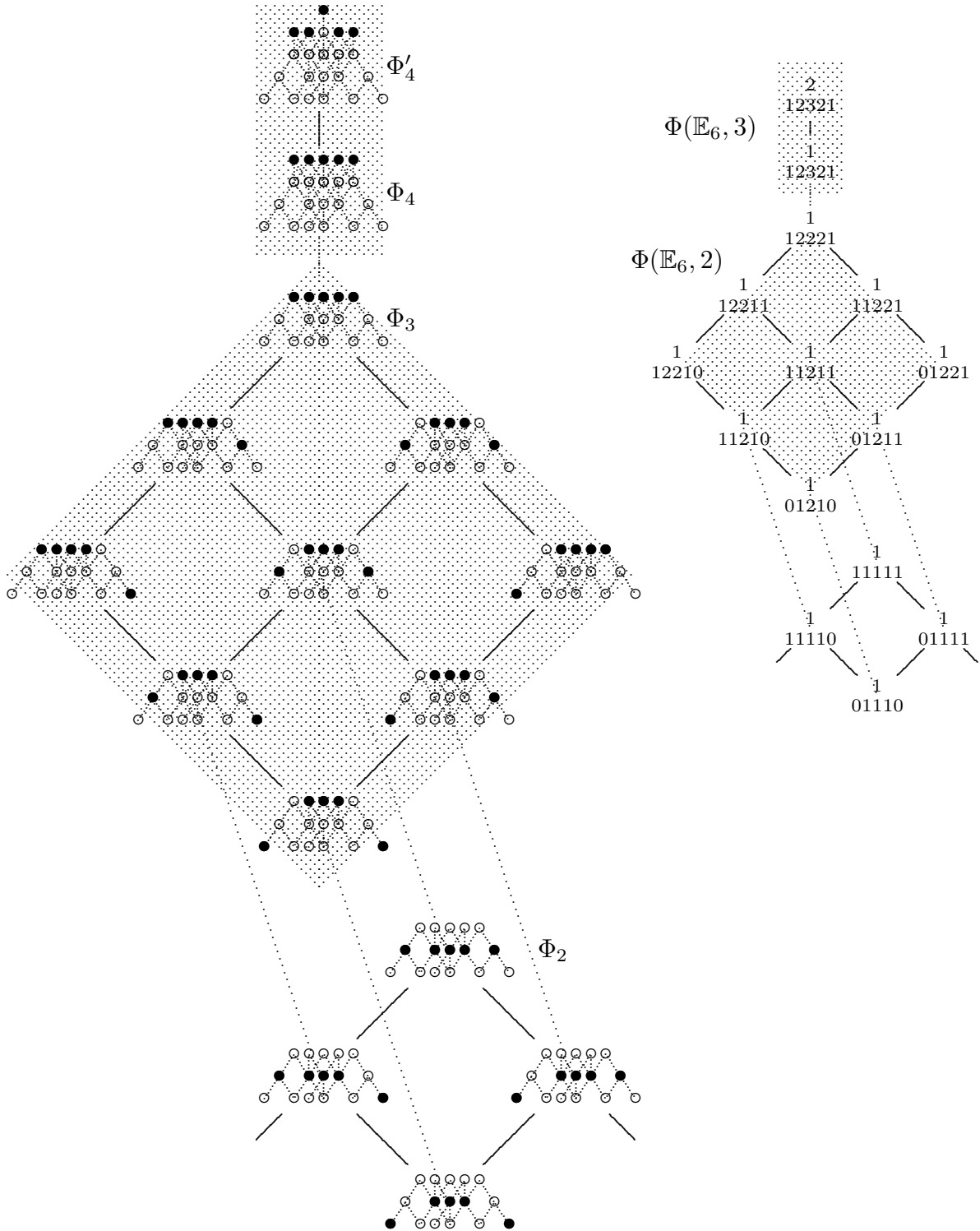
{\beginpicture
\setcoordinatesystem units <.4cm,.4cm>
\multiput{$\bullet$} at 0 0  0 1  -1 2  1 2  0 3  0 4 /
\plot 0 0  0 1  -1 2  0 3  0 4 /
\plot 0 1  1 2  0 3 /
\endpicture}
$$

\vfill\eject
{\bf The case $\Bbb E_6.$}

$$
{\beginpicture
\setcoordinatesystem units <1cm,1cm>
\put{\beginpicture
\setcoordinatesystem units <.9cm,.9cm>
\put{$1\atop01110$} at 1 -3
\put{$1\atop11110$} at 0 -2
\put{$1\atop01111$} at 2 -2
\put{$1\atop11111$} at 1 -1

\plot 0.3 -2.3  0.7 -2.7 /
\plot 1.3 -2.7  1.7 -2.3 /
\plot .3 -1.7  .7 -1.3 /
\plot 1.3 -1.3  1.7 -1.7 /

\plot -.3 -2.3 -.5 -2.5 /
\plot 2.3 -2.3 2.5 -2.5 /

\setdots <1mm>
\plot -1 1  0 -2 /
\plot 0 2  1 -1 /
\plot 0 0  1 -3 /
\plot 1 1  2 -2 /

\setsolid

\put{$1\atop01210$} at 0 0
\put{$1\atop11210$} at -1 1
\put{$1\atop01211$} at 1 1
\put{$1\atop12210$} at -2 2
\put{$1\atop11211$} at 0 2
\put{$1\atop01221$} at 2 2
\put{$1\atop12211$} at -1 3
\put{$1\atop11221$} at 1 3
\put{$1\atop12221$} at 0 4

\put{$1\atop12321$} at 0 5
\put{$2\atop12321$} at 0 6

\plot 0.3 0.3  0.7 0.7 /
\plot -0.3 0.3 -0.7 0.7 /
\plot 1.3 1.3  1.7 1.7 /
\plot -1.3 1.3 -1.7 1.7 /
\plot  0.3 1.7  0.7 1.3 /
\plot -0.3 1.7  -.7 1.3 /
\plot  1.7 2.3  1.3 2.7 /
\plot  -1.7 2.3  -1.3 2.7 /
\plot 0.3 2.3  0.7 2.7 /
\plot -.3 2.3  -.7 2.7 /
\plot 0.3 3.7  0.7 3.3 /
\plot -.3 3.7  -.7 3.3 /

\setdots <.5mm>
\plot 0 4.35  0 4.65 /
\setsolid
\plot 0 5.4  0 5.6 /

\setshadegrid span <.6mm>
\vshade -2 2 2 <z,z,,> 0 0 4  <z,z,,> 2 2 2    /
\vshade -.5 4.5 6.5 <z,z,,> .5 4.5 6.5    /

\put{$\Phi(\Bbb E_6,3)$} at -1.5 5.5
\put{$\Phi(\Bbb E_6,2)$} at -2 3.5

\endpicture} at 6.5 3

\put{\beginpicture
\setcoordinatesystem units <1.7cm,1.7cm>

\put{
\beginpicture
 \setcoordinatesystem units <.2cm,.3cm>
 \multiput{$\circ$} at
   0 0  2 0  4 0  6 0  8 0  1 1  3 1  5 1  7 1   2 2  4 2  6 2
      3 2  5 2  4 1  3 0 /
       \setdots <.5mm>
        \plot 0 0  1 1  2 0  3 1  4 0  5 1  6 0  7 1  8 0 /
 \plot  1 1  2 2  3 1  4 2  5 1  6 2  7 1 /
  \plot 3 0  5 2  5 1  /
   \plot 4 0  4 1  3 2  3 1  /
   \multiput{$\bullet$} at 0 0  3 1  4 1  5 1  8 0  /
\endpicture} at 1 -3

\put{
\beginpicture
 \setcoordinatesystem units <.2cm,.3cm>
 \multiput{$\circ$} at
   0 0  2 0  4 0  6 0  8 0  1 1  3 1  5 1  7 1   2 2  4 2  6 2
      3 2  5 2  4 1  3 0 /
       \setdots <.5mm>
        \plot 0 0  1 1  2 0  3 1  4 0  5 1  6 0  7 1  8 0 /
 \plot  1 1  2 2  3 1  4 2  5 1  6 2  7 1 /
  \plot 3 0  5 2  5 1  /
   \plot 4 0  4 1  3 2  3 1  /
   \multiput{$\bullet$} at 1 1  3 1  4 1  5 1  8 0  /
\endpicture} at 0 -2

\put{
\beginpicture
 \setcoordinatesystem units <.2cm,.3cm>
 \multiput{$\circ$} at
   0 0  2 0  4 0  6 0  8 0  1 1  3 1  5 1  7 1   2 2  4 2  6 2
      3 2  5 2  4 1  3 0 /
       \setdots <.5mm>
        \plot 0 0  1 1  2 0  3 1  4 0  5 1  6 0  7 1  8 0 /
 \plot  1 1  2 2  3 1  4 2  5 1  6 2  7 1 /
  \plot 3 0  5 2  5 1  /
   \plot 4 0  4 1  3 2  3 1  /
   \multiput{$\bullet$} at 0 0  3 1  4 1  5 1  7 1  /
\endpicture} at 2 -2

\put{
\beginpicture
 \setcoordinatesystem units <.2cm,.3cm>
 \multiput{$\circ$} at
   0 0  2 0  4 0  6 0  8 0  1 1  3 1  5 1  7 1   2 2  4 2  6 2
      3 2  5 2  4 1  3 0 /
       \setdots <.5mm>
        \plot 0 0  1 1  2 0  3 1  4 0  5 1  6 0  7 1  8 0 /
 \plot  1 1  2 2  3 1  4 2  5 1  6 2  7 1 /
  \plot 3 0  5 2  5 1  /
   \plot 4 0  4 1  3 2  3 1  /
   \multiput{$\bullet$} at 1 1  3 1  4 1  5 1  7 1  /
\endpicture} at 1 -1

\plot 0.3 -2.3  0.7 -2.7 /
\plot 1.3 -2.7  1.7 -2.3 /
\plot .3 -1.7  .7 -1.3 /
\plot 1.3 -1.3  1.7 -1.7 /

\plot -.3 -2.3 -.5 -2.5 /
\plot 2.3 -2.3 2.5 -2.5 /

\setdots <1mm>
\plot -1 1  0 -2 /
\plot 0 2  1 -1 /
\plot 0 0  1 -3 /
\plot 1 1  2 -2 /

\setsolid

\put{
\beginpicture
 \setcoordinatesystem units <.2cm,.3cm>
 \multiput{$\circ$} at
   0 0  2 0  4 0  6 0  8 0  1 1  3 1  5 1  7 1   2 2  4 2  6 2
      3 2  5 2  4 1  3 0 /
       \setdots <.5mm>
        \plot 0 0  1 1  2 0  3 1  4 0  5 1  6 0  7 1  8 0 /
 \plot  1 1  2 2  3 1  4 2  5 1  6 2  7 1 /
  \plot 3 0  5 2  5 1  /
   \plot 4 0  4 1  3 2  3 1  /
   \multiput{$\bullet$} at 0 0  3 2  4 2  5 2  8 0  /
\endpicture} at 0 0

\put{
\beginpicture
 \setcoordinatesystem units <.2cm,.3cm>
 \multiput{$\circ$} at
   0 0  2 0  4 0  6 0  8 0  1 1  3 1  5 1  7 1   2 2  4 2  6 2
      3 2  5 2  4 1  3 0 /
       \setdots <.5mm>
        \plot 0 0  1 1  2 0  3 1  4 0  5 1  6 0  7 1  8 0 /
 \plot  1 1  2 2  3 1  4 2  5 1  6 2  7 1 /
  \plot 3 0  5 2  5 1  /
   \plot 4 0  4 1  3 2  3 1  /
   \multiput{$\bullet$} at 1 1  3 2  4 2  5 2  8 0  /
\endpicture} at -1 1

\put{
\beginpicture
 \setcoordinatesystem units <.2cm,.3cm>
 \multiput{$\circ$} at
   0 0  2 0  4 0  6 0  8 0  1 1  3 1  5 1  7 1   2 2  4 2  6 2
      3 2  5 2  4 1  3 0 /
       \setdots <.5mm>
        \plot 0 0  1 1  2 0  3 1  4 0  5 1  6 0  7 1  8 0 /
 \plot  1 1  2 2  3 1  4 2  5 1  6 2  7 1 /
  \plot 3 0  5 2  5 1  /
   \plot 4 0  4 1  3 2  3 1  /
   \multiput{$\bullet$} at 0 0  3 2  4 2  5 2  7 1  /
\endpicture} at 1 1

\put{
\beginpicture
 \setcoordinatesystem units <.2cm,.3cm>
 \multiput{$\circ$} at
   0 0  2 0  4 0  6 0  8 0  1 1  3 1  5 1  7 1   2 2  4 2  6 2
      3 2  5 2  4 1  3 0 /
       \setdots <.5mm>
        \plot 0 0  1 1  2 0  3 1  4 0  5 1  6 0  7 1  8 0 /
 \plot  1 1  2 2  3 1  4 2  5 1  6 2  7 1 /
  \plot 3 0  5 2  5 1  /
   \plot 4 0  4 1  3 2  3 1  /
   \multiput{$\bullet$} at 2 2  3 2  4 2  5 2  8 0  /
\endpicture} at -2 2

\put{
\beginpicture
 \setcoordinatesystem units <.2cm,.3cm>
 \multiput{$\circ$} at
   0 0  2 0  4 0  6 0  8 0  1 1  3 1  5 1  7 1   2 2  4 2  6 2
      3 2  5 2  4 1  3 0 /
       \setdots <.5mm>
        \plot 0 0  1 1  2 0  3 1  4 0  5 1  6 0  7 1  8 0 /
 \plot  1 1  2 2  3 1  4 2  5 1  6 2  7 1 /
  \plot 3 0  5 2  5 1  /
   \plot 4 0  4 1  3 2  3 1  /
   \multiput{$\bullet$} at 1 1  3 2  4 2  5 2  7 1  /
\endpicture} at 0 2

\put{
\beginpicture
 \setcoordinatesystem units <.2cm,.3cm>
 \multiput{$\circ$} at
   0 0  2 0  4 0  6 0  8 0  1 1  3 1  5 1  7 1   2 2  4 2  6 2
      3 2  5 2  4 1  3 0 /
       \setdots <.5mm>
        \plot 0 0  1 1  2 0  3 1  4 0  5 1  6 0  7 1  8 0 /
 \plot  1 1  2 2  3 1  4 2  5 1  6 2  7 1 /
  \plot 3 0  5 2  5 1  /
   \plot 4 0  4 1  3 2  3 1  /
   \multiput{$\bullet$} at 0 0  3 2  4 2  5 2  6 2  /
\endpicture} at 2 2

\put{
\beginpicture
 \setcoordinatesystem units <.2cm,.3cm>
 \multiput{$\circ$} at
   0 0  2 0  4 0  6 0  8 0  1 1  3 1  5 1  7 1   2 2  4 2  6 2
      3 2  5 2  4 1  3 0 /
       \setdots <.5mm>
        \plot 0 0  1 1  2 0  3 1  4 0  5 1  6 0  7 1  8 0 /
 \plot  1 1  2 2  3 1  4 2  5 1  6 2  7 1 /
  \plot 3 0  5 2  5 1  /
   \plot 4 0  4 1  3 2  3 1  /
   \multiput{$\bullet$} at 2 2  3 2  4 2  5 2  7 1  /
\endpicture} at -1 3

\put{
\beginpicture
 \setcoordinatesystem units <.2cm,.3cm>
 \multiput{$\circ$} at
   0 0  2 0  4 0  6 0  8 0  1 1  3 1  5 1  7 1   2 2  4 2  6 2
      3 2  5 2  4 1  3 0 /
       \setdots <.5mm>
        \plot 0 0  1 1  2 0  3 1  4 0  5 1  6 0  7 1  8 0 /
 \plot  1 1  2 2  3 1  4 2  5 1  6 2  7 1 /
  \plot 3 0  5 2  5 1  /
   \plot 4 0  4 1  3 2  3 1  /
   \multiput{$\bullet$} at 1 1  3 2  4 2  5 2  7 1  /
\endpicture} at 1 3

\put{
\beginpicture
 \setcoordinatesystem units <.2cm,.3cm>
 \multiput{$\circ$} at
   0 0  2 0  4 0  6 0  8 0  1 1  3 1  5 1  7 1   2 2  4 2  6 2
      3 2  5 2  4 1  3 0 /
       \setdots <.5mm>
        \plot 0 0  1 1  2 0  3 1  4 0  5 1  6 0  7 1  8 0 /
 \plot  1 1  2 2  3 1  4 2  5 1  6 2  7 1 /
  \plot 3 0  5 2  5 1  /
   \plot 4 0  4 1  3 2  3 1  /
   \multiput{$\bullet$} at 2 2  3 2  4 2  5 2  6 2  /
\endpicture} at 0 4

\put{
\beginpicture
 \setcoordinatesystem units <.2cm,.3cm>
 \multiput{$\circ$} at
   0 0  2 0  4 0  6 0  8 0  1 1  3 1  5 1  7 1   2 2  4 2  6 2
      3 2  5 2  4 1  3 0 /
       \setdots <.5mm>
        \plot 0 0  1 1  2 0  3 1  4 0  5 1  6 0  7 1  8 0 /
 \plot  1 1  2 2  3 1  4 2  5 1  6 2  7 1 /
  \plot 3 0  5 2  5 1  /
   \plot 4 0  4 1  3 2  3 1  /
   \multiput{$\circ$} at 2 2  3 2  4 2  5 2  6 2  /
\multiput{$\bullet$} at 2 3  3 3  4 3  5 3  6 3  /
 \plot 2 2  2 3  3 2  4 3  5 2  6 3  6 2  5 3  4 2  3 3  2 2 /
  \plot 4 2  4 3 /
  \endpicture} at 0 5

\put{
\beginpicture
 \setcoordinatesystem units <.2cm,.3cm>
 \multiput{$\circ$} at
   0 0  2 0  4 0  6 0  8 0  1 1  3 1  5 1  7 1   2 2  4 2  6 2
      3 2  5 2  4 1  3 0 /
       \setdots <.5mm>
        \plot 0 0  1 1  2 0  3 1  4 0  5 1  6 0  7 1  8 0 /
 \plot  1 1  2 2  3 1  4 2  5 1  6 2  7 1 /
  \plot 3 0  5 2  5 1  /
   \plot 4 0  4 1  3 2  3 1  /
   \multiput{$\circ$} at 2 2  3 2  4 2  5 2  6 2  4 3  /
\multiput{$\bullet$} at 2 3  3 3  4 4  5 3  6 3  /
 \plot 2 2  2 3  3 2  4 3  5 2  6 3  6 2  5 3  4 2  3 3  2 2 /
  \plot 4 2  4 4 /
  \endpicture} at 0 6.1

\plot 0.3 0.3  0.7 0.7 /
\plot -0.3 0.3 -0.7 0.7 /
\plot 1.3 1.3  1.7 1.7 /
\plot -1.3 1.3 -1.7 1.7 /
\plot  0.3 1.7  0.7 1.3 /
\plot -0.3 1.7  -.7 1.3 /
\plot  1.7 2.3  1.3 2.7 /
\plot  -1.7 2.3  -1.3 2.7 /
\plot 0.3 2.3  0.7 2.7 /
\plot -.3 2.3  -.7 2.7 /
\plot 0.3 3.7  0.7 3.3 /
\plot -.3 3.7  -.7 3.3 /

\setdots <.5mm>
\plot 0 4.35  0 4.65 /
\setsolid
\plot 0 5.38  0 5.65 /

\put{$\Phi_2$} at  1.85 -1
\put{$\Phi_3$} at  0.65 4
\put{$\Phi_4$} at  0.65 5
\put{$\Phi_4'$} at  0.65 6

\setshadegrid span <.6mm>
\vshade -2.5 2 2 <z,z,,> 0 -.5 4.5  <z,z,,> 2.5 2 2    /
\vshade -.5 4.5 6.5 <z,z,,> .5 4.5 6.5     /
\endpicture} at 0 0 

\endpicture}
$$
We see that there is a poset isomorphism between $\Phi_+(\Bbb E_5)$ and 
$\Cal R(\Bbb E_6).$

\vfill\eject
$$
{\beginpicture
\setcoordinatesystem units <1cm,1cm>
\put{\beginpicture
\setcoordinatesystem units <.9cm,.63cm>
\put{$2\ \,\atop123321$} at 0 9
\put{$2\ \,\atop123221$} at -1 8
\put{$2\ \,\atop123211$} at -2 7
\put{$2\ \,\atop123210$} at  -3 6
\put{$1\ \,\atop123321$} at  1 8
\put{$1\ \,\atop123221$} at  0 7
\put{$1\ \,\atop123211$} at -1 6
\put{$1\ \,\atop123210$} at -2 5
\setdots <1mm>
\plot 0 7  1 6 /
\plot -1 6  0 5 /
\plot -2 5  -1 4 /

\setsolid

\plot .7 8.3  0.3 8.7 /
\plot -.3 7.3  -.7 7.7 /
\plot -1.3 6.3  -1.7 6.7  /
\plot -2.3 5.3  -2.7 5.7  /
\plot -1.7 5.3  -1.3 5.7 /
\plot -.7 6.3  -.3 6.7 /
\plot .3 7.3  .7 7.7 /
\plot -2.7 6.3  -2.3 6.7 /
\plot -1.7 7.3  -1.3 7.7 /
\plot -.7 8.3  -.3 8.7 /

\put{$1\ \,\atop012100$} at 0 0
\put{$1\ \,\atop112100$} at -1 1
\put{$1\ \,\atop012110$} at  1 1
\put{$1\ \,\atop122100$} at -2 2
\put{$1\ \,\atop112110$} at 0 2
\put{$1\ \,\atop012111$} at 2 2
\put{$1\ \,\atop122110$} at -1 3
\put{$1\ \,\atop112111$} at 1 3
\put{$1\ \,\atop122111$} at 0 4

\put{$1\ \,\atop012210$} at 1 2
\put{$1\ \,\atop112210$} at 0 3
\put{$1\ \,\atop012211$} at 2 3
\put{$1\ \,\atop122210$} at -1 4  
\put{$1\ \,\atop112211$} at 1 4
\put{$1\ \,\atop012221$} at 3 4
\put{$1\ \,\atop122211$} at 0 5
\put{$1\ \,\atop112221$} at 2 5
\put{$1\ \,\atop122221$} at 1 6

\plot 0.3 0.3  0.7 0.7 /
\plot -0.3 0.3 -0.7 0.7 /
\plot 1.3 1.3  1.7 1.7 /
\plot -1.3 1.3 -1.7 1.7 /
\plot  0.3 1.7  0.7 1.3 /
\plot -0.3 1.7  -.7 1.3 /
\plot  1.7 2.3  1.3 2.7 /
\plot  -1.7 2.3  -1.3 2.7 /
\plot 0.3 2.3  0.7 2.7 /
\plot -.3 2.3  -.7 2.7 /
\plot 0.3 3.7  0.7 3.3 /
\plot -.3 3.7  -.7 3.3 /

\plot 1.3 2.3  1.7 2.7 /
\plot 0.7 2.3 0.3 2.7 /
\plot 2.3 3.3  2.7 3.7 /
\plot -.3 3.3 -.7 3.7 /
\plot  1.3 3.7  1.7 3.3 /
\plot  .7 3.7  .3 3.3 /
\plot  2.7 4.3  2.3 4.7 /
\plot  -.7 4.3  -.3 4.7 /
\plot 1.3 4.3  1.7 4.7 /
\plot .7 4.3  .3 4.7 /
\plot 1.3 5.7  1.7 5.3 /
\plot .7 5.7  .3 5.3 /

\plot 1 1.3  1 1.7 /
\plot 0 2.3  0 2.7  /
\plot 2 2.3  2 2.7 /
\plot -1 3.3 -1 3.7 /
\plot 1 3.3  1 3.7 /
\plot 0 4.3  0 4.7 /

\put{$1\ \,\atop011100$} at 1 -3
\put{$1\ \,\atop111100$} at 0 -2
\put{$1\ \,\atop011110$} at 2 -2
\put{$1\ \,\atop111110$} at 1 -1
\put{$1\ \,\atop011111$} at 3 -1
\put{$1\ \,\atop111111$} at 2 0

\plot 0.3 -2.3  0.7 -2.7 /
\plot 1.3 -1.3  1.7 -1.7 /
\plot 2.3 -.3  2.7 -.7 /
\plot 1.3 -2.7  1.7 -2.3 /
\plot 2.3 -1.7  2.7 -1.3 /
\plot .3 -1.7  .7 -1.3 /
\plot 1.3 -0.7  1.7 -.3 /

\setdots <1mm>
\plot -1 1  0 -2 /
\plot 0 2  1 -1 /
\plot 1 3  2 0 /
\plot 0 0  1 -3 /
\plot 1 1  2 -2 /
\plot 2 2  3 -1 /

\setshadegrid span <.6mm>
\vshade -3.5 6 6 <z,z,,> -2 4.5 7.5 <z,z,,> 0 6.5 9.5 <z,z,,> 1.5 8 8   /
\vshade -2.5 2 2  <z,z,,> -1.4 0.5 3 <z,z,,> -1.3 0.6 4 <z,z,,> 0 -0.5 5.5   
   <z,z,,> 1 0.5 6.5  <z,z,,> 2.3 1.8 5.2  <z,z,,> 2.4 2.8 5.1 <z,z,,> 3.5 4 4  /
\put{$\Phi(\Bbb E_7,3)$} at -2.5 8.5 
\put{$\Phi(\Bbb E_7,2)$} at -2 -.2

\endpicture} at 5.9 5.8

\put{\beginpicture
\put{\bf The case $\Bbb E_7.$} at -4 14
\setcoordinatesystem units <1.7cm,1.25cm>
\put{$\Phi_2$} at 2.8 0
\put{$\Phi_3$} at 1.7 6
\put{$\Phi_4$} at 1.7 8
\put{$\Phi_4'$} at 0.7 9
\put
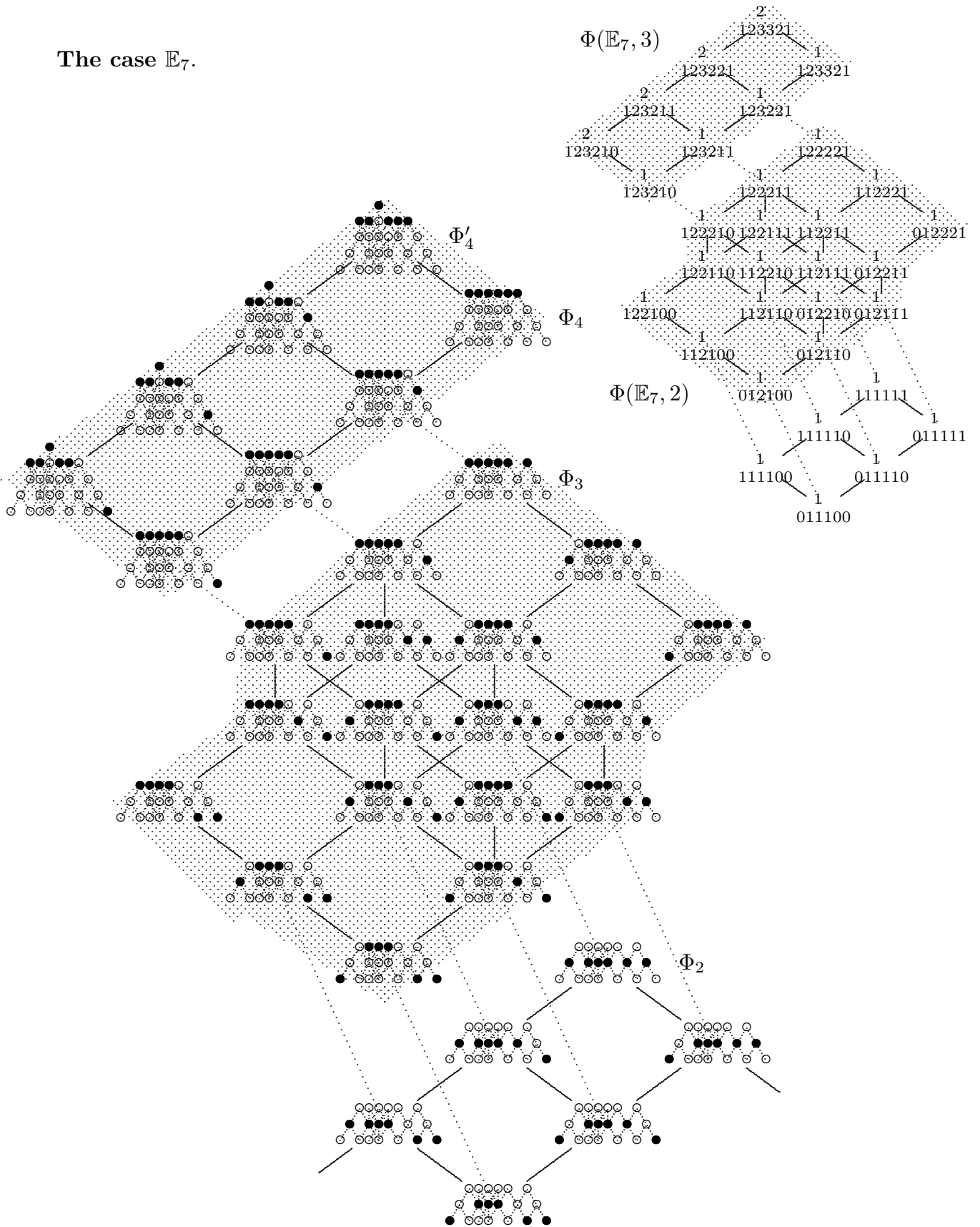
{
\beginpicture
 \setcoordinatesystem units <.15cm,.25cm>
 \multiput{$\circ$} at
   0 0  2 0  4 0  6 0  8 0  10 0
      1 1  3 1  5 1  7 1  9 1
         2 2  4 2  6 2  8 2
	    3 2  5 2  4 1  3 0 /
	     \setdots <.5mm>
	      \plot 0 0  1 1  2 0  3 1  4 0  5 1  6 0  7 1  8 0  9 1 /
 \plot  1 1  2 2  3 1  4 2  5 1  6 2  7 1  8 2  10 0 /
  \plot 3 0  5 2  5 1  /
   \plot 4 0  4 1  3 2  3 1  /
  \multiput{$\circ$} at    2 3  3 3  4 3  5 3  6 3  7 3 /
  \plot 2 2  3 3  4 2  4 3  3 2  2 3  2 2 /
  \plot 4 2  5 3  6 2  7 3  8 2 /
  \plot 6 2  6 3 /
  \put{$\circ$} at 4 4
  \plot 4 3  4 4 /
  \multiput{$\bullet$} at 2 3  3 3  4 4  5 3  6 3  7 3   /
\endpicture} at 0 9

\put{
\beginpicture
 \setcoordinatesystem units <.15cm,.25cm>
 \multiput{$\circ$} at
   0 0  2 0  4 0  6 0  8 0  10 0
      1 1  3 1  5 1  7 1  9 1
         2 2  4 2  6 2  8 2
	    3 2  5 2  4 1  3 0 /
	     \setdots <.5mm>
	      \plot 0 0  1 1  2 0  3 1  4 0  5 1  6 0  7 1  8 0  9 1 /
 \plot  1 1  2 2  3 1  4 2  5 1  6 2  7 1  8 2  10 0 /
  \plot 3 0  5 2  5 1  /
   \plot 4 0  4 1  3 2  3 1  /
  \multiput{$\circ$} at    2 3  3 3  4 3  5 3  6 3  7 3 /
  \plot 2 2  3 3  4 2  4 3  3 2  2 3  2 2 /
  \plot 4 2  5 3  6 2  7 3  8 2 /
  \plot 6 2  6 3 /
  \put{$\circ$} at 4 4
  \plot 4 3  4 4 /
  \multiput{$\bullet$} at 2 3  3 3  4 4  5 3  6 3  8 2   /
\endpicture} at -1 8

\put{
\beginpicture
 \setcoordinatesystem units <.15cm,.25cm>
 \multiput{$\circ$} at
   0 0  2 0  4 0  6 0  8 0  10 0
      1 1  3 1  5 1  7 1  9 1
         2 2  4 2  6 2  8 2
	    3 2  5 2  4 1  3 0 /
	     \setdots <.5mm>
	      \plot 0 0  1 1  2 0  3 1  4 0  5 1  6 0  7 1  8 0  9 1 /
 \plot  1 1  2 2  3 1  4 2  5 1  6 2  7 1  8 2  10 0 /
  \plot 3 0  5 2  5 1  /
   \plot 4 0  4 1  3 2  3 1  /
  \multiput{$\circ$} at    2 3  3 3  4 3  5 3  6 3  7 3 /
  \plot 2 2  3 3  4 2  4 3  3 2  2 3  2 2 /
  \plot 4 2  5 3  6 2  7 3  8 2 /
  \plot 6 2  6 3 /
  \put{$\circ$} at 4 4
  \plot 4 3  4 4 /
  \multiput{$\bullet$} at 2 3  3 3  4 4  5 3  6 3  9 1   /
\endpicture} at -2 7

\put{
\beginpicture
 \setcoordinatesystem units <.15cm,.25cm>
 \multiput{$\circ$} at
   0 0  2 0  4 0  6 0  8 0  10 0
      1 1  3 1  5 1  7 1  9 1
         2 2  4 2  6 2  8 2
	    3 2  5 2  4 1  3 0 /
	     \setdots <.5mm>
	      \plot 0 0  1 1  2 0  3 1  4 0  5 1  6 0  7 1  8 0  9 1 /
 \plot  1 1  2 2  3 1  4 2  5 1  6 2  7 1  8 2  10 0 /
  \plot 3 0  5 2  5 1  /
   \plot 4 0  4 1  3 2  3 1  /
  \multiput{$\circ$} at    2 3  3 3  4 3  5 3  6 3  7 3 /
  \plot 2 2  3 3  4 2  4 3  3 2  2 3  2 2 /
  \plot 4 2  5 3  6 2  7 3  8 2 /
  \plot 6 2  6 3 /
  \put{$\circ$} at 4 4
  \plot 4 3  4 4 /
  \multiput{$\bullet$} at 2 3  3 3  4 4  5 3  6 3  10 0   /
\endpicture} at -3 6

\put{
\beginpicture
 \setcoordinatesystem units <.15cm,.25cm>
 \multiput{$\circ$} at
   0 0  2 0  4 0  6 0  8 0  10 0
      1 1  3 1  5 1  7 1  9 1
         2 2  4 2  6 2  8 2
	    3 2  5 2  4 1  3 0 /
	     \setdots <.5mm>
	      \plot 0 0  1 1  2 0  3 1  4 0  5 1  6 0  7 1  8 0  9 1 /
 \plot  1 1  2 2  3 1  4 2  5 1  6 2  7 1  8 2  10 0 /
  \plot 3 0  5 2  5 1  /
   \plot 4 0  4 1  3 2  3 1  /
  \multiput{$\circ$} at    2 3  3 3  4 3  5 3  6 3  7 3 /
  \plot 2 2  3 3  4 2  4 3  3 2  2 3  2 2 /
  \plot 4 2  5 3  6 2  7 3  8 2 /
  \plot 6 2  6 3 /
  \multiput{$\bullet$} at 2 3  3 3  4 3  5 3  6 3  7 3  /
\endpicture} at  1 8

\put{
\beginpicture
 \setcoordinatesystem units <.15cm,.25cm>
 \multiput{$\circ$} at
   0 0  2 0  4 0  6 0  8 0  10 0
      1 1  3 1  5 1  7 1  9 1
         2 2  4 2  6 2  8 2
	    3 2  5 2  4 1  3 0 /
	     \setdots <.5mm>
	      \plot 0 0  1 1  2 0  3 1  4 0  5 1  6 0  7 1  8 0  9 1 /
 \plot  1 1  2 2  3 1  4 2  5 1  6 2  7 1  8 2  10 0 /
  \plot 3 0  5 2  5 1  /
   \plot 4 0  4 1  3 2  3 1  /
  \multiput{$\circ$} at    2 3  3 3  4 3  5 3  6 3  7 3 /
  \plot 2 2  3 3  4 2  4 3  3 2  2 3  2 2 /
  \plot 4 2  5 3  6 2  7 3  8 2 /
  \plot 6 2  6 3 /
  \multiput{$\bullet$} at 2 3  3 3  4 3  5 3  6 3  8 2  /
\endpicture} at 0 7

\put{
\beginpicture
 \setcoordinatesystem units <.15cm,.25cm>
 \multiput{$\circ$} at
   0 0  2 0  4 0  6 0  8 0  10 0
      1 1  3 1  5 1  7 1  9 1
         2 2  4 2  6 2  8 2
	    3 2  5 2  4 1  3 0 /
	     \setdots <.5mm>
	      \plot 0 0  1 1  2 0  3 1  4 0  5 1  6 0  7 1  8 0  9 1 /
 \plot  1 1  2 2  3 1  4 2  5 1  6 2  7 1  8 2  10 0 /
  \plot 3 0  5 2  5 1  /
   \plot 4 0  4 1  3 2  3 1  /
  \multiput{$\circ$} at    2 3  3 3  4 3  5 3  6 3  7 3 /
  \plot 2 2  3 3  4 2  4 3  3 2  2 3  2 2 /
  \plot 4 2  5 3  6 2  7 3  8 2 /
  \plot 6 2  6 3 /
  \multiput{$\bullet$} at 2 3  3 3  4 3  5 3  6 3  9 1  /
\endpicture} at -1 6

\put{
\beginpicture
 \setcoordinatesystem units <.15cm,.25cm>
 \multiput{$\circ$} at
   0 0  2 0  4 0  6 0  8 0  10 0
      1 1  3 1  5 1  7 1  9 1
         2 2  4 2  6 2  8 2
	    3 2  5 2  4 1  3 0 /
	     \setdots <.5mm>
	      \plot 0 0  1 1  2 0  3 1  4 0  5 1  6 0  7 1  8 0  9 1 /
 \plot  1 1  2 2  3 1  4 2  5 1  6 2  7 1  8 2  10 0 /
  \plot 3 0  5 2  5 1  /
   \plot 4 0  4 1  3 2  3 1  /
  \multiput{$\circ$} at    2 3  3 3  4 3  5 3  6 3  7 3 /
  \plot 2 2  3 3  4 2  4 3  3 2  2 3  2 2 /
  \plot 4 2  5 3  6 2  7 3  8 2 /
  \plot 6 2  6 3 /
  \multiput{$\bullet$} at 2 3  3 3  4 3  5 3  6 3  10 0  /
\endpicture} at -2 5

\setdots <1mm>
\plot 0 7  1 6 /
\plot -1 6  0 5 /
\plot -2 5  -1 4 /

\setsolid

\plot .7 8.3  0.3 8.7 /
\plot -.3 7.3  -.7 7.7 /
\plot -1.3 6.3  -1.7 6.7  /
\plot -2.3 5.3  -2.7 5.7  /
\plot -1.7 5.3  -1.3 5.7 /
\plot -.7 6.3  -.3 6.7 /
\plot .3 7.3  .7 7.7 /
\plot -2.7 6.3  -2.3 6.7 /
\plot -1.7 7.3  -1.3 7.7 /
\plot -.7 8.3  -.3 8.7 /

\put{
\beginpicture
 \setcoordinatesystem units <.15cm,.25cm>
 \multiput{$\circ$} at
   0 0  2 0  4 0  6 0  8 0  10 0
      1 1  3 1  5 1  7 1  9 1
         2 2  4 2  6 2  8 2
	    3 2  5 2  4 1  3 0 /
	     \setdots <.5mm>
	      \plot 0 0  1 1  2 0  3 1  4 0  5 1  6 0  7 1  8 0  9 1 /
 \plot  1 1  2 2  3 1  4 2  5 1  6 2  7 1  8 2  10 0 /
  \plot 3 0  5 2  5 1  /
   \plot 4 0  4 1  3 2  3 1  /
   \multiput{$\bullet$} at 0 0  3 2  4 2  5 2  8 0  10 0  /
\endpicture} at 0 0

\put{
\beginpicture
 \setcoordinatesystem units <.15cm,.25cm>
 \multiput{$\circ$} at
   0 0  2 0  4 0  6 0  8 0  10 0
      1 1  3 1  5 1  7 1  9 1
         2 2  4 2  6 2  8 2
	    3 2  5 2  4 1  3 0 /
	     \setdots <.5mm>
	      \plot 0 0  1 1  2 0  3 1  4 0  5 1  6 0  7 1  8 0  9 1 /
 \plot  1 1  2 2  3 1  4 2  5 1  6 2  7 1  8 2  10 0 /
  \plot 3 0  5 2  5 1  /
   \plot 4 0  4 1  3 2  3 1  /
   \multiput{$\bullet$} at 1 1  3 2  4 2  5 2  8 0  10 0  /
\endpicture} at -1 1

\put{
\beginpicture
 \setcoordinatesystem units <.15cm,.25cm>
 \multiput{$\circ$} at
   0 0  2 0  4 0  6 0  8 0  10 0
      1 1  3 1  5 1  7 1  9 1
         2 2  4 2  6 2  8 2
	    3 2  5 2  4 1  3 0 /
	     \setdots <.5mm>
	      \plot 0 0  1 1  2 0  3 1  4 0  5 1  6 0  7 1  8 0  9 1 /
 \plot  1 1  2 2  3 1  4 2  5 1  6 2  7 1  8 2  10 0 /
  \plot 3 0  5 2  5 1  /
   \plot 4 0  4 1  3 2  3 1  /
   \multiput{$\bullet$} at 0 0  3 2  4 2  5 2  7 1  10 0  /
\endpicture} at 1 1

\put{
\beginpicture
 \setcoordinatesystem units <.15cm,.25cm>
 \multiput{$\circ$} at
   0 0  2 0  4 0  6 0  8 0  10 0
      1 1  3 1  5 1  7 1  9 1
         2 2  4 2  6 2  8 2
	    3 2  5 2  4 1  3 0 /
	     \setdots <.5mm>
	      \plot 0 0  1 1  2 0  3 1  4 0  5 1  6 0  7 1  8 0  9 1 /
 \plot  1 1  2 2  3 1  4 2  5 1  6 2  7 1  8 2  10 0 /
  \plot 3 0  5 2  5 1  /
   \plot 4 0  4 1  3 2  3 1  /
   \multiput{$\bullet$} at 2 2  3 2  4 2  5 2  8 0  10 0  /
\endpicture} at -2 2

\put{
\beginpicture
 \setcoordinatesystem units <.15cm,.25cm>
 \multiput{$\circ$} at
   0 0  2 0  4 0  6 0  8 0  10 0
      1 1  3 1  5 1  7 1  9 1
         2 2  4 2  6 2  8 2
	    3 2  5 2  4 1  3 0 /
	     \setdots <.5mm>
	      \plot 0 0  1 1  2 0  3 1  4 0  5 1  6 0  7 1  8 0  9 1 /
 \plot  1 1  2 2  3 1  4 2  5 1  6 2  7 1  8 2  10 0 /
  \plot 3 0  5 2  5 1  /
   \plot 4 0  4 1  3 2  3 1  /
   \multiput{$\bullet$} at 1 1  3 2  4 2  5 2  7 1  10 0  /
\endpicture} at 0 2

\put{
\beginpicture
 \setcoordinatesystem units <.15cm,.25cm>
 \multiput{$\circ$} at
   0 0  2 0  4 0  6 0  8 0  10 0
      1 1  3 1  5 1  7 1  9 1
         2 2  4 2  6 2  8 2
	    3 2  5 2  4 1  3 0 /
	     \setdots <.5mm>
	      \plot 0 0  1 1  2 0  3 1  4 0  5 1  6 0  7 1  8 0  9 1 /
 \plot  1 1  2 2  3 1  4 2  5 1  6 2  7 1  8 2  10 0 /
  \plot 3 0  5 2  5 1  /
   \plot 4 0  4 1  3 2  3 1  /
   \multiput{$\bullet$} at 0 0  3 2  4 2  5 2   7 1  9 1  /
\endpicture} at 2 2

\put{
\beginpicture
 \setcoordinatesystem units <.15cm,.25cm>
 \multiput{$\circ$} at
   0 0  2 0  4 0  6 0  8 0  10 0
      1 1  3 1  5 1  7 1  9 1
         2 2  4 2  6 2  8 2
	    3 2  5 2  4 1  3 0 /
	     \setdots <.5mm>
	      \plot 0 0  1 1  2 0  3 1  4 0  5 1  6 0  7 1  8 0  9 1 /
 \plot  1 1  2 2  3 1  4 2  5 1  6 2  7 1  8 2  10 0 /
  \plot 3 0  5 2  5 1  /
   \plot 4 0  4 1  3 2  3 1  /
   \multiput{$\bullet$} at 2 2  3 2  4 2  5 2  7 1  10 0   /
\endpicture} at -1 3

\put{
\beginpicture
 \setcoordinatesystem units <.15cm,.25cm>
 \multiput{$\circ$} at
   0 0  2 0  4 0  6 0  8 0  10 0
      1 1  3 1  5 1  7 1  9 1
         2 2  4 2  6 2  8 2
	    3 2  5 2  4 1  3 0 /
	     \setdots <.5mm>
	      \plot 0 0  1 1  2 0  3 1  4 0  5 1  6 0  7 1  8 0  9 1 /
 \plot  1 1  2 2  3 1  4 2  5 1  6 2  7 1  8 2  10 0 /
  \plot 3 0  5 2  5 1  /
   \plot 4 0  4 1  3 2  3 1  /
   \multiput{$\bullet$} at 1 1  3 2  4 2  5 2  7 1  9 1  /
\endpicture} at 1 3

\put{
\beginpicture
 \setcoordinatesystem units <.15cm,.25cm>
 \multiput{$\circ$} at
   0 0  2 0  4 0  6 0  8 0  10 0
      1 1  3 1  5 1  7 1  9 1
         2 2  4 2  6 2  8 2
	    3 2  5 2  4 1  3 0 /
	     \setdots <.5mm>
	      \plot 0 0  1 1  2 0  3 1  4 0  5 1  6 0  7 1  8 0  9 1 /
 \plot  1 1  2 2  3 1  4 2  5 1  6 2  7 1  8 2  10 0 /
  \plot 3 0  5 2  5 1  /
   \plot 4 0  4 1  3 2  3 1  /
   \multiput{$\bullet$} at 2 2  3 2  4 2  5 2   7 1  9 1  /
\endpicture} at 0 4

\put{
\beginpicture
 \setcoordinatesystem units <.15cm,.25cm>
 \multiput{$\circ$} at
   0 0  2 0  4 0  6 0  8 0  10 0
      1 1  3 1  5 1  7 1  9 1
         2 2  4 2  6 2  8 2
	    3 2  5 2  4 1  3 0 /
	     \setdots <.5mm>
	      \plot 0 0  1 1  2 0  3 1  4 0  5 1  6 0  7 1  8 0  9 1 /
 \plot  1 1  2 2  3 1  4 2  5 1  6 2  7 1  8 2  10 0 /
  \plot 3 0  5 2  5 1  /
   \plot 4 0  4 1  3 2  3 1  /
   \multiput{$\bullet$} at 1 1  3 2  4 2  5 2  6 2  10 0  /
\endpicture} at 1 2

\put{
\beginpicture
 \setcoordinatesystem units <.15cm,.25cm>
 \multiput{$\circ$} at
   0 0  2 0  4 0  6 0  8 0  10 0
      1 1  3 1  5 1  7 1  9 1
         2 2  4 2  6 2  8 2
	    3 2  5 2  4 1  3 0 /
	     \setdots <.5mm>
	      \plot 0 0  1 1  2 0  3 1  4 0  5 1  6 0  7 1  8 0  9 1 /
 \plot  1 1  2 2  3 1  4 2  5 1  6 2  7 1  8 2  10 0 /
  \plot 3 0  5 2  5 1  /
   \plot 4 0  4 1  3 2  3 1  /
   \multiput{$\bullet$} at 1 1  3 2  4 2  5 2  6 2  10 0  /
\endpicture} at 0 3

\put{
\beginpicture
 \setcoordinatesystem units <.15cm,.25cm>
 \multiput{$\circ$} at
   0 0  2 0  4 0  6 0  8 0  10 0
      1 1  3 1  5 1  7 1  9 1
         2 2  4 2  6 2  8 2
	    3 2  5 2  4 1  3 0 /
	     \setdots <.5mm>
	      \plot 0 0  1 1  2 0  3 1  4 0  5 1  6 0  7 1  8 0  9 1 /
 \plot  1 1  2 2  3 1  4 2  5 1  6 2  7 1  8 2  10 0 /
  \plot 3 0  5 2  5 1  /
   \plot 4 0  4 1  3 2  3 1  /
   \multiput{$\bullet$} at 0 0  3 2  4 2  5 2  6 2  9 1  /
\endpicture} at 2 3

\put{
\beginpicture
 \setcoordinatesystem units <.15cm,.25cm>
 \multiput{$\circ$} at
   0 0  2 0  4 0  6 0  8 0  10 0
      1 1  3 1  5 1  7 1  9 1
         2 2  4 2  6 2  8 2
	    3 2  5 2  4 1  3 0 /
	     \setdots <.5mm>
	      \plot 0 0  1 1  2 0  3 1  4 0  5 1  6 0  7 1  8 0  9 1 /
 \plot  1 1  2 2  3 1  4 2  5 1  6 2  7 1  8 2  10 0 /
  \plot 3 0  5 2  5 1  /
   \plot 4 0  4 1  3 2  3 1  /
   \multiput{$\bullet$} at 2 2  3 2  4 2  5 2  6 2  10 0  /
\endpicture} at -1 4

\put{
\beginpicture
 \setcoordinatesystem units <.15cm,.25cm>
 \multiput{$\circ$} at
   0 0  2 0  4 0  6 0  8 0  10 0
      1 1  3 1  5 1  7 1  9 1
         2 2  4 2  6 2  8 2
	    3 2  5 2  4 1  3 0 /
	     \setdots <.5mm>
	      \plot 0 0  1 1  2 0  3 1  4 0  5 1  6 0  7 1  8 0  9 1 /
 \plot  1 1  2 2  3 1  4 2  5 1  6 2  7 1  8 2  10 0 /
  \plot 3 0  5 2  5 1  /
   \plot 4 0  4 1  3 2  3 1  /
   \multiput{$\bullet$} at 1 1  3 2  4 2  5 2  6 2  9 1  /
\endpicture} at 1 4

\put{
\beginpicture
 \setcoordinatesystem units <.15cm,.25cm>
 \multiput{$\circ$} at
   0 0  2 0  4 0  6 0  8 0  10 0
      1 1  3 1  5 1  7 1  9 1
         2 2  4 2  6 2  8 2
	    3 2  5 2  4 1  3 0 /
	     \setdots <.5mm>
	      \plot 0 0  1 1  2 0  3 1  4 0  5 1  6 0  7 1  8 0  9 1 /
 \plot  1 1  2 2  3 1  4 2  5 1  6 2  7 1  8 2  10 0 /
  \plot 3 0  5 2  5 1  /
   \plot 4 0  4 1  3 2  3 1  /
   \multiput{$\bullet$} at 0 0  3 2  4 2  5 2  6 2  8 2  /
\endpicture} at 3 4

\put{
\beginpicture
 \setcoordinatesystem units <.15cm,.25cm>
 \multiput{$\circ$} at
   0 0  2 0  4 0  6 0  8 0  10 0
      1 1  3 1  5 1  7 1  9 1
         2 2  4 2  6 2  8 2
	    3 2  5 2  4 1  3 0 /
	     \setdots <.5mm>
	      \plot 0 0  1 1  2 0  3 1  4 0  5 1  6 0  7 1  8 0  9 1 /
 \plot  1 1  2 2  3 1  4 2  5 1  6 2  7 1  8 2  10 0 /
  \plot 3 0  5 2  5 1  /
   \plot 4 0  4 1  3 2  3 1  /
   \multiput{$\bullet$} at 2 2  3 2  4 2  5 2  6 2  9 1  /
\endpicture} at 0 5

\put{
\beginpicture
 \setcoordinatesystem units <.15cm,.25cm>
 \multiput{$\circ$} at
   0 0  2 0  4 0  6 0  8 0  10 0
      1 1  3 1  5 1  7 1  9 1
         2 2  4 2  6 2  8 2
	    3 2  5 2  4 1  3 0 /
	     \setdots <.5mm>
	      \plot 0 0  1 1  2 0  3 1  4 0  5 1  6 0  7 1  8 0  9 1 /
 \plot  1 1  2 2  3 1  4 2  5 1  6 2  7 1  8 2  10 0 /
  \plot 3 0  5 2  5 1  /
   \plot 4 0  4 1  3 2  3 1  /
   \multiput{$\bullet$} at 1 1  3 2  4 2  5 2  6 2  8 2  /
\endpicture} at 2 5

\put{
\beginpicture
 \setcoordinatesystem units <.15cm,.25cm>
 \multiput{$\circ$} at
   0 0  2 0  4 0  6 0  8 0  10 0
      1 1  3 1  5 1  7 1  9 1
         2 2  4 2  6 2  8 2
	    3 2  5 2  4 1  3 0 /
	     \setdots <.5mm>
	      \plot 0 0  1 1  2 0  3 1  4 0  5 1  6 0  7 1  8 0  9 1 /
 \plot  1 1  2 2  3 1  4 2  5 1  6 2  7 1  8 2  10 0 /
  \plot 3 0  5 2  5 1  /
   \plot 4 0  4 1  3 2  3 1  /
   \multiput{$\bullet$} at 2 2  3 2  4 2  5 2  6 2  8 2  /
\endpicture} at 1 6

\plot 0.3 0.3  0.7 0.7 /
\plot -0.3 0.3 -0.7 0.7 /
\plot 1.3 1.3  1.7 1.7 /
\plot -1.3 1.3 -1.7 1.7 /
\plot  0.3 1.7  0.7 1.3 /
\plot -0.3 1.7  -.7 1.3 /
\plot  1.7 2.3  1.3 2.7 /
\plot  -1.7 2.3  -1.3 2.7 /
\plot 0.3 2.3  0.7 2.7 /
\plot -.3 2.3  -.7 2.7 /
\plot 0.3 3.7  0.7 3.3 /
\plot -.3 3.7  -.7 3.3 /

\plot 1.3 2.3  1.7 2.7 /
\plot 0.7 2.3 0.3 2.7 /
\plot 2.3 3.3  2.7 3.7 /
\plot -.3 3.3 -.7 3.7 /
\plot  1.3 3.7  1.7 3.3 /
\plot  .7 3.7  .3 3.3 /
\plot  2.7 4.3  2.3 4.7 /
\plot  -.7 4.3  -.3 4.7 /
\plot 1.3 4.3  1.7 4.7 /
\plot .7 4.3  .3 4.7 /
\plot 1.3 5.7  1.7 5.3 /
\plot .7 5.7  .3 5.3 /

\plot 1 1.3  1 1.7 /
\plot 0 2.3  0 2.7  /
\plot 2 2.3  2 2.7 /
\plot -1 3.3 -1 3.7 /
\plot 1 3.3  1 3.7 /
\plot 0 4.3  0 4.7 /

\put{
\beginpicture
 \setcoordinatesystem units <.15cm,.25cm>
 \multiput{$\circ$} at
   0 0  2 0  4 0  6 0  8 0  10 0
      1 1  3 1  5 1  7 1  9 1
         2 2  4 2  6 2  8 2
	    3 2  5 2  4 1  3 0 /
	     \setdots <.5mm>
	      \plot 0 0  1 1  2 0  3 1  4 0  5 1  6 0  7 1  8 0  9 1 /
 \plot  1 1  2 2  3 1  4 2  5 1  6 2  7 1  8 2  10 0 /
  \plot 3 0  5 2  5 1  /
   \plot 4 0  4 1  3 2  3 1  /
  \multiput{$\bullet$} at 1 1  3 1  4 1  5 1  8 0  10 0 /
\endpicture} at 0 -2

\put{
\beginpicture
 \setcoordinatesystem units <.15cm,.25cm>
 \multiput{$\circ$} at
   0 0  2 0  4 0  6 0  8 0  10 0
      1 1  3 1  5 1  7 1  9 1
         2 2  4 2  6 2  8 2
	    3 2  5 2  4 1  3 0 /
	     \setdots <.5mm>
	      \plot 0 0  1 1  2 0  3 1  4 0  5 1  6 0  7 1  8 0  9 1 /
 \plot  1 1  2 2  3 1  4 2  5 1  6 2  7 1  8 2  10 0 /
  \plot 3 0  5 2  5 1  /
   \plot 4 0  4 1  3 2  3 1  /
  \multiput{$\bullet$} at 0 0  3 1  4 1  5 1  7 1  10 0  /
\endpicture} at 2 -2

\put{
\beginpicture
 \setcoordinatesystem units <.15cm,.25cm>
 \multiput{$\circ$} at
   0 0  2 0  4 0  6 0  8 0  10 0
      1 1  3 1  5 1  7 1  9 1
         2 2  4 2  6 2  8 2
	    3 2  5 2  4 1  3 0 /
	     \setdots <.5mm>
	      \plot 0 0  1 1  2 0  3 1  4 0  5 1  6 0  7 1  8 0  9 1 /
 \plot  1 1  2 2  3 1  4 2  5 1  6 2  7 1  8 2  10 0 /
  \plot 3 0  5 2  5 1  /
   \plot 4 0  4 1  3 2  3 1  /
   \multiput{$\bullet$} at 1 1  3 1  4 1  5 1  7 1  10 0  /
\endpicture} at 1 -1

\put{
\beginpicture
 \setcoordinatesystem units <.15cm,.25cm>
 \multiput{$\circ$} at
   0 0  2 0  4 0  6 0  8 0  10 0
      1 1  3 1  5 1  7 1  9 1
         2 2  4 2  6 2  8 2
	    3 2  5 2  4 1  3 0 /
	     \setdots <.5mm>
	      \plot 0 0  1 1  2 0  3 1  4 0  5 1  6 0  7 1  8 0  9 1 /
 \plot  1 1  2 2  3 1  4 2  5 1  6 2  7 1  8 2  10 0 /
  \plot 3 0  5 2  5 1  /
   \plot 4 0  4 1  3 2  3 1  /
  \multiput{$\bullet$} at 0 0  3 1  4 1  5 1  7 1  9 1 /

\endpicture} at 3 -1

\put{
\beginpicture
 \setcoordinatesystem units <.15cm,.25cm>
 \multiput{$\circ$} at
   0 0  2 0  4 0  6 0  8 0  10 0
      1 1  3 1  5 1  7 1  9 1
         2 2  4 2  6 2  8 2
	    3 2  5 2  4 1  3 0 /
	     \setdots <.5mm>
	      \plot 0 0  1 1  2 0  3 1  4 0  5 1  6 0  7 1  8 0  9 1 /
 \plot  1 1  2 2  3 1  4 2  5 1  6 2  7 1  8 2  10 0 /
  \plot 3 0  5 2  5 1  /
   \plot 4 0  4 1  3 2  3 1  /
  \multiput{$\bullet$} at 1 1  3 1  4 1  5 1  7 1  9 1 /
\endpicture} at 2 0 

\put{
\beginpicture
 \setcoordinatesystem units <.15cm,.25cm>
 \multiput{$\circ$} at
   0 0  2 0  4 0  6 0  8 0  10 0
      1 1  3 1  5 1  7 1  9 1
         2 2  4 2  6 2  8 2
	    3 2  5 2  4 1  3 0 /
	     \setdots <.5mm>
	      \plot 0 0  1 1  2 0  3 1  4 0  5 1  6 0  7 1  8 0  9 1 /
 \plot  1 1  2 2  3 1  4 2  5 1  6 2  7 1  8 2  10 0 /
  \plot 3 0  5 2  5 1  /
   \plot 4 0  4 1  3 2  3 1  /
  \multiput{$\bullet$} at 0 0  3 1  4 1  5 1  8 0  10 0 /
\endpicture} at 1 -3

\plot 0.3 -2.3  0.7 -2.7 /
\plot 1.3 -1.3  1.7 -1.7 /
\plot 2.3 -.3  2.7 -.7 /
\plot 1.3 -2.7  1.7 -2.3 /
\plot 2.3 -1.7  2.7 -1.3 /
\plot .3 -1.7  .7 -1.3 /
\plot 1.3 -0.7  1.7 -.3 /

\setdots <1mm>
\plot -1 1  0 -2 /
\plot 0 2  1 -1 /
\plot 1 3  2 0 /
\plot 0 0  1 -3 /
\plot 1 1  2 -2 /
\plot 2 2  3 -1 /

\setsolid
\plot -0.3 -2.3  -0.6  -2.6 /
\plot 3.3 -1.3  3.6  -1.6 /
\setshadegrid span <.6mm>
\vshade -3.5 6 6 <z,z,,> -2 4.5 7.5 <z,z,,> 0 6.5 9.5 <z,z,,> 1.5 8 8   /
\vshade -2.5 2 2  <z,z,,> -1.4 0.5 3 <z,z,,> -1.3 0.6 4 <z,z,,> 0 -0.5 5.5   
   <z,z,,> 1 0.5 6.5  <z,z,,> 2.3 1.8 5.2  <z,z,,> 2.4 2.8 5.1 <z,z,,> 3.5 4 4  /
\endpicture} at 0 0
\endpicture}
$$

\vfill\eject
It remains to present the rich antichains between $\Phi_4'$ and $\Phi_5$. They are
contained in $\Phi_{4,5}$:
$$
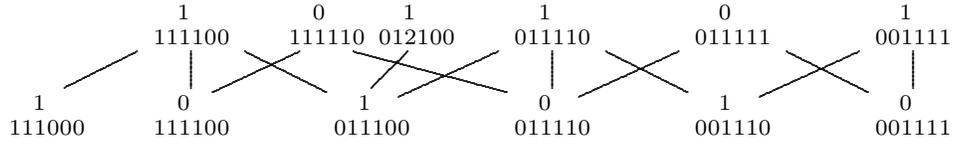
{\beginpicture
\setcoordinatesystem units <2cm,1cm>
\put{$1\ \,\atop111000$} at 0.2 0 
\put{$0\ \,\atop111100$} at 1 0
\put{$1\ \,\atop011100$} at 2 0
\put{$0\ \,\atop011110$} at 3 0
\put{$1\ \,\atop001110$} at 4 0
\put{$0\ \,\atop001111$} at 5 0

\put{$1\ \,\atop111100$} at 1 1 
\put{$0\ \,\atop111110$} at 1.75 1
\put{$1\ \,\atop012100$} at 2.25 1
\put{$1\ \,\atop011110$} at 3 1
\put{$0\ \,\atop011111$} at 4 1
\put{$1\ \,\atop001111$} at 5 1
\plot 0.3 0.3  0.7 0.7 /
\plot 1 0.3  1 0.7 /
\plot 2 0.3  2.2 0.7 /
\plot 3 0.3  3 0.7 /
\plot 5 0.3  5 0.7 /
\plot 1.12 .25  1.6 .7 /
\plot 1.3 .7  1.75 0.25 /
\plot 1.9 0.7  2.75 0.25 /
\plot 2.15 0.2  2.7 0.7 /
\plot 3.15 0.2  3.7 0.7 /
\plot 3.3 0.7  3.75 0.25 /
\plot 4.15 0.2  4.7 0.7 /
\plot 4.3 0.7  4.75 0.25 /
\endpicture}
$$
	\medskip
\noindent
Here are the rich antichains between $\Phi'_4$ and $\Phi_5$:
	\medskip
$$
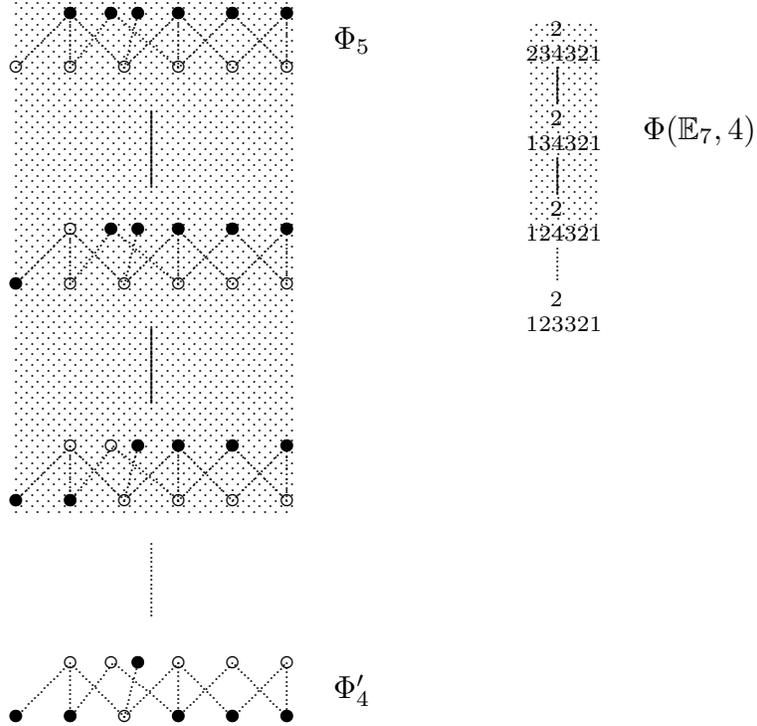
{\beginpicture
\setcoordinatesystem units <.6cm,.6cm>
\put{\beginpicture
\setcoordinatesystem units <1cm,1cm>
\put{$2\ \,\atop234321$} at 0 0
\put{$2\ \,\atop134321$} at 0 -1
\put{$2\ \,\atop124321$} at 0 -2
\put{$2\ \,\atop123321$} at 0 -3
\plot -.07 -.3  -.07 -.7 /
\plot -.07 -1.3  -.07 -1.7 /
\setdots <.5mm>
\plot -.07 -2.3  -.07 -2.7 /

\setshadegrid span <.6mm>
\vshade -.4 -2.2 0.2 <z,z,,> .4 -2.2 0.2  /
\put{$\Phi(\Bbb E_7,4)$} at 1.5 -1 

\endpicture} at 9 -2.5
\put{\beginpicture
\multiput{$\circ$} at 0 0  1 0  2 0  3 0  4 0  5 0
  1 1  1.75 1  2.25 1  3 1  4 1  5 1 /
\setdots <.5mm>
\plot 0 0  1 1  2 0  3 1  4 0  5 1  5 0  4 1  3 0  1.8 1  1 0  1 1  /
\plot 2 0  2.25 1 /
\plot 3 0  3 1 /
\multiput{$\bullet$} at 1 1  1.75 1  2.25 1  3 1  4 1  5 1 /
\endpicture} at 0 0
\put{$\Phi_5$} at 3.7 0
\put{\beginpicture
\multiput{$\circ$} at 0 0  1 0  2 0  3 0  4 0  5 0
  1 1  1.75 1  2.25 1  3 1  4 1  5 1 /
\setdots <.5mm>
\plot 0 0  1 1  2 0  3 1  4 0  5 1  5 0  4 1  3 0  1.8 1  1 0  1 1  /
\plot 2 0  2.25 1 /
\plot 3 0  3 1 /
\multiput{$\bullet$} at 0 0  1.75 1  2.25 1  3 1  4 1  5 1 /
\endpicture} at 0 -4
\put{\beginpicture
\multiput{$\circ$} at 0 0  1 0  2 0  3 0  4 0  5 0
  1 1  1.75 1  2.25 1  3 1  4 1  5 1 /
\setdots <.5mm>
\plot 0 0  1 1  2 0  3 1  4 0  5 1  5 0  4 1  3 0  1.8 1  1 0  1 1  /
\plot 2 0  2.25 1 /
\plot 3 0  3 1 /
\multiput{$\bullet$} at 0 0  1 0  2.25 1  3 1  4 1  5 1 /
\endpicture} at 0 -8
\put{\beginpicture
\multiput{$\circ$} at 0 0  1 0  2 0  3 0  4 0  5 0
  1 1  1.75 1  2.25 1  3 1  4 1  5 1 /
\setdots <.5mm>
\plot 0 0  1 1  2 0  3 1  4 0  5 1  5 0  4 1  3 0  1.8 1  1 0  1 1  /
\plot 2 0  2.25 1 /
\plot 3 0  3 1 /
\multiput{$\bullet$} at 0 0  1 0  2.25 1  3 0  4 0  5 0 /
\endpicture} at 0 -12
\plot 0 -1.3  0 -2.7 /
\plot 0 -5.3  0 -6.7 /
\setdots <.5mm>
\plot 0 -9.3  0 -10.7 /
\put{$\Phi_4'$} at 3.7 -12 
\setshadegrid span <.6mm>
\vshade -2.55 -8.8 .8 <z,z,,> 2.6 -8.8 .8  /
\endpicture}
$$
	\medskip 
Altogether we see that there is a poset isomorphism between $\Phi_+(\Bbb E_7)$ and 
$\Cal R(\Bbb E_7).$

\vfill\eject
{\bf Case} $\Bbb E_8.$
Let $\Phi_+ = \Phi_*(\Bbb E_8)$.
Let us exhibit a solid subposet of $\Cal R(\Bbb E_8)$ consisting of rich antichains in
$\Phi_{67}$:
$$
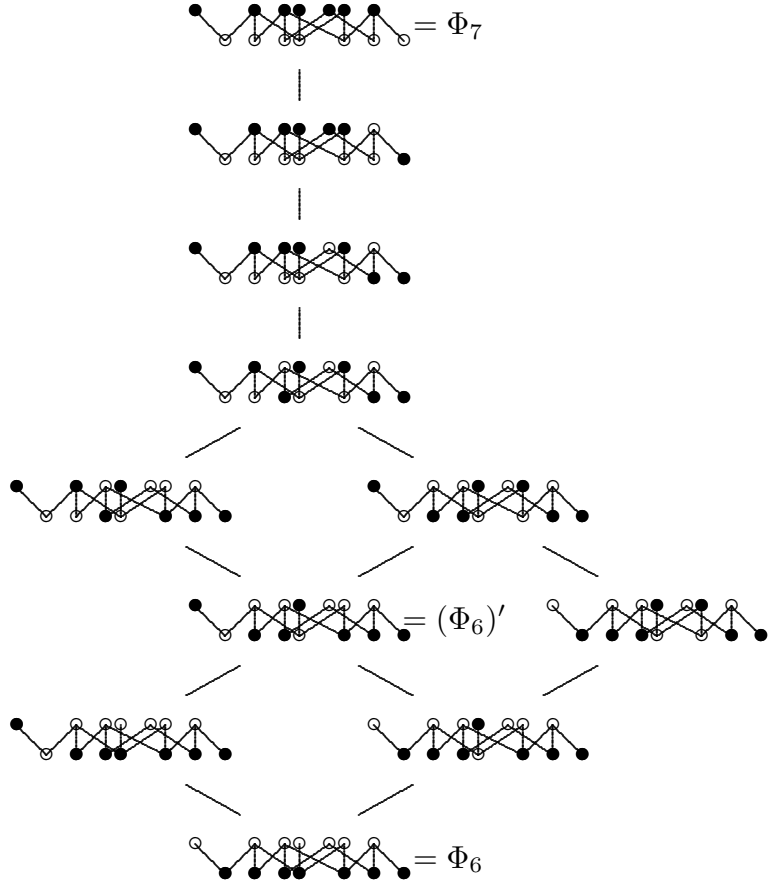
{\beginpicture
\setcoordinatesystem units <.33cm,.33cm>
\put{\beginpicture
\multiput{$\circ$} at 0 1  1 0  2 1  3.5 0  5 1
    2 0  3 1  5 0  6 1  7 0    3 0  3.5 1  4.5 1  6 0
        /
	\plot  0 1  1 0  2 1  3.5 0  5 1  5 0 /
	\plot  2 1  2 0  3 1   3 0   4.5 1  6 0  6 1 /
	\plot 3 1  5 0  6 1  7 0 /  
	\plot 3.5 0  3.5 1 /
	\multiput{$\bullet$} at 0 1  2 1  3 1  3.5 1  4.5 1  5 1  6 1 /
	\endpicture} at 0 0
	\put{$=\Phi_7$} at 5 0
\put{\beginpicture
	\multiput{$\circ$} at 0 1  1 0  2 1  3.5 0  5 1
	    2 0  3 1  5 0  6 1  7 0    3 0  3.5 1  4.5 1  6 0
	        /
		\plot  0 1  1 0  2 1  3.5 0  5 1  5 0 /
		\plot  2 1  2 0  3 1   3 0   4.5 1  6 0  6 1 /
		\plot 3 1  5 0  6 1  7 0 /
		\plot 3.5 0  3.5 1 /
		\multiput{$\bullet$} at 0 1  2 1  3 1  3.5 1  4.5 1  5 1  7 0 /
		\endpicture} at 0 -4
\put{\beginpicture
		\multiput{$\circ$} at 0 1  1 0  2 1  3.5 0  5 1
		    2 0  3 1  5 0  6 1  7 0    3 0  3.5 1  4.5 1  6 0
		        /
			\plot  0 1  1 0  2 1  3.5 0  5 1  5 0 /
			\plot  2 1  2 0  3 1   3 0   4.5 1  6 0  6 1 /
			\plot 3 1  5 0  6 1  7 0 /
			\plot 3.5 0  3.5 1 /
			\multiput{$\bullet$} at 0 1  2 1  3 1  3.5 1    5 1  6 0  7 0 /
			\endpicture} at 0 -8
\put{\beginpicture
			\multiput{$\circ$} at 0 1  1 0  2 1  3.5 0  5 1
			    2 0  3 1  5 0  6 1  7 0    3 0  3.5 1  4.5 1  6 0
			        /
				\plot  0 1  1 0  2 1  3.5 0  5 1  5 0 /
				\plot  2 1  2 0  3 1   3 0   4.5 1  6 0  6 1 /
				\plot 3 1  5 0  6 1  7 0 /
				\plot 3.5 0  3.5 1 /
				\multiput{$\bullet$} at 0 1  2 1  3 0  3.5 1    5 1  6 0  7 0 /
				\endpicture} at 0 -12
\put{\beginpicture
				\multiput{$\circ$} at 0 1  1 0  2 1  3.5 0  5 1
				    2 0  3 1  5 0  6 1  7 0    3 0  3.5 1  4.5 1  6 0
				        /
					\plot  0 1  1 0  2 1  3.5 0  5 1  5 0 /
					\plot  2 1  2 0  3 1   3 0   4.5 1  6 0  6 1 /
					\plot 3 1  5 0  6 1  7 0 /
					\plot 3.5 0  3.5 1 /
					\multiput{$\bullet$} at 0 1  2 1  3 0  3.5 1    5 0  6 0  7 0 /
					\endpicture} at -6 -16
\put{\beginpicture
					\multiput{$\circ$} at 0 1  1 0  2 1  3.5 0  5 1
					    2 0  3 1  5 0  6 1  7 0    3 0  3.5 1  4.5 1  6 0
					        /
						\plot  0 1  1 0  2 1  3.5 0  5 1  5 0 /
						\plot  2 1  2 0  3 1   3 0   4.5 1  6 0  6 1 /
						\plot 3 1  5 0  6 1  7 0 /
						\plot 3.5 0  3.5 1 /
						\multiput{$\bullet$} at 0 1  2 0  3 0  3.5 1    5 1  6 0  7 0 /
						\endpicture} at 6 -16
\put{\beginpicture
						\multiput{$\circ$} at 0 1  1 0  2 1  3.5 0  5 1
						    2 0  3 1  5 0  6 1  7 0    3 0  3.5 1  4.5 1  6 0
						        /
							\plot  0 1  1 0  2 1  3.5 0  5 1  5 0 /
							\plot  2 1  2 0  3 1   3 0   4.5 1  6 0  6 1 /
							\plot 3 1  5 0  6 1  7 0 /
							\plot 3.5 0  3.5 1 /
							\multiput{$\bullet$} at 0 1  2 0  3 0  3.5 1    5 0  6 0  7 0 /
							\endpicture} at 0 -20
							\put{$=(\Phi_6)'$} at 5.2 -20

\put{\beginpicture
\multiput{$\circ$} at 0 1  1 0  2 1  3.5 0  5 1
    2 0  3 1  5 0  6 1  7 0    3 0  3.5 1  4.5 1  6 0
        /
	\plot  0 1  1 0  2 1  3.5 0  5 1  5 0 /
	\plot  2 1  2 0  3 1   3 0   4.5 1  6 0  6 1 /
	\plot 3 1  5 0  6 1  7 0 /
	\plot 3.5 0  3.5 1 /
	\multiput{$\bullet$} at 1 0  2 0  3 0  3.5 1    5 1  6 0  7 0 /
	\endpicture} at 12 -20

\put{\beginpicture
\multiput{$\circ$} at 0 1  1 0  2 1  3.5 0  5 1
    2 0  3 1  5 0  6 1  7 0    3 0  3.5 1  4.5 1  6 0
        /
	\plot  0 1  1 0  2 1  3.5 0  5 1  5 0 /
	\plot  2 1  2 0  3 1   3 0   4.5 1  6 0  6 1 /
	\plot 3 1  5 0  6 1  7 0 /
	\plot 3.5 0  3.5 1 /
	\multiput{$\bullet$} at 0 1  2 0  3 0  3.5 0    5 0  6 0  7 0 /
	\endpicture} at -6 -24
	\put{\beginpicture
	\multiput{$\circ$} at 0 1  1 0  2 1  3.5 0  5 1
	    2 0  3 1  5 0  6 1  7 0    3 0  3.5 1  4.5 1  6 0
	        /
		\plot  0 1  1 0  2 1  3.5 0  5 1  5 0 /
		\plot  2 1  2 0  3 1   3 0   4.5 1  6 0  6 1 /
		\plot 3 1  5 0  6 1  7 0 /
		\plot 3.5 0  3.5 1 /
		\multiput{$\bullet$} at 1 0  2 0  3 0  3.5 1    5 0  6 0  7 0 /
		\endpicture} at 6 -24
\put{\beginpicture
		\multiput{$\circ$} at 0 1  1 0  2 1  3.5 0  5 1
		    2 0  3 1  5 0  6 1  7 0    3 0  3.5 1  4.5 1  6 0
		        /
			\plot  0 1  1 0  2 1  3.5 0  5 1  5 0 /
			\plot  2 1  2 0  3 1   3 0   4.5 1  6 0  6 1 /
			\plot 3 1  5 0  6 1  7 0 /
			\plot 3.5 0  3.5 1 /
			\multiput{$\bullet$} at 1 0  2 0  3 0  3.5 0    5 0  6 0  7 0 /
			\endpicture} at 0 -28
			\put{$=\Phi_6$} at 5 -28
			\plot 0 -1.5  0 -2.5  /
			\plot 0 -5.5  0 -6.5  /
			\plot 0 -9.5  0 -10.5  /
			\plot -2 -13.5  -3.8 -14.5 /
			\plot  2 -13.5   3.8 -14.5 /
			\plot -3.8 -17.5  -2 -18.5 /
			\plot  3.8 -17.5   2 -18.5 /
			\plot  8.2 -17.5   10 -18.5 /
			\plot -2 -21.5  -3.8 -22.5 /
			\plot  2 -21.5   3.8 -22.5 /
			\plot  10 -21.5  8.2 -22.5 /
			\plot -3.8 -25.5  -2 -26.5 /
			\plot  3.8 -25.5   2 -26.5 /
			\endpicture}
			$$
According to Theorem 1, $\Phi_7$ is the maximal element of $\Cal R(\Bbb E_8)$.
Thus, we see that the dual poset $(\Cal R(\Bbb E_8))^*$ has two different elements of height 5, whereas
$(\Phi_+)^*$ has just one element of height 5. This shows that $\Cal R(\Bbb E_8)$ is not
isomorphic to $\Phi_+$. 
        \bigskip
{\bf Case} $\Bbb F_4.$
     Let $\Phi_+ = \Phi_*(\Bbb F_4)$.
     Let us exhibit a solid subposet of $\Cal R(\Bbb F_4)$ consisting of rich antichains in
     $\Phi_{45}$:
     $$
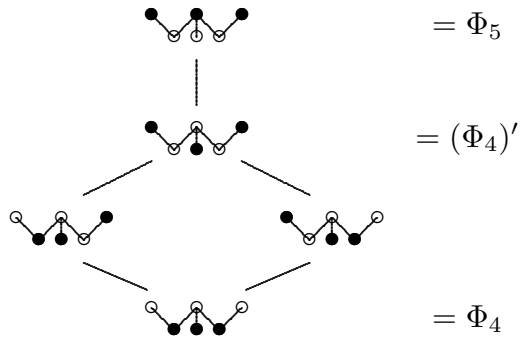
{\beginpicture
\setcoordinatesystem units <.25cm,.25cm>
\put{\beginpicture
\multiput{$\circ$} at 0 1  1 0  2 0  2 1  3 0 4 1 /
\plot 0 1  1 0  2 1  3 0  4 1 /
\plot 2 0  2 1 /
\multiput{$\bullet$} at 0 1  2 1  4 1 /
\endpicture} at 0 1
\put{$=\Phi_5$} at 12 1
\put{\beginpicture
\multiput{$\circ$} at 0 1  1 0  2 0  2 1  3 0 4 1 /
\plot 0 1  1 0  2 1  3 0  4 1 /
\plot 2 0  2 1 /
\multiput{$\bullet$} at 0 1  2 0  4 1 /
\endpicture} at 0 -4
\put{$=(\Phi_4)'$} at 12 -4
\put{\beginpicture
\multiput{$\circ$} at 0 1  1 0  2 0  2 1  3 0 4 1 /
\plot 0 1  1 0  2 1  3 0  4 1 /
\plot 2 0  2 1 /
\multiput{$\bullet$} at 1 0  2 0  4 1 /
\endpicture} at -6 -8
\put{\beginpicture
\multiput{$\circ$} at 0 1  1 0  2 0  2 1  3 0 4 1 /
\plot 0 1  1 0  2 1  3 0  4 1 /
\plot 2 0  2 1 /
\multiput{$\bullet$} at 0 1  2 0  3 0 /
\endpicture} at 6 -8
\put{\beginpicture
\multiput{$\circ$} at 0 1  1 0  2 0  2 1  3 0 4 1 /
\plot 0 1  1 0  2 1  3 0  4 1 /
\plot 2 0  2 1 /
\multiput{$\bullet$} at 1 0  2 0  3 0 /
\endpicture} at 0 -12
\put{$=\Phi_4$} at 12 -12
\plot 0 -.5  0 -2.5 /
\plot -2 -5  -5 -6.5 /
\plot 2 -5  5 -6.5 /
\plot -2.2 -10.7  -5 -9.5 /
\plot 2.2 -10.7  5 -9.5 /
\endpicture}
$$
According to Theorem 1, $\Phi_5$ is the maximal element of $\Cal R(\Bbb F_4)$.
Thus, we see that the dual poset $(\Cal R(\Bbb F_4))^*$ has two different elements of height 3, whereas
$(\Phi_+)^*$ has just one element of height 3. This shows that $\Cal R(\Bbb F_4)$ is not
isomorphic to $\Phi_+$. 
        
This completes the proof of Theorem 3. \hfill$\square$

\vfill\eject

{\gross Appendix: Some pictures of the root posets.}
	\medskip 
We exhibit a visualization of the root posets $\Bbb D_n, \Bbb E_6, \Bbb E_7,
\Bbb E_8$ and $\Bbb F_4$ which we found useful when 
preparing this note. Indeed, 
the pictures which we found in the literature and in the
web did not seem to be quite convincing (for us). 
The pictures shown below draw the attention to the fact
that for any pair of elements $x,z\in \Phi_+$, the
interval $\{y\in \Phi_+\mid x\le y \le z\}$ 
is a distributive
lattice which can be constructed in a convenient way using segments, squares and 
(3-dimensional) cubes.
A drawback of our visualization is that it does not take into
account the diversity of the various edges
(which otherwise could be indicated by using different slopes).
As a remedy, we label the edges by the corresponding basis
vectors (it is sufficient to do this at the boundary,
since the modularity transfers this information to the
remaining edges).

For the sake of completeness, we also include the cases
$\Bbb A_n, \Bbb B_n, \Bbb C_n$ and $\Bbb G_2$. It is well-known that the
root posets $\Bbb B_n$ and $\Bbb C_n$
do not differ as long as we do not refer
to the labels of the edges. But since we present the cases
 $\Bbb D_n, \Bbb E_m$ and $\Bbb F_4$ with labels, we do it also for
$\Bbb B_n$ and $\Bbb C_n$, thus we have to present
these types separately. 
      
The   cubical pictures
for  $\Bbb D_n, \Bbb E_m$ and $\Bbb F_4$ stress a division of the
positive roots into ``levels'' which seems to be a kind of measure
of the complexity of a positive root.
The roots belonging to a fixed
level form a planar graph, often a rectangle. 
All levels have a unique maximal element.
Level 1 
has $n-1$ minimal elements, all other levels have a unique 
minimal element. We list the
minimal and the maximal elements, as well as the number of roots
belonging to the level. A further analysis of the root posets 
will be provided in [R3]. 

	  \bigskip 
$$
\hbox{\beginpicture
\setcoordinatesystem units <.7cm,.55cm>
\multiput{} at 0 0  14 0 /
\put{\bf The root poset $\Bbb A_5$.} at -4 6 
\put{\beginpicture
\setcoordinatesystem units <1cm,1cm>
\multiput{$\circ$} at 0 0  1 0  2 0  3 0  4 0 /
\plot 0.1 0  0.9 0 /
\plot 1.1 0  1.9 0 /
\plot 2.1 0  2.9 0 /
\plot 3.1 0  3.9 0 /
\put{$a$\strut} at  0 -.3
\put{$b$\strut} at  1 -.3
\put{$c$\strut} at  2 -.3
\put{$d$\strut} at  3 -.3
\put{$e$\strut} at 4 -.3
\endpicture} at 9.5 5.5
\put{\beginpicture
\setcoordinatesystem units <.7cm,.65cm>
\put{$a$\strut} at  0 -.5
\put{$b$\strut} at  2 -.5
\put{$c$\strut} at  4 -.5
\put{$d$\strut} at  6 -.5
\put{$e$\strut} at 8 -.5
\multiput{$\ssize a$\strut} at 4.7 3.8   /
\multiput{$\ssize b$\strut} at 0.3 0.8  5.7 2.8 /
\multiput{$\ssize c$\strut} at  1.3  1.8  6.7  1.8 /
\multiput{$\ssize d$\strut} at  2.3 2.8  7.7 0.8 /
\multiput{$\ssize e$\strut} at  3.3  3.8 /
\multiput{$\bullet$} at  
  0 0  2 0  4 0  6 0  8 0 
  1 1  3 1  5 1  7 1
  2 2  4 2  6 2
  3 3  5 3
  4 4 /
\plot 0 0  4 4  8 0 /
\plot 1 1  2 0  5 3 /
\plot 2 2  4 0  6 2 /
\plot 3 3  6 0  7 1 /    
\endpicture} at 2 2
\endpicture}
$$
	\bigskip 
$$
\hbox{\beginpicture
\setcoordinatesystem units <.7cm,.55cm>
\multiput{} at 0 0  14 0 /
\put{\bf The root poset $\Bbb A_6$.} at -4 6 
\put{\beginpicture
\setcoordinatesystem units <1cm,1cm>
\multiput{$\circ$} at 0 0  1 0  2 0  3 0  4 0  5 0 /
\plot 0.1 0  0.9 0 /
\plot 1.1 0  1.9 0 /
\plot 2.1 0  2.9 0 /
\plot 3.1 0  3.9 0 /
\plot 4.1 0  4.9 0 /
\put{$a$\strut} at  0 -.3
\put{$b$\strut} at  1 -.3
\put{$c$\strut} at  2 -.3
\put{$d$\strut} at  3 -.3
\put{$e$\strut} at 4 -.3
\put{$f$\strut} at 5 -.3
\endpicture} at 9.5 5.5

\put{\beginpicture
\setcoordinatesystem units <.7cm,.65cm>
\put{$a$\strut} at  0 -.5
\put{$b$\strut} at  2 -.5
\put{$c$\strut} at  4 -.5
\put{$d$\strut} at  6 -.5
\put{$e$\strut} at 8 -.5
\put{$f$\strut} at 10 -.5
\multiput{$\ssize a$\strut} at 5.7 4.8   /
\multiput{$\ssize b$\strut} at 0.3 0.8  6.7 3.8 /
\multiput{$\ssize c$\strut} at  1.3  1.8  7.7  2.8 /
\multiput{$\ssize d$\strut} at  2.3 2.8  8.7 1.8 /
\multiput{$\ssize e$\strut} at  3.3  3.8  9.7  .8 /
\multiput{$\ssize f$\strut} at  4.3  4.8 /

\multiput{$\bullet$} at 
  0 0  2 0  4 0  6 0  8 0  10 0
  1 1  3 1  5 1  7 1  9 1
  2 2  4 2  6 2  8 2
  3 3  5 3  7 3
  4 4  6 4 
  5 5 /
\plot 0 0  5 5   10 0 /
\plot 1 1  2 0  6 4 /
\plot 2 2  4 0  7 3 /
\plot 3 3  6 0  8 2 /
\plot 4 4  8 0  9 1 /    
\endpicture} at 2 1
\endpicture}
$$
\vfill\eject

$$
\hbox{\beginpicture
\setcoordinatesystem units <1cm,1cm>
\put{\bf The root poset $\Bbb B_5$.} at -1 5
\put{\bf The root poset $\Bbb C_5$.} at 5 5
\put{\beginpicture
\setcoordinatesystem units <.8cm,.6cm>
\multiput{$\circ$} at 0 0  -1 0  -2 0  -3 0  -4 0 /
\plot -0.1 0  -0.9 0 /
\plot -1.1 0  -1.9 0 /
\plot -2.1 0  -2.9 0 /
\plot -3.1 0.1  -3.9 0.1 /
\plot -3.1 -.1  -3.9 -.1 /
\plot -3.35 .25  -3.6 0  -3.35 -.25 /
\put{$g$\strut} at 0 -0.6
\put{$f$\strut} at -1 -0.6
\put{$e$\strut} at -2 -0.6
\put{$d$\strut} at -3 -0.6
\put{$c$\strut} at -4 -0.6
\endpicture} at 0 4
\put{\beginpicture
\setcoordinatesystem units <.8cm,.6cm>
\multiput{$\circ$} at 0 0  -1 0  -2 0  -3 0  -4 0 /
\plot -0.1 0  -0.9 0 /
\plot -1.1 0  -1.9 0 /
\plot -2.1 0  -2.9 0 /
\plot -3.1 0.1  -3.9 0.1 /
\plot -3.1 -.1  -3.9 -.1 /
\plot -3.65 .25  -3.4 0  -3.65 -.25 /
\put{$f$\strut} at 0 -0.6
\put{$e$\strut} at -1 -0.6
\put{$d$\strut} at -2 -0.6
\put{$c$\strut} at -3 -0.6
\put{$u$\strut} at -4 -0.6
\endpicture} at 6 4
\put{\beginpicture
\setcoordinatesystem units <.55cm,.55cm>
\multiput{$\bullet$} at  
   0 0  -2 0  -4 0  -6 0  -8 0  
      -1 1   -3 1  -5 1  -7 1  
        -2 2  -4 2  -6 2  -8 2
            -3 3  -5 3  -7 3
             -4 4  -6 4  -8 4
                -5 5  -7 5
                  -6 6  -8 6
                    -7 7
                       -8 8 /
\plot 0 0  -8 8 /
\plot -1 1  -2 0  -8 6  -7 7 /
\plot -2 2  -4 0  -8 4  -6 6 /
\plot -3 3  -6 0  -8 2  -5 5 /
\plot -4 4  -8 0 /
\put{$g\strut$} at 0 -.6
\put{$f\strut$} at -2 -.6
\put{$e\strut$} at -4 -.6
\put{$d\strut$} at -6 -.6
\put{$c\strut$} at -8 -.6
\multiput{$\ssize g$} at           -7.7 6.8 /
\multiput{$\ssize f$} at  -0.3 0.9  -7.7 4.8  -7.3 7.9 /
\multiput{$\ssize e$} at  -1.3 1.9  -7.7 2.8  -6.3 6.9 /
\multiput{$\ssize d$} at  -2.3 2.9  -7.7 0.8  -5.3 5.9 /
\multiput{$\ssize c$} at  -3.3 3.9           -4.3 4.9 /

\endpicture} at 0 1
\put{\beginpicture
\setcoordinatesystem units <.55cm,.55cm>
\multiput{$\bullet$} at  
   0 0  -2 0  -4 0  -6 0  -8 0  
      -1 1   -3 1  -5 1  -7 1  
        -2 2  -4 2  -6 2  -8 2
            -3 3  -5 3  -7 3
             -4 4  -6 4  -8 4
                -5 5  -7 5
                  -6 6  -8 6
                     -7 7
                       -8 8 /
\plot 0 0  -8 8 /
\plot -1 1  -2 0  -8 6  -7 7 /
\plot -2 2  -4 0  -8 4  -6 6 /
\plot -3 3  -6 0  -8 2  -5 5 /
\plot -4 4  -8 0 /
\put{$f\strut$} at 0 -.6
\put{$e\strut$} at -2 -.6
\put{$d\strut$} at -4 -.6
\put{$c\strut$} at -6 -.6
\put{$u\strut$} at -8 -.6
\multiput{$\ssize f$} at           -7.7 6.8  -7.3 7.9 /
\multiput{$\ssize e$} at  -0.3 0.9  -7.7 4.8  -6.3 6.9 /
\multiput{$\ssize d$} at  -1.3 1.9  -7.7 2.8  -5.3 5.9 /
\multiput{$\ssize c$} at  -2.3 2.9  -7.7 0.8  -4.3 4.9 /
\multiput{$\ssize u$} at  -3.3 3.9  /

\endpicture} at 6 1
\endpicture}
$$
\bigskip\bigskip

$$
\hbox{\beginpicture
\setcoordinatesystem units <1cm,1cm>
\put{\bf The root poset $\Bbb B_6$.} at -1 5
\put{\bf The root poset $\Bbb C_6$.} at 5 5

\put{\beginpicture
\setcoordinatesystem units <.8cm,.6cm>
\multiput{$\circ$} at 1 0  0 0  -1 0  -2 0  -3 0  -4 0 /
\plot .1 0  .9 0 /
\plot -.1 0  -.9 0 /
\plot -1.1 0  -1.9 0 /
\plot -2.1 0  -2.9 0 /
\plot -3.1 0.1  -3.9 0.1 /
\plot -3.1 -.1  -3.9 -.1 /
\plot -3.35 .25  -3.6 0  -3.35 -.25 /
\put{$h$\strut} at 1 -0.6
\put{$g$\strut} at 0 -0.6
\put{$f$\strut} at -1 -0.6
\put{$e$\strut} at -2 -0.6
\put{$d$\strut} at -3 -0.6
\put{$c$\strut} at -4 -0.6
\endpicture} at 0 4

\put{\beginpicture
\setcoordinatesystem units <.5cm,.5cm>
\multiput{$\bullet$} at  
   0 0  -2 0  -4 0  -6 0  -8 0  -10 0
      -1 1   -3 1  -5 1  -7 1  -9 1
        -2 2  -4 2  -6 2  -8 2  -10 2
            -3 3  -5 3  -7 3  -9 3
             -4 4  -6 4  -8 4  -10 4
                -5 5  -7 5  -9 5
                  -6 6  -8 6  -10 6
                     -7 7   -9 7
                       -8 8  -10 8
                           -9 9
                             -10 10 /
\plot 0 0  -10 10 /
\plot -1 1  -2 0  -10 8  -9 9 /
\plot -2 2  -4 0  -10 6  -8 8 /
\plot -3 3  -6 0  -10 4  -7 7 /
\plot -4 4  -8 0  -10 2  -6 6 /
\plot -5 5  -10 0 /
\put{$h\strut$} at 0 -.6
\put{$g\strut$} at -2 -.6
\put{$f\strut$} at -4 -.6
\put{$e\strut$} at -6 -.6
\put{$d\strut$} at -8 -.6
\put{$c\strut$} at -10 -.6
\multiput{$\ssize h$} at           -9.7 8.8  /
\multiput{$\ssize g$} at  -0.3 0.9  -9.7 6.8  -9.3 9.9 /
\multiput{$\ssize f$} at  -1.3 1.9  -9.7 4.8  -8.3 8.9 /
\multiput{$\ssize e$} at  -2.3 2.9  -9.7 2.8  -7.3 7.9 /
\multiput{$\ssize d$} at  -3.3 3.9  -9.7 0.8  -6.3 6.9 /
\multiput{$\ssize c$} at  -4.3 4.9           -5.3 5.9 /

\endpicture} at 0 0.5
\put{\beginpicture
\setcoordinatesystem units <.8cm,.6cm>
\multiput{$\circ$} at 1 0  0 0  -1 0  -2 0  -3 0  -4 0 /
\plot .1 0  .9 0 /
\plot -.1 0  -.9 0 /
\plot -1.1 0  -1.9 0 /
\plot -2.1 0  -2.9 0 /
\plot -3.1 0.1  -3.9 0.1 /
\plot -3.1 -.1  -3.9 -.1 /
\plot -3.65 .25  -3.4 0  -3.65 -.25 /
\put{$g$\strut} at 1 -0.6
\put{$f$\strut} at 0 -0.6
\put{$e$\strut} at -1 -0.6
\put{$d$\strut} at -2 -0.6
\put{$c$\strut} at -3 -0.6
\put{$u$\strut} at -4 -0.6
\endpicture} at 6 4

\put{\beginpicture
\setcoordinatesystem units <.5cm,.5cm>
\multiput{$\bullet$} at  
   0 0  -2 0  -4 0  -6 0  -8 0  -10 0
      -1 1   -3 1  -5 1  -7 1  -9 1
        -2 2  -4 2  -6 2  -8 2  -10 2
            -3 3  -5 3  -7 3  -9 3
             -4 4  -6 4  -8 4  -10 4
                -5 5  -7 5  -9 5
                  -6 6  -8 6  -10 6
                     -7 7   -9 7
                       -8 8  -10 8
                           -9 9
                             -10 10 /
\plot 0 0  -10 10 /
\plot -1 1  -2 0  -10 8  -9 9 /
\plot -2 2  -4 0  -10 6  -8 8 /
\plot -3 3  -6 0  -10 4  -7 7 /
\plot -4 4  -8 0  -10 2  -6 6 /
\plot -5 5  -10 0 /
\put{$g\strut$} at 0 -.6
\put{$f\strut$} at -2 -.6
\put{$e\strut$} at -4 -.6
\put{$d\strut$} at -6 -.6
\put{$c\strut$} at -8 -.6
\put{$u\strut$} at -10 -.6
\multiput{$\ssize g$} at           -9.7 8.8  -9.3 9.9 /
\multiput{$\ssize f$} at  -0.3 0.9  -9.7 6.8  -8.3 8.9 /
\multiput{$\ssize e$} at  -1.3 1.9  -9.7 4.8  -7.3 7.9 /
\multiput{$\ssize d$} at  -2.3 2.9  -9.7 2.8  -6.3 6.9 /
\multiput{$\ssize c$} at  -3.3 3.9  -9.7 0.8  -5.3 5.9  /
\multiput{$\ssize u$} at  -4.3 4.9          /

\endpicture} at 6 0.5
\endpicture}
$$

\vfill\eject

$$
\hbox{\beginpicture
\setcoordinatesystem units <.7cm,.55cm>
\put{\bf The root poset $\Bbb D_5.$} at -5 8.5 

\put{\beginpicture
\setcoordinatesystem units <.8cm,.65cm>
\multiput{$\circ$} at  -1 0  -2 0  -3 0  -4 0  
   -3 1 /
\plot -1.1 0  -1.9 0 /
\plot -2.1 0  -2.9 0 /
\plot -3.1 0  -3.9 0 /
\plot -3 0.1  -3 0.9 /
\put{$b$\strut} at -4 -.5 
\put{$c$\strut} at -3 -.5 
\put{$d$\strut} at -2 -.5 
\put{$e$\strut} at -1 -.5 
\put{$u$\strut} at -3.3 1.2 

\put{$u$\strut} at -3.4 1.2
\endpicture} at 4 8.5

\multiput{$\bullet$} at  
  0 0  -2 0  -4 0  -5 0  -6 0
  -1 1  -3 1  -4 1  -5 1
  -2 2  -3 2  -5 2
  -2 3  -4 3
  -3 4  -5 4
  -4 5
  -5 6 /
               
\plot 0 0  -2 2  -4 0  -5 1  -6 0 /
\plot -1 1  -2 0  -3 1  /
\plot -2 2  -2 3  -5 6 /
\plot -3 1  -3 2  -5 4 /
\plot -4 0  -4 1  -5 2 /
\plot -5 1  -5 2  -3 4 / 
\plot -5 4  -4 5 /
\plot -5 0  -2 3 /

\setdashes <1mm>
\plot -3 4  -3 3  -5 1 /
\plot -3 1  -4 2  -4 3 /  
\plot -2 2  -3 3 /

\multiput{$\circ$} at -3 3  -4 2 /

\setshadegrid span <.6mm>
\vshade -5 1 2  <,z,,> -4 0 1    <z,,,> -2 2 3 /

\put{$b\strut$} at -6 -.6
\put{$c\strut$} at -4 -.6
\put{$d\strut$} at -2 -.6
\put{$e\strut$} at 0 -.6
\put{$u\strut$} at -5 -.6
\multiput{$\ssize b$} at -2.3 3.8 /
\multiput{$\ssize c$} at -1.3 1.8  -3.3 4.8  -5.95 0.4  
   -4.95 .4 /
\multiput{$\ssize d$} at -0.3 0.8  -4.3 5.8  -4.7 2.8 /
\multiput{$\ssize e$} at                   -4.7 4.8 /
\multiput{$\ssize u$} at -1.7 2.5     /

\endpicture}
$$
\bigskip
$$
\hbox{\beginpicture
\setcoordinatesystem units <.7cm,.55cm>
\put{\bf The root poset $\Bbb D_6.$} at -6 10.5 

\put{\beginpicture
\setcoordinatesystem units <.8cm,.65cm>
\multiput{$\circ$} at 0 0  -1 0  -2 0  -3 0  -4 0  
   -3 1 /
\plot -0.1 0  -0.9 0 /
\plot -1.1 0  -1.9 0 /
\plot -2.1 0  -2.9 0 /
\plot -3.1 0  -3.9 0 /
\plot -3 0.1  -3 0.9 /
\put{$b$\strut} at -4 -.5 
\put{$c$\strut} at -3 -.5 
\put{$d$\strut} at -2 -.5 
\put{$e$\strut} at -1 -.5 
\put{$f$\strut} at -0 -.5
\put{$u$\strut} at -3.3 1.2
\endpicture} at 3 10.5

\multiput{$\bullet$} at 0 0  -2 0  -4 0  -6 0  -8 0  
     -1 1  -3 1  -5 1  -7 1
     -2 2  -4 2 
     -3 3 
      -6 1
     -5 2  -7 2
     -4 3  -6 3   
     -3 4  -5 4  -7 4 
     -4 5  -6 5 
     -5 6  -7 6 
      -6 7 
     -7 8 
     -7 0 /
               
\plot 0 0  -3 3  -6 0  -7 1  -8 0 /
\plot -1 1  -2 0  -4 2 /
\plot -2 2  -4 0  -5 1 /

 \plot -7 0  -3 4  -7 8 /
\plot -6 1  -7 2  -4 5 /
\plot -5 2  -7 4  -5 6 /
\plot -4 3  -7 6  -6 7 /
\plot -3 3  -3 4 /
\plot -4 2  -4 3 /
\plot -5 1  -5 2 /
\plot -6 0  -6 1 /
\plot -7 1  -7 2 /

\setdashes <1mm>
\plot  -3 3  -4 4  -7 1 /
\plot -4 4  -4 5 /
\plot -4 2  -5 3  -5 4 /
\plot -5 1  -6 2  -6 3 /
   
\multiput{$\circ$} at -4 4  -5 3  -6 2 /
\setshadegrid span <.6mm>
\vshade   -7 1 2  <,z,,> -6 0 1  <z,,,> -3 3 4 /

\put{$b\strut$} at -8 -.6
\put{$c\strut$} at -6 -.6
\put{$d\strut$} at -4 -.6
\put{$e\strut$} at -2 -.6
\put{$f\strut$} at 0 -.6
\put{$u\strut$} at -7 -.6
\multiput{$\ssize b$} at -3.3 4.8 /
\multiput{$\ssize c$} at -2.3 2.8  -4.3 5.8 
   -7.95 0.4  
   -6.95 .4 /
\multiput{$\ssize d$} at  -1.3 1.8  -5.3 6.8  -6.7 2.8 /
\multiput{$\ssize e$} at  -0.3 0.8  -6.3 7.8  -6.7 4.8 /
\multiput{$\ssize f$} at                    -6.7 6.8 /
\multiput{$\ssize u$} at  -2.7 3.5  /
\endpicture}
$$ 
\bigskip\bigskip 
{\bf The levels for $\Bbb D_n$.}

$$
\hbox{\beginpicture
\setcoordinatesystem units <.95cm,.7cm>
\multiput{} at 0 1  13 0 /
\put{\bf Level} at 1 1
\put{ 1} at 1 0
\put{ 2} at 1 -1
\put{Conditions}  at 3.5 1
\put{$u = 0$}  at 3.5 0
\put{$u = 1$}  at 3.5 -1
\put{Minimal elements}  at 7 1
\put{Maximal element}  at 10.7 1 
\put{$n\!-\!1$ simple roots} at 7 0
\put{$\smallmatrix & 0  \cr
                   1&1 & \cdots & 1 & 1 & 1 \endsmallmatrix$} at 10.7 0 
\put{$\smallmatrix  & 1  \cr
                   0& 0 & \cdots & 0 & 0 & 0 \endsmallmatrix$} at 7  -1 
\put{$\smallmatrix  & 1  \cr
                   1& 2 & \cdots & 2 & 2 & 1 \endsmallmatrix$} at 10.7  -1 
\plot 0.6 0.5  14.5 0.5 /
\plot 2 1.3  2 -1.4 /
\plot 5 1.3  5 -1.4 /
\plot 8.8 1.3  8.8 -1.4 /
\plot 12.5 1.3  12.5 -1.4 /

\put{number}  at 13.5 1
\put{$\binom n2$} [l] at 13 0
\put{$\binom n2$} [l] at 13 -1

\endpicture}
$$                           
\vfill\eject
Unfortunately, the cubical pictures hides the symmetry given by 
the automorphism of order 2. Thus, one may want to ``squeeze''  the picture slightly. 
Here is the case $n=6$:
$$
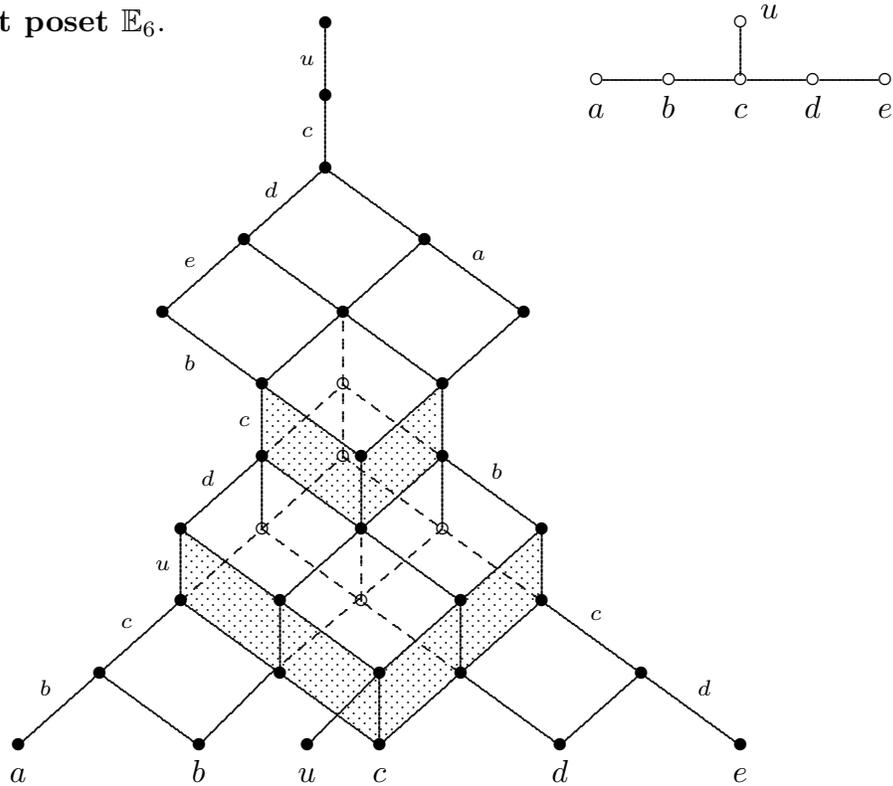
{\beginpicture
\setcoordinatesystem units <.6cm,.6cm>
\multiput{$\bullet$} at 0 0  -2 0  -4 0  -6 0  
   -1 1  -3 1  -5 1  
   -2 2  -4 2 
   -3 3 
   -5 5  -6 6  -7 7  -8 8
   -6 4  -7 5  -8 6
   -7 3  -8 4
   -8 2 
    -3.8 4  -4,8 3  -5.8 2  -6.8 1  -7.8 0 /
\plot 0 0  -3 3  -6 0  -6.8 1  -8 2  /
\plot -1 1  -2 0  -4 2  -4.8 3  -6 4  -8 6  -7 5    /
\plot -2 2  -4 0  -5 1  -5.8 2  -7 3  -8 4  -6 6 /
\plot -3 3  -3.8 4  -5 5  -8 8 /
\plot -5 5  -8 2 /
\plot -7.8 0  -3.8 4 /
\plot -7 7  -8 6 /

\setdashes <1mm>
\plot -3 3  -4.2 4  -5 5 /
\plot -4 2  -5.2 3  -6 4 /
\plot -5 1  -6.2 2  -7 3 /
\plot -6 0  -7.2 1  -8 2 /
\plot -8.2 0  -4.2 4 /
\multiput{$\circ$} at -8.2 0  -7.2 1  -6.2 2  -5.2 3  -4.2 4 /
\endpicture}
$$
\medskip
\hrule 
$$
\hbox{\beginpicture
\setcoordinatesystem units <1cm,.8cm>
\put{\bf The root poset $\Bbb E_6.$} at 0 10
\put{\beginpicture
\setcoordinatesystem units <.8cm,.65cm>
\multiput{$\circ$} at 0 0  1 0  2 0  3 0  4 0  
   2 1 /
\plot 0.1 0  0.9 0 /
\plot 1.1 0  1.9 0 /
\plot 2.1 0  2.9 0 /
\plot 3.1 0  3.9 0 /
\plot 2 0.1  2 0.9 /
\put{$a$\strut} at 0 -.5 
\put{$b$\strut} at 1 -.5 
\put{$c$\strut} at 2 -.5 
\put{$d$\strut} at 3 -.5 
\put{$e$\strut} at 4 -.5
\put{$u$\strut} at 2.4 1.2
\endpicture} at 8 9.5

\multiput{$\circ$} at 3.8 2
       2.7 3  4.7 3
       3.6 4  
       3.6 5 /

\multiput{$\bullet$} at
  0 0  2 0  4 0  6 0  8 0
  0.9 1  2.9 1  4.9 1  6.9 1
  1.8 2  5.8 2  / 
\plot 0 0  1.8 2  /
\plot 5.8 2   8 0 /
\plot 0.9 1  2 0  2.8 1 /
\plot 1.8 2  4 0  5.8 2 /
\plot 4.9 1  6 0 6.9 1 /

\setdashes <1mm>
\plot 1.8 2  3.6 4  5.8 2 / 
\plot 2.8 1  4.7 3 /
\plot  2.7 3  4.9 1 /
\plot 3.8 2   3.8 3  /
\plot 3.6 4  3.6 5 /
\plot 2.7 4   3.6 5  4.7 4 /
\plot 3.6 5  3.6 6  /

\setsolid 
\multiput{$\bullet$} at 4 1
  2.9 2  4.9 2 
  1.8 3  3.8 3  5.8 3
  2.7 4  4.7 4  / %
\plot 4 0  4 1 /
\plot   2.9 1   2.9 2 /
\plot 4.9 1  4.9 2 /
\plot  1.8 2    1.8 3  /
\plot 5.8 2  5.8 3 /
\plot  2.7 3    2.7 4  /
\plot 4.7 3  4.7 4 /
\plot   4 1  5.8 3  4.7 4 /
\plot 2.7 4  
    1.8 3  4 1 / 
\plot 2.9 2  4.7 4 /
\plot  4.9 2   2.7 4  /

\multiput{$\bullet$} at 3.8 4
    2.7 5  4.7 5
    3.6 6 /
\plot  3.8 3  3.8 4 /
\plot 2.7 4  2.7 5  /
\plot 4.7 4  4.7 5 /
\plot 3.8 4
    2.7 5   3.6 6  4.7 5
   3.8 4 /

\setshadegrid span <.6mm>
\vshade 1.8 2 3  <,z,,> 4 0 1  <z,,,> 5.8 2 3 /
\vshade 2.7 4 5  <,z,,> 3.8 3 4  <z,,,> 4.7 4 5 /

\multiput{$\bullet$} at 1.6 6  5.6 6  2.5 7  4.5 7  3.4 8  / 
\plot 2.7 5  1.6 6  3.4 8  5.6 6  4.7 5 /
\plot 2.5 7  3.6 6  4.5 7 / 

\multiput{$\bullet$} at 3.4 9  3.4 10  / 
\plot 3.4 8  3.4 10 /
  
\put{$\bullet$} at 3.2 0 
\plot 3.2 0  4 1 /

\put{$a\strut$} at 0 -.4
\put{$b\strut$} at 2 -.4
\put{$c\strut$} at 4 -.4
\put{$d\strut$} at 6 -.4
\put{$e\strut$} at 8 -.4
\put{$u\strut$} at 3.2 -.4
\multiput{$\ssize a$} at  5.1 6.8 /
\multiput{$\ssize b$} at  0.3 0.8  1.9 5.3  5.3 3.8 /
\multiput{$\ssize c$} at  1.2 1.7  2.5 4.5  3.2 8.5
         6.4 1.8  /
\multiput{$\ssize d$} at  2.1 3.7  2.8 7.7
         7.6 0.8 /
\multiput{$\ssize e$} at  1.9 6.7 /
\multiput{$\ssize u$} at  1.6 2.5  3.2 9.5 /

\endpicture}
$$

\bigskip\bigskip
$$
\hbox{\beginpicture
\setcoordinatesystem units <.95cm,.7cm>
\multiput{} at 0 1  13 0 /
\put{\bf Level} at 1 1
\put{ 1} at 1 0
\put{ 2} at 1 -1
\put{ 3} at 1 -2
\put{Conditions}  at 3.5 1
\put{$u = 0$}  at 3.5 0
\put{$u = 1$}  at 3.5 -1
\put{$c\ge 2$}  at 3.5 -2
\put{Minimal elements}  at 7 1
\put{Maximal element}  at 10.7 1 
\put{5 simple roots} at 7 0
\put{$\smallmatrix   &        &  0  \cr
                   1 & 1 & 1 & 1 & 1 \endsmallmatrix$} at 10.7 0 
\put{$\smallmatrix   &   & 1  \cr
                   0 & 0 & 0 & 0 & 0 \endsmallmatrix$} at 7  -1 
\put{$\smallmatrix   &   & 1  \cr
                   1 & 1 & 1 & 1 & 1 \endsmallmatrix$} at 10.7  -1 
\put{$\smallmatrix   &   & 1  \cr
                   0 & 1 & 2 & 1 & 0 \endsmallmatrix$} at 7  -2 
\put{$\smallmatrix   &   & 1  \cr
                   1 & 2 & 3 & 2 & 1 \endsmallmatrix$} at 10.7  -2 
\plot 0.6 0.5  14.5 0.5 /
\plot 2 1.3  2 -2.4 /
\plot 5 1.3  5 -2.4 /
\plot 8.8 1.3  8.8 -2.4 /
\plot 12.5 1.3  12.5 -2.4 /

\put{number}  at 13.5 1
\put{$15$} [l] at 13 0
\put{$3\times3+1$} [l] at 13 -1
\put{$3\times3+2$} [l] at 13 -2

\endpicture}
$$

\vfill\eject
$$
\hbox{\beginpicture
\setcoordinatesystem units <1cm,.8cm>
\put{\bf The root poset $\Bbb E_7.$} at -2 15.5 
\put
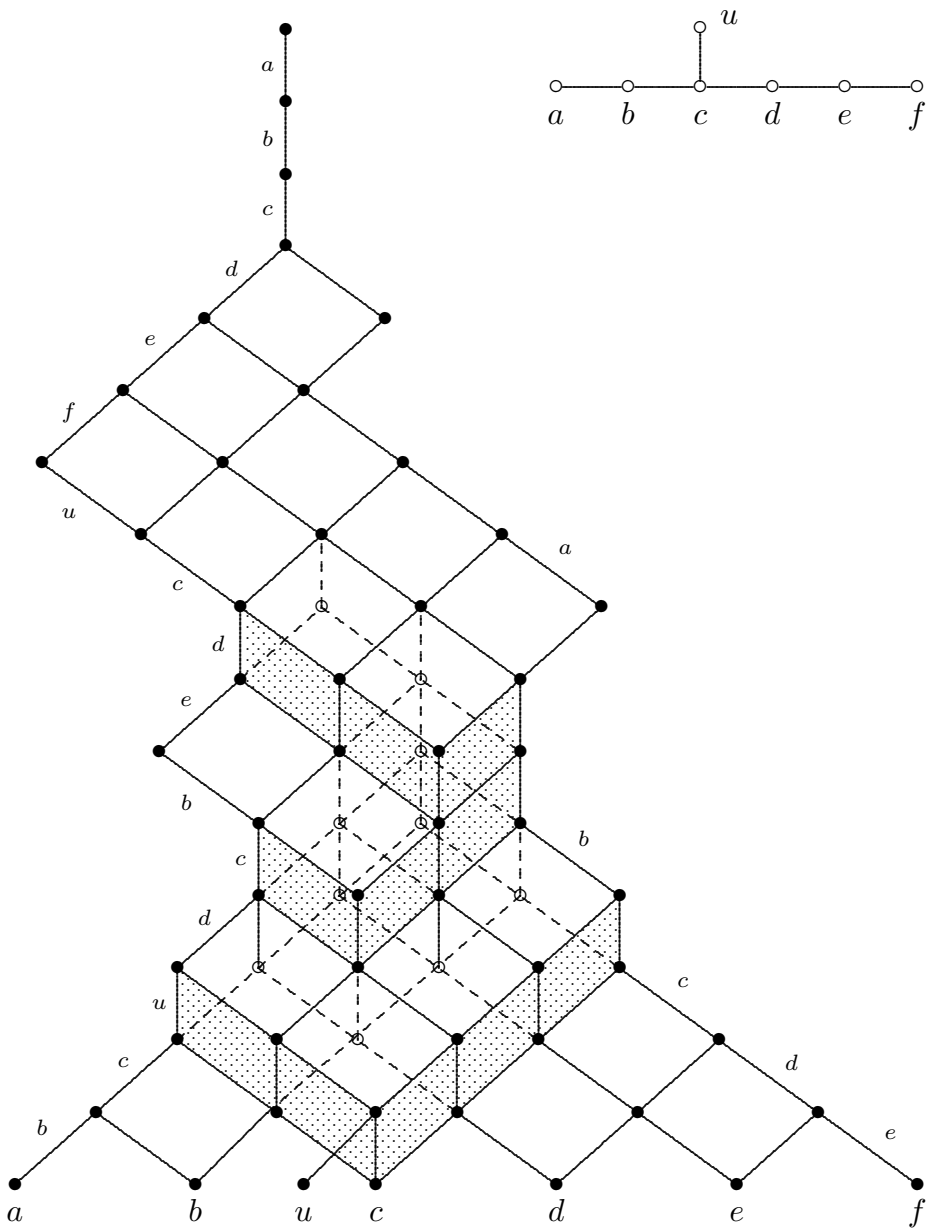
{\beginpicture
\setcoordinatesystem units <.8cm,.65cm>
\multiput{$\circ$} at 0 0  1 0  2 0  3 0  4 0  5 0 
   2 1 /
\plot 0.1 0  0.9 0 /
\plot 1.1 0  1.9 0 /
\plot 2.1 0  2.9 0 /
\plot 3.1 0  3.9 0 /
\plot 4.1 0  4.9 0 /
\plot 2 0.1  2 0.9 /
\put{$a$\strut} at 0 -.5 
\put{$b$\strut} at 1 -.5 
\put{$c$\strut} at 2 -.5 
\put{$d$\strut} at 3 -.5 
\put{$e$\strut} at 4 -.5 
\put{$f$\strut} at 5 -.5 
\put{$u$\strut} at 2.4 1.2
\endpicture} at 8 15.5

\multiput{$\circ$} at 3.8 2
       2.7 3  4.7 3
       3.6 4  5.6 4
       4.5 5   4.5 6 
       3.6 5 
       4.5 7  3.4 8 /

\multiput{$\bullet$} at
  0 0  2 0  4 0  6 0  8 0  10 0  
  0.9 1  2.9 1  4.9 1  6.9 1  8.9 1
  1.8 2  5.8 2  7.8 2
  6.7 3 /  
\plot 0 0  1.8 2  /
\plot 5.8 2   8 0  8.9 1 /
\plot 0.9 1  2 0  2.8 1 /
\plot 1.8 2  4 0  5.8 2 /
\plot 4.9 1  6 0  7.8 2 /
\plot 5.8 2  6.7 3  10 0 /

\setdashes <1mm>
\plot 1.8 2  3.6 4  5.8 2 / 
\plot 3.6 4  4.5 5  6.7 3 /
\plot 2.8 1  5.6  4 /
\plot  2.7 3  4.9 1 /
\plot 3.8 2   3.8 3  /
\plot 3.6 4  3.6 5 /
\plot 2.7 4   3.6 5  4.7 4 /
\plot 3.6 5  3.6 6  /
\plot 5.6 4  5.6 5 /
\plot 4.5 5  4.5 6 /
\plot 3.6 5  4.5 6  5.6 5 / 
\plot  4.5 6  4.5 7 /
\plot 2.5 7  3.4 8  5.6 6 /
\plot 3.6 6  4.5 7 / 
\plot 4.5 7  4.5 8 /

\plot 3.4 8  3.4 9 /

\setsolid 
\multiput{$\bullet$} at 4 1
  2.9 2  4.9 2 
  1.8 3  3.8 3  5.8 3
  2.7 4  4.7 4 
   6.7 4  5.6 5 / 
\plot 4 0  4 1 /
\plot 2.9 1   2.9 2 /
\plot 4.9 1  4.9 2 /
\plot 1.8 2    1.8 3  /
\plot 5.8 2  5.8 3 /
\plot 2.7 3    2.7 4  /
\plot 4.7 3  4.7 4 /
\plot 4 1  5.8 3  4.7 4 /
\plot 2.7 4  
    1.8 3  4 1 / 
\plot 2.9 2  4.7 4 /
\plot  4.9 2   2.7 4  /
\plot 6.7 3  6.7 4  5.8 3  /
\plot 6.7 4  5.6 5  4.7 4 /
\plot 5.6 5  5.6 6 /

\multiput{$\bullet$} at 3.8 4
    2.7 5  4.7 5
    3.6 6 /
\plot  3.8 3  3.8 4 /
\plot 2.7 4  2.7 5  /
\plot 4.7 4  4.7 5 /
\plot 3.8 4
    2.7 5   3.6 6  4.7 5
   3.8 4 /

\plot  5.6 6  4.7 5 /

\multiput{$\bullet$} at   0.3 10  1.4 9  2.5 8   3.6 7  4.7 6  
     1.2 11   2.3 10  3.4 9  4.5 8   5.6 7 
     2.1 12  3.2 11  4.3 10  5.4 9  6.5 8   
     3. 13  4.1 12
 /
\plot 2.1 12  3 13  4.1 12  3.2 11 /
\plot 0.3 10  4.7 6  6.5 8  2.1 12  0.3 10 /
\plot 1.2 11  5.6 7 /

\plot 1.4 9  3.2 11 /
\plot 2.5 8  4.3 10 /
\plot 3.6 7  5.4 9 /

\plot 4.7 5  4.7 6 /
\plot 3.6 6  3.6 7 /
\plot 5.6 6  5.6 7 /
\plot 2.5 7  2.5 8 /

\multiput{$\bullet$} at  3  14  3 15  3 16  /
\plot 3  13  3 16 /

\setshadegrid span <.6mm>
\vshade 1.8 2 3  <,z,,> 4 0 1  <z,,,> 6.7 3 4 /
\vshade 2.7 4 5  <,z,,> 3.8 3 4  <z,,,> 5.6  5 6 /
\vshade 2.5 7 8  <,z,,> 4.7 5 6  <z,,,> 5.6   6 7 /

\multiput{$\bullet$} at 1.6 6   5.6 6  2.5 7  /

\plot 2.7 5  1.6 6   2.5 7 /

\plot 2.5 7  3.6 6 /

\put{$\bullet$} at 3.2 0 
\plot 3.2 0  4 1 /

\put{$a\strut$} at 0 -.4
\put{$b\strut$} at 2 -.4
\put{$c\strut$} at 4 -.4
\put{$d\strut$} at 6 -.4
\put{$e\strut$} at 8 -.4
\put{$f\strut$} at 10 -.4
\put{$u\strut$} at 3.2 -.4
\multiput{$\ssize a$} at  6.1 8.8 
         2.8 15.5 /
\multiput{$\ssize b$} at  0.3 0.8  1.9 5.3  6.3 4.8
         2.8 14.5  /
\multiput{$\ssize c$} at  1.2 1.7  2.5 4.5  1.8 8.3
         2.8 13.5 
         7.4 2.8  /
\multiput{$\ssize d$} at  2.1 3.7  2.25 7.5  2.4 12.7 
         8.6 1.7 /
\multiput{$\ssize e$} at  1.9 6.7  1.5 11.7 
         9.7 .7 /
\multiput{$\ssize f$} at  0.6 10.7 /
\multiput{$\ssize u$} at  1.6 2.5  0.6 9.3 /

\endpicture}
$$

\bigskip\bigskip
$$
\hbox{\beginpicture
\setcoordinatesystem units <.95cm,.7cm>
\multiput{} at 0 1  13 0 /
\put{\bf Level} at 1 1
\put{ 1} at 1 0
\put{ 2} at 1 -1
\put{ 3} at 1 -2
\put{ 4} at 1 -3
\put{Conditions}  at 3.5 1
\put{$u = 0$}  at 3.5 0
\put{$u = 1,\ c\le 1$}  at 3.5 -1
\put{$c=2,\ d=1$}  at 3.5 -2
\put{$d\ge 2$}  at 3.5 -3
\put{Minimal elements}  at 7 1
\put{Maximal element}  at 10.7 1 
\put{6 simple roots} at 7 0
\put{$\smallmatrix   &        &  0  \cr
                   1 & 1 & 1 & 1 & 1 & 1 \endsmallmatrix$} at 10.7 0 
\put{$\smallmatrix   &   & 1  \cr
                   0 & 0 & 0 & 0 & 0 & 0 \endsmallmatrix$} at 7  -1 
\put{$\smallmatrix   &   & 1  \cr
                   1 & 1 & 1 & 1 & 1 & 1\endsmallmatrix$} at 10.7  -1 
\put{$\smallmatrix   &   & 1  \cr
                   0 & 1 & 2 & 1 & 0 & 0 \endsmallmatrix$} at 7  -2 
\put{$\smallmatrix   &   & 1  \cr
                   1 & 2 & 2 & 1 & 1 & 1 \endsmallmatrix$} at 10.7  -2 
\put{$\smallmatrix   &   & 1  \cr
                   0 & 1 & 2 & 2 & 1 & 0 \endsmallmatrix$} at 7  -3 
\put{$\smallmatrix   &   & 2  \cr
                   2 & 3 & 4 & 3 & 2 & 1 \endsmallmatrix$} at 10.7  -3 
\plot 0.6 0.5  14.5 0.5 /
\plot 2 1.3  2 -3.4 /
\plot 5 1.3  5 -3.4 /
\plot 8.8 1.3  8.8 -3.4 /
\plot 12.5 1.3  12.5 -3.4 /

\put{number}  at 13.5 1
\put{$21$} [l] at 13 0
\put{$3\times4+1$} [l] at 13 -1
\put{$3\times3$} [l] at 13 -2
\put{$5\times3+5$} [l] at 13 -3

\endpicture}
$$

\vfill\eject

	
$$
\hbox{\beginpicture
\setcoordinatesystem units <.8cm,.65cm>
\put{\bf The root poset $\Bbb E_8$.} at -2 27.8 

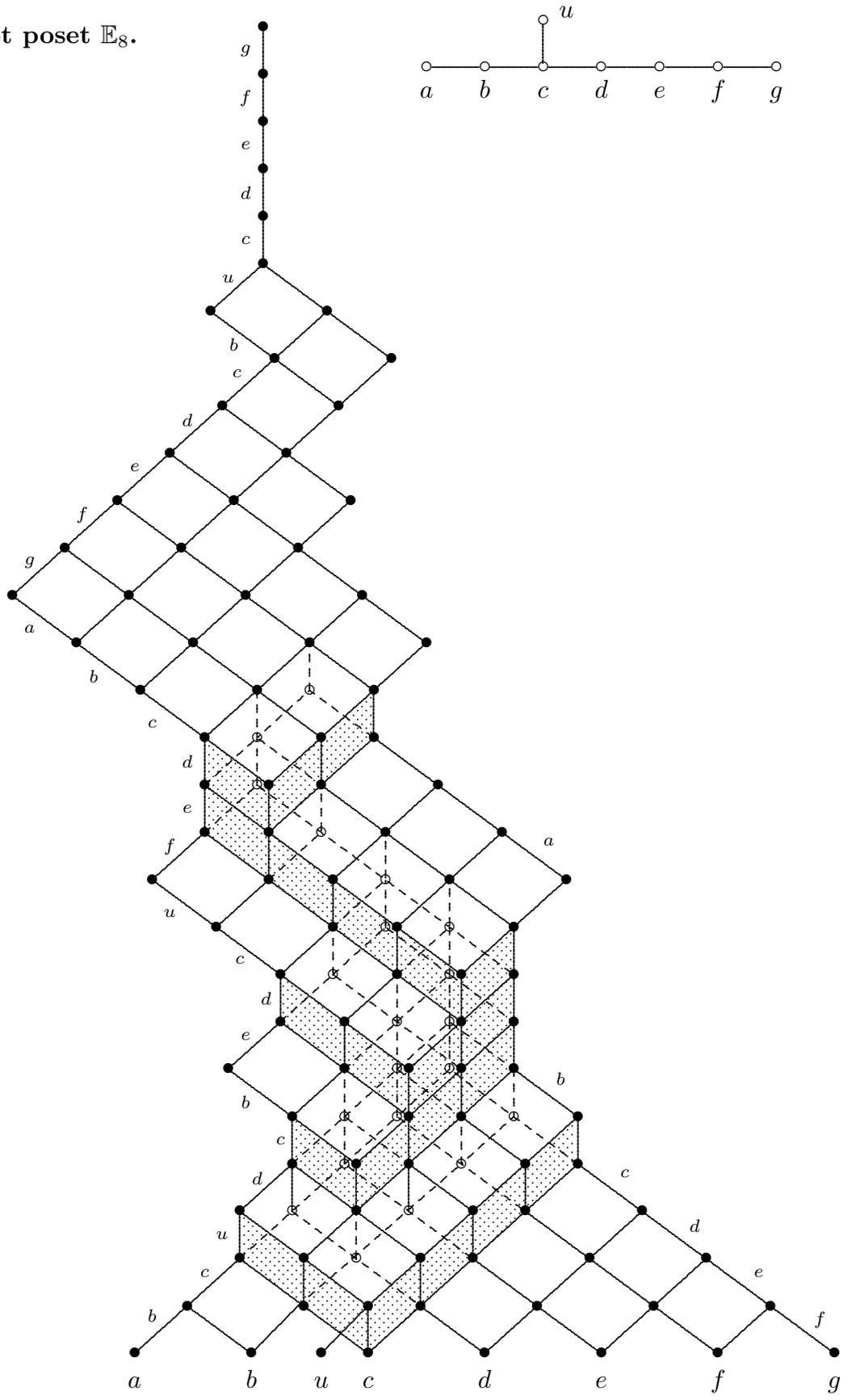
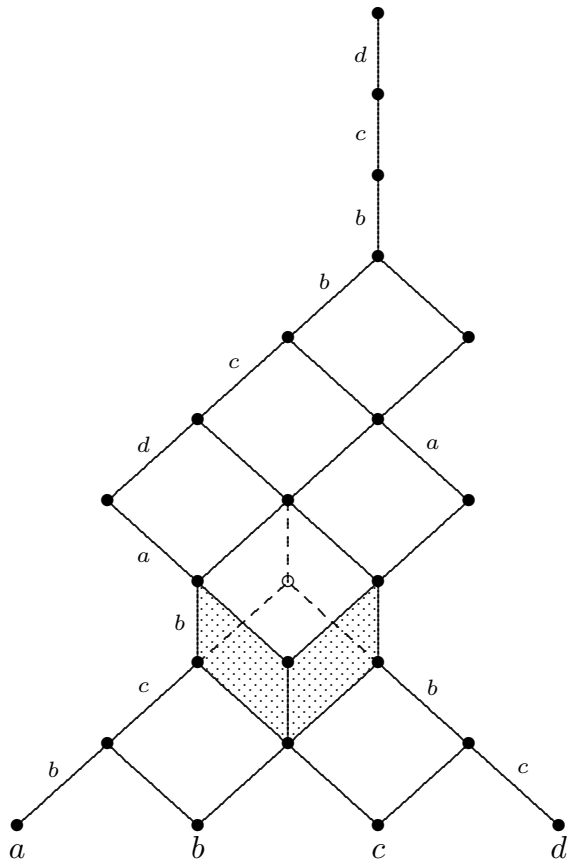
\put{\beginpicture
\multiput{$\circ$} at 0 0  1 0  2 0  3 0  4 0  5 0  6 0
   2 1 /
\plot 0.1 0  0.9 0 /
\plot 1.1 0  1.9 0 /
\plot 2.1 0  2.9 0 /
\plot 3.1 0  3.9 0 /
\plot 4.1 0  4.9 0 /
\plot 5.1 0  5.9 0 /
\plot 2 0.1  2 0.9 /
\put{$a$\strut} at 0 -.5 
\put{$b$\strut} at 1 -.5 
\put{$c$\strut} at 2 -.5 
\put{$d$\strut} at 3 -.5 
\put{$e$\strut} at 4 -.5 
\put{$f$\strut} at 5 -.5 
\put{$g$\strut} at 6 -.5 
\put{$u$\strut} at 2.4 1.2
\endpicture} at 8 27.5

\multiput{$\circ$} at 3.8 2
       2.7 3  4.7 3
       3.6 4  5.6 4
       4.5 5  6.5 5
       5.4 6  
       3.6 5  
       4.5 6   
       5.4 7 
       4.5 7  3.4 8 
        5.4 8 
       4.3 9 
      2.1 12  3.2 11  4.3 10  5.4 9   /

\multiput{$\bullet$} at
  0 0  2 0  4 0  6 0  8 0  10 0  
  0.9 1  2.9 1  4.9 1  6.9 1  8.9 1
  1.8 2  5.8 2  7.8 2
  6.7 3 
  7.6 4  8.7 3  9.8 2  10.9 1  12 0 
   7.6 5  6.5 6  6.5 7 
/  
\plot 0 0  1.8 2  /
\plot 5.8 2   8 0  9.8 2 /
\plot 0.9 1  2 0  2.8 1 /
\plot 1.8 2  4 0  5.8 2 /
\plot 4.9 1  6 0  8.7 3 /
\plot 5.8 2  6.7 3  10 0  10.9 1 /
\plot  7.6 4  12 0 /
\plot 6.7 3  7.6 4  7.6 5  6.5 6  5.6 5 /
\plot 6.5 7  5.6 6 /
\plot 6.5 6  6.5 8 /
\plot 6.7 4 7.6 5 /

\setdashes <1mm>
\plot 1.8 2  3.6 4  5.8 2 / 
\plot 3.6 4  4.5 5  6.7 3 /
\plot 2.8 1  5.6  4 /
\plot  2.7 3  4.9 1 /
\plot 3.8 2   3.8 3  /
\plot 3.6 4  3.6 5 /
\plot 2.7 4   3.6 5  4.7 4 /
\plot 3.6 5  3.6 6  /
\plot 5.6 4  5.6 5 /
\plot 4.5 5  4.5 6 /
\plot 3.6 5  4.5 6  5.6 5 / 
\plot  4.5 6  4.5 7 /
\plot 2.5 7  3.4 8  5.6 6 /
\plot 3.6 6  4.5 7 / 
\plot 4.5 7  4.5 8 /
\plot 3.4 8  3.4 9 /

\plot 7.6 4  5.4 6  4.5 5 /
\plot 6.5 6  6.5 5  5.6 4 /
\plot 4.5 6  5.4 7  6.5 6 /
\plot 5.4 7  5.4 6 /

\plot 5.4 7  5.4 9 /
\plot 4.5 7  5.4 8   6.5 7 /
\plot 3.4 8  4.3 9  5.4 8 /
\plot 4.3 9  4.3 10 /

\setsolid 
\multiput{$\bullet$} at 4 1
  2.9 2  4.9 2 
  1.8 3  3.8 3  5.8 3
  2.7 4  4.7 4 
   6.7 4  5.6 5 / 
\plot 4 0  4 1 /
\plot 2.9 1   2.9 2 /
\plot 4.9 1  4.9 2 /
\plot 1.8 2    1.8 3  /
\plot 5.8 2  5.8 3 /
\plot 2.7 3    2.7 4  /
\plot 4.7 3  4.7 4 /
\plot 4 1  5.8 3  4.7 4 /
\plot 2.7 4  
    1.8 3  4 1 / 
\plot 2.9 2  4.7 4 /
\plot  4.9 2   2.7 4  /
\plot 6.7 3  6.7 4  5.8 3  /
\plot 6.7 4  5.6 5  4.7 4 /
\plot 5.6 5  5.6 6 /

\multiput{$\bullet$} at 3.8 4
    2.7 5  4.7 5
    3.6 6 /
\plot  3.8 3  3.8 4 /
\plot 2.7 4  2.7 5  /
\plot 4.7 4  4.7 5 /
\plot 3.8 4
    2.7 5   3.6 6  4.7 5
   3.8 4 /

\plot  5.6 6  4.7 5 /

\multiput{$\bullet$} at   0.3 10  1.4 9  2.5 8   3.6 7  4.7 6  
     1.2 11   2.3 10  3.4 9  4.5 8   5.6 7  6.5 8  /
\plot   1.2 11 0.3 10  4.7 6  6.5 8 /

\plot 1.2 11  5.6 7 /

\plot 1.4 9  2.3 10 /
\plot 2.5 8  3.4 9  /
\plot 3.6 7  4.5 8  /

\plot 4.7 5  4.7 6 /
\plot 3.6 6  3.6 7 /
\plot 5.6 6  5.6 7 /
\plot 2.5 7  2.5 8 /
\setdashes <1mm>
\plot  1.2 11  2.1 12  6.5 8  /
\plot  2.3 10  3.2 11  3.2 12  /
\plot  3.4 9   4.3 10  4.3 11 /
\plot  4.5 8   5.4 9   5.4 10 /
\plot 2.1 12  2.1 13 /
\plot 1.2 12  2.1 13  3.2 12 /

\plot 3 14  3 15 /
\plot 2.1 14  2.1 13  3 14  4.1 13 /
\setsolid 

\multiput{$\circ$} at 2.1 13  3.0 14  /
\multiput{$\bullet$} at    
    1.2 12  2.3 11  3.4 10  4.5 9  5.6 8 
    3.2 12  4.3 11  5.4 10  6.5 9 
    4.1 13  5.2 12  6.3 11  7.4 10 /
\plot 1.2 13  1.2 12  5.6 8  7.4 10  4.1 13  2.3 11  2.3 12  /
\plot 3.2 13  3.2 12  6.5 9 /
\plot  3.4 10  5.2 12 /
\plot 4.5 9  6.3 11 /

\plot 1.2 12  1.2 11 /
\plot 2.3 11  2.3 10 /
\plot 3.4 10  3.4 9 /
\plot 4.5  9  4.5 8 /
\plot 5.6  8  5.6 7 /
\plot 6.5  9  6.5 8 /

\plot 4.1 13  4.1 14 /

\setshadegrid span <.6mm>
\vshade 1.8 2 3  <,z,,> 4 0 1  <z,,,>  7.6 4 5 /
\vshade 2.7 4 5  <,z,,> 3.8 3 4  <z,,,> 6.5   6 7 /
\vshade 2.5 7 8  <,z,,> 4.7 5 6  <z,,,> 6.5   7 8 /
\vshade 1.2 11 12  <,z,,> 5.6 7 8  <z,,,> 6.5 8 9 /
\vshade 1.2 12 13  <,z,,> 2.3 11 12  <z,,,> 4.1 13 14 /

\multiput{$\bullet$} at 1.6 6   5.6 6  2.5 7  / 
\plot 2.7 5  1.6 6   2.5 7 /
\plot 2.5 7  3.6 6 /

\put{$\bullet$} at 3.2 0 
\plot 3.2 0  4 1 /

\multiput{$\bullet$} at 
  2.3 12  3.2 13  4.1 14  5 15 
  1.2 13  2.1 14  3 15  3.9 16 
  0.1 14  1 15  1.9 16  2.8 17  3.7 18 
  -1 15  -.1 16    .8 17  1.7 18  2.6 19  3.5 20  4.4 21 
  -2.1 16 -1.2 17 -.3 18  0.6 19  1.5 20    2.4 21  3.3 22               
                             1.3 22  2.2 23 /
\plot 3.3 22  -2.1 16  2.3 12  5 15  0.6 19 /

\multiput{$\bullet$} at  2.2 24  2.2  25  2.2 26  2.2 27  2.2 28 /
\plot 2.2 23  2.2 28  /
\plot 3.5 20  1.3 22  2.2 23  4.4 21  -1 15 /
\plot 0.1 14  3.7 18  1.5 20 /
\plot 1.2 13  3.9 16 /
\plot 3.2 13 -1.2 17 /
\plot 4.1 14  -.3 18 /

\put{$a\strut$} at 0 -.6
\put{$b\strut$} at 2 -.6
\put{$c\strut$} at 4 -.6
\put{$d\strut$} at 6 -.6
\put{$e\strut$} at 8 -.6
\put{$f\strut$} at 10 -.6
\put{$g\strut$} at 12 -.6
\put{$u\strut$} at 3.2 -.6
\multiput{$\ssize a$} at 
         -1.8 15.3 
         7.1 10.8 /
\multiput{$\ssize b$} at  0.3 0.8  1.9 5.3 
         -.7  14.3  1.7 21.3
         7.3 5.8 /
\multiput{$\ssize c$} at  1.2 1.7  2.5 4.5  1.8 8.3
         .3 13.3  1.75 20.7  1.9 23.5 
         8.4 3.8  /
\multiput{$\ssize d$} at  2.1 3.7  2.25 7.5  .9 12.5 
          0.9 19.7  1.9 24.5
         9.6 2.7 /
\multiput{$\ssize e$} at  1.9 6.7  .9 11.5 
           0 18.7  1.9 25.5
         10.7 1.7 /
\multiput{$\ssize f$} at  0.6 10.7 
           -.9 17.7   1.9 26.5
          11.75 0.7 /
\multiput{$\ssize g$} at   -1.8 16.7   1.9 27.5 /
\multiput{$\ssize u$} at  1.5 2.5  0.6 9.3 
         1.6 22.7 /
\endpicture}
$$

\vfill\eject
$$
\hbox{\beginpicture
\setcoordinatesystem units <.92cm,.7cm>
\multiput{} at 0 1  13 0 /
\put{\bf Level} at 1 1
\put{ 1} at 1 0
\put{ 2} at 1 -1
\put{ 3} at 1 -2
\put{ 4} at 1 -3
\put{ 5} at 1 -4
\put{ 6} at 1 -5
\put{Conditions}  at 3.5 1
\put{$u = 0$}  at 3.5 0
\put{$u = 1,\ c\le 1$}  at 3.5 -1
\put{$c=2,\ d=1$}  at 3.5 -2
\put{$d=2,\ e=1$}  at 3.5 -3
\put{$d=2,\ e=3$}  at 3.5 -4
\put{$d\ge 3$}  at 3.5 -5
\put{Minimal elements}  at 7 1
\put{Maximal element}  at 10.7 1 
\put{7 simple roots} at 7 0
\put{$\smallmatrix  & & 0 \cr
         1 & 1 & 1 & 1 & 1 & 1 & 1\endsmallmatrix$} at 10.7 0 
\put{$\smallmatrix  & & 1 \cr
         0 & 0 & 0 & 0 & 0 & 0 & 0\endsmallmatrix$} at 7  -1 
\put{$\smallmatrix  & & 1 \cr
         1 & 1 & 1 & 1 & 1 & 1 & 1\endsmallmatrix$} at 10.7  -1 
\put{$\smallmatrix  & & 1 \cr
         0 & 1 & 2 & 1 & 0 & 0 & 0\endsmallmatrix$} at 7  -2 
\put{$\smallmatrix  & & 1 \cr
         1 & 2 & 2 & 1 & 1 & 1 & 1\endsmallmatrix$} at 10.7  -2 
\put{$\smallmatrix  & & 1 \cr
         0 & 1 & 2 & 2 & 1 & 0 & 0\endsmallmatrix$} at 7  -3 
\put{$\smallmatrix  & & 2 \cr
         1 & 2 & 3 & 2 & 1 & 1 & 1\endsmallmatrix$} at 10.7  -3 
\put{$\smallmatrix  & & 1 \cr
         0 & 1 & 2 & 2 & 2 & 1 & 0\endsmallmatrix$} at 7  -4 
\put{$\smallmatrix  & & 2 \cr
         1 & 2 & 3 & 2 & 2 & 2 & 1\endsmallmatrix$} at 10.7  -4 
\put{$\smallmatrix  & & 1 \cr
         1 & 2 & 3 & 3 & 2 & 1 & 0\endsmallmatrix$} at 7  -5 
\put{$\smallmatrix  & & 3 \cr
         2 & 4 & 6 & 5 & 4 & 3 & 12\endsmallmatrix$} at 10.7  -5 
\plot 0.6 0.5  14.5 0.5 /
\plot 2 1.3  2 -5.4 /
\plot 5 1.3  5 -5.4 /
\plot 8.8 1.3  8.8 -5.4 /
\plot 12.5 1.3  12.5 -5.4 /

\put{number}  at 13.5 1
\put{$28$} [l] at 13 0
\put{$3\times5+1$} [l] at 13 -1
\put{$3\times4$} [l] at 13 -2
\put{$5\times3$} [l] at 13 -3
\put{$5\times3$} [l] at 13 -4
\put{$5\times4+14$} [l] at 13 -5

\endpicture}
$$                           

\bigskip
\hrule
\bigskip
$$
\hbox{\beginpicture
\setcoordinatesystem units <1cm,.9cm>
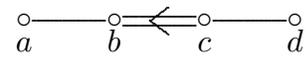
\put{\bf The root poset $\Bbb F_4.$} at 0 10.5 
\put
{\beginpicture
\setcoordinatesystem units <1cm,.5cm>

\multiput{$\circ$} at 0 0  1 0  2 0  3 0 /
\put{$a\strut$} at 0 -.5
\put{$b\strut$} at 1 -.5
\put{$c\strut$} at 2 -.5
\put{$d\strut$} at 3 -.5
\plot 0.1 0  0.9 0 /
\plot 1.1 0.1  1.9 0.1 /
\plot 1.1 -.1  1.9 -.1 /
\plot 2.1 0  2.9 0 /
\plot 1.6 -.3 1.4 0  1.6 .3 /
\endpicture} at 8 10.5

\multiput{$\bullet$} at 
     0 0  2 0  4 0  6 0 
     1 1  3 1  5 1
     2 2  3 2  4 2 
     2 3  4 3
     1 4  3 4  5 4
     2 5  4 5
     3 6  5 6
     4 7  4 8  4 9  4 10  /
\multiput{$\circ$} at 3 3 /
\plot  0 0  2 2  4 0 5 1  /
\plot 1 1  2 0  4 2  6 0 /
\plot 3 1  3 2 /
\plot 2 2  2 3  4 5 /
\plot 4 2  4 3  2 5 /
\plot 3 2  5 4  3 6  1 4  3 2 /
\plot 3 6  4 7  5 6  4 5 /
\plot 4 7  4 10 /

\setdashes <1mm>
\plot 2 2  3 3  4 2 /
\plot 3 3  3 4 /

\setshadegrid span <.6mm>
\vshade 2 2 3  <,z,,> 3 1 2  <z,,,>  4 2 3  /
\put{$a\strut$} at 0 -.3
\put{$b\strut$} at 2 -.3
\put{$c\strut$} at 4 -.3
\put{$d\strut$} at 6 -.3
\multiput{$\ssize a$} at  1.4 3.3  4.6 4.7 /
\multiput{$\ssize b$} at  0.4 0.7  1.8 2.5  3.4 6.7  3.8 7.5 
              4.6 1.7 /
\multiput{$\ssize c$} at  1.4 1.7  2.4 5.7  3.8 8.5 
              5.6 .7 /
\multiput{$\ssize d$} at  1.4 4.7  3.8 9.5 /
\endpicture}
$$
	\bigskip  
$$
\hbox{\beginpicture
\setcoordinatesystem units <.95cm,.7cm>
\multiput{} at 0 1  13 0 /
\put{\bf Level} at 1 1
\put{ 1} at 1 0
\put{ 2} at 1 -1
\put{Conditions}  at 3.5 1
\put{$b\le 1$}  at 3.5 0
\put{$b\ge 2$}  at 3.5 -1
\put{Minimal elements}  at 7 1
\put{Maximal element}  at 10.7 1 
\put{4 simple roots} at 7 0
\put{$\smallmatrix  1 & 1 & 1 & 1 \endsmallmatrix$} at 10.7 0 
\put{$\smallmatrix  0 & 2 & 1 & 0 \endsmallmatrix$} at 7  -1 
\put{$\smallmatrix  2 & 4 & 3 & 2 \endsmallmatrix$} at 10.7  -1 
\plot 0.6 0.5  14.5 0.5 /
\plot 2 1.3  2 -1.4 /
\plot 5 1.3  5 -1.4 /
\plot 8.8 1.3  8.8 -1.4 /
\plot 12.5 1.3  12.5 -1.4 /

\put{number}  at 13.5 1
\put{$10$} [l] at 13 0
\put{$3\times3+5$} [l] at 13 -1

\endpicture}
$$                           

\vfill\eject 

$$
\hbox{\beginpicture
\setcoordinatesystem units <1cm,.8cm>
\put{\bf The root poset $\Bbb G_2.$} at 0 0 

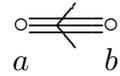
\put{\beginpicture
\put{} at -7 0 
\multiput{$\circ$} at 0 0  1 0 /
\put{$a\strut$} at 0 -.5
\put{$b\strut$} at 1 -.5
\plot 0.1 0  0.9 0 /
\plot 0.1 0.1  0.9 0.1 /
\plot 0.1 -.1  0.9 -.1 /
\plot 0.6 -.3  0.4 0  0.6 .3 /
\endpicture} at 4 0
\put{} at 11 0
\endpicture}
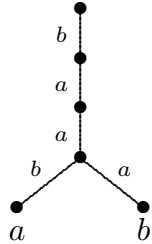
$$
$$
\hbox{\beginpicture
\setcoordinatesystem units <.7cm,.55cm>
\put{$a\strut$} at 0 -.5
\put{$b\strut$} at 2 -.5
\multiput{$\ssize a\strut$} at 1.7 0.8  0.7 1.5  0.7 2.5 /
\multiput{$\ssize b\strut$} at 0.3 0.8  0.7 3.5 /

\multiput{$\bullet$} at 0 0  2 0
  1 1   1 2  1 3  1 4 /
\plot 0 0  1 1  2 0 /
\plot 1 1  1 4 /
\endpicture}
$$

   \bigskip\bigskip 
{\bf References.}
     \medskip 
\item{[Ar]} Armstrong, D.: Generalized noncrossing
 partitions and combinatorics of Coxeter groups.
 Memoirs of the Amer\. Math\. Soc\. 949 (2009).
\item{[At]} Athanasiadis, C\. A\.: On a refinement of the generalized
   Catalan numbers for Weyl groups. Trans\. Amer\. Math\. Soc\. 357 (2004), 179-196.
\item{[B]} Bourbaki, N.: Groupes et Algebres de Lie. 
 IV-VI. Hermann (1968).
\item{[H]} Humphreys, J\. E\.: Reflection groups and
 Coxeter groups. Cambridge University Press (1990).
\item{[S]} Richard Stanley, Enumerative Combinatorics, vol.1, 
    Cambridge Studies in Advanced Mathematics 49, 
    Cambridge University Press, 1995, 
\item{[R1]} Ringel, C\. M\.: 
    The $(n-1)$-antichains in a root poset of width $n$.  arXiv:1306.1593v1
\item{[R2]} Ringel, C\. M\.: 
   The Catalan combinatorics of the hereditary artin algebras.
    In: Recent Developments in Representation Theory, Contemp. Math., 673, 
    Amer. Math. Soc.,  Providence, RI, 2016. 
\item{[R3]} Ringel, C\. M\.: Root posets and hammocks. In preparation. 
	\bigskip

 \bigskip\bigskip

\parindent=0truecm
Claus Michael Ringel \par
Fakult\"at f\"ur Mathematik, Universit\"at Bielefeld\par
POBox 100 131, 33 501 Bielefeld, Germany \par
{\tt e-mail: ringel\@math.uni-bielefeld.de}

\bye